\documentclass[10pt]{article}
\usepackage{pstricks}
\usepackage{pst-node,pst-coil}
\usepackage{color}
\usepackage[normalem]{ulem}
\usepackage{pictexwd}

\def\twocoaln#1#2#3{
\ncarc{-}{#1}{#3}
\ncarc{-}{#2}{#3}
}
\usepackage{latexsym,amssymb, amsthm}
\usepackage[centertags]{amsmath}
\newcommand{\Rand}[1]{\marginpar{#1}}
     \renewcommand{\Rand}[1]{}
\newcommand{\be}[1]{\Rand{\vspace{0,6cm}\tt #1}\begin{equation}\label{#1}}
\newcommand{\bew}[1]{\Rand{\vspace{0,6cm}\tt #1}\begin{equation*}\label{#1}}
\newcommand{\bea}[1]{\Rand{\vspace{0,6cm}\tt #1}\begin{eqnarray}\label{#1}}

\newcommand{\beL}[2]{\Rand{\vspace{0,6cm}\tt #1}\begin{lemma}[#2]\label{#1}}
\newcommand{\beD}[2]{\Rand{\vspace{0,6cm}\tt #1}\begin{definition}[#2]\label{#1}}
\newcommand{\beT}[2]{\Rand{\vspace{0,6cm}\tt #1}\begin{theorem}[#2]\label{#1}}
\newcommand{\beP}[2]{\Rand{\vspace{0,6cm}\tt #1}\begin{proposition}[#2]\label{#1}}
\newcommand{\beC}[2]{\Rand{\vspace{0,6cm}\tt #1}\begin{corollary}[#2]\label{#1}}

\newcommand{\sm}{\smallskip}
\newcommand{\bi}{\bigskip}
\newcommand{\wh}{\widehat}

\newcommand{\bean}{\begin{eqnarray*}}
\newcommand{\ee}{\end{equation}}
\newcommand{\eea}{\end{eqnarray}}
\newcommand{\eean}{\end{eqnarray*}}
\newcommand{\el}{\end{lemma}}
\newfont{\graf}{eufm10}% scaled\magstep2}

\newcommand{\N}{\mathbb{N}}
\newcommand{\R}{\mathbb{R}}

\newcommand{\Z}{\mathbb{Z}}

\newcommand{\suml}{\sum\limits}

\newcommand{\noi}{\noindent}
\newcommand{\ed}{\stackrel{d}{=}}

\newcommand{\ve}{\varepsilon}

\newcommand{\CA}{{\cal A}}

\newcommand{\CI}{{\cal I}}

\newcommand{\CP}{{\cal P}}
\newcommand{\CL}{{\cal L}}

\newcommand{\CM}{{\cal M}}
\newcommand{\CN}{{\cal N}}

\newcommand{\CE}{{\cal E}}

\newtheorem{proposition}{Proposition}[section]
\newtheorem{corollary}{Corollary}[section]
\newtheorem{lemma}{Lemma}[section]
\newtheorem{definition}{Definition}[section]
\newtheorem{remark}{Remark}[section]
\newtheorem{theorem}{Theorem}
\newtheorem{example}{Example}[section]

\newcommand{\D}{\displaystyle}

\newcommand{\la}{\longrightarrow}
\newcommand{\La}{\Longrightarrow}

\newcommand{\rhoto}{{_{\D \longrightarrow \atop \rho \to \infty}}}

\newcommand{\nto}{{_{\D \longrightarrow \atop n \to \infty}}}
\newcommand{\tto}{{_{\D \longrightarrow \atop t \to \infty}}}
\newcommand{\Tto}{{_{\D \Longrightarrow \atop t \to \infty}}}
\newcommand{\Tno}{{_{\D \Longrightarrow \atop n \to \infty}}}
\newcommand{\TuO}{{_{\D \Longrightarrow \atop u \to 0}}}
\newcommand{\Tmo}{{_{\D \Longrightarrow \atop m \to \infty}}}

\newcommand{\Nto}{{_{\D \Longrightarrow \atop N \to \infty}}}

\renewcommand{\ln}{\log}
\begin{document}
\numberwithin{equation}{section}
%\numberwithin{proposition}{section}
%\numberwithin{definition}{section}
%\numberwithin{remark}{section}

\title{ {\Large Coalescent processes arising in
a study of diffusive clustering}}
\author{Andreas Greven$^{(1)}$ \mbox{ }
Vlada Limic$^{(2)}$ \mbox{ } Anita Winter$^{(1)}$}
\date{July 12, 2009 \\
{\small \thanks{Acknowledgment: The research was in part supported by
the DFG Forschergruppe 498, through grant GR-876/13,1-2, by NSERC, and
by the Alfred P.~Sloan Research Fellowship of Vlada Limic and in part
at the Technion by a fellowship of the
Aly Kaufman Foundation of Anita Winter.}}\\
{\small (script/winter/coalclu/090603GLW2-rev10.tex)}}
\maketitle

\begin{abstract}
This paper studies the spatial coalescent on $\Z^2$.
In our setting, the partition elements are located at the sites of
$\Z^2$ and undergo local delayed coalescence and migration. That is,
pairs of partition elements located at the same site coalesce
into one partition element
after exponential waiting times. In addition, the partition
elements perform independent random walks.
The system starts in either locally finite configurations or in
configurations containing countably
many partition elements  per site.
These two situations are relevant if the coalescent is used to
study the scaling limits for genealogies in Moran models
respectively interacting Fisher-Wright diffusions (or
Fleming-Viot processes), which is the key application of
the present work.

Our goal is to determine the longtime behavior with an initial
population of countably many individuals per site restricted to a box
$[-t^{\alpha/2}, t^{\alpha/2}]^2 \cap \Z^2$
and observed at time $t^\beta$ with $1 \geq \beta \geq \alpha\ge 0$.
We study both asymptotics, as $t\to\infty$,
for a fixed value of $\alpha$
as the parameter $\beta\in[\alpha,1]$ varies,
and for a fixed $\beta$,
as the parameter $\alpha\in [0,\beta]$ varies.
This exhibits the genealogical structure of the mono-type clusters
arising in 2-dimensional Moran and Fisher-Wright systems.

A new random object, the so-called {\em coalescent with rebirth},
is constructed via look-down and shown to arise in the limit. For sake of
completeness,
and in view of
future applications we introduce the spatial coalescent with
rebirth and study its longtime asymptotics as well.

The present paper is the basis for forthcoming works \cite{glw2}
and \cite{GSW}, where the
genealogies in interacting Moran models and
Fisher-Wright diffusions on $\Z^2$ are studied, and
where the spatial continuum limit of the Moran
model on $\Z^1$ (Brownian web) is developed, respectively.
There the coalescent with rebirth is needed to describe
the ``complete'' genealogical forests, i.e., the genealogical
structures which include also the ``fossils''.
\end{abstract}
%\vspace{0.5cm}

\noi
\sm

\noi
AMS-Subject classification:  60K35, 60G09, 92D25, 60K05
\vspace{0.5cm}

{\footnoterule
\noi
\hspace*{0.3cm}{\footnotesize $^{(1)}$
        {Mathematisches Institut, Universit\"at Erlangen--N\"urnberg,
        Bismarckstra\ss{}e $1\frac12$, 91054 Erlangen, \\
        \hspace*{0.3cm} Germany}}\\
        \hspace*{0.3cm} {\footnotesize $^{(2)}$
        {CR1, CNRS,
	Universit\'{e} de Provence, Technop\^{o}le de Ch\^{a}teau-Gombert
	UMR 6632, LATP, CMI\\
	39, rue F. Joliot Curie, 13453 Marseille, cedex 13, France}}}

\newpage

\tableofcontents

\newpage

\section{Introduction}
\label{S:Intro}

%\subsection{Background}
%\label{Sub:backg}
The spatial (delayed) coalescent processes on $\Z^d$
and their space-time scale behavior
are the key mathematical tools for the analysis of the
%($d=2$)
asymptotic behavior of a certain class of neutral population  models, namely of
interacting particle models,
known as the interacting Moran models, and their diffusion limit,
the interacting Fisher-Wright diffusions, respectively Fleming-Viot diffusions. These models describe
populations in which individuals have a type and a geographic location
evolving by resampling and migration. (See Shiga \cite{Sh80}
and Durrett \cite{Durrett02}).
The coalescent process allows to construct the genealogies of the
current population in these models explicitly. In particular if one attempts to
understand the scaling behavior of the genealogical trees generated by
the population in critical dimension the spatial coalescent process is
the key tool.

We believe, however, that a number of the results on the spatial
coalescent are of independent interest and have
possible applications outside of the context of Moran and
Fleming-Viot models.
For this reason we present and prove them here separately,
and refer the reader for example to \cite{GreLimWin05,glw2,glw3}
for population model applications.
In this paper the results are formulated for individuals, types and locations only and do not
involve continuum constructions using $\R$-trees etc., which will be
necessary in forthcoming work \cite{glw2,glw3} that builds on the
results presented here and in fact motivates many constructions
in the form given here.

A class of spatial stochastic systems on $\Z^d$
that combine migration between the sites and
a stochastic mechanism acting at each site
(including the voter model,
branching random walks
or interacting diffusions, see, for example,
Liggett \cite{Lig85}, Dawson \cite{Daw93}, Shiga \cite{Sh80}
and Cox and Greven \cite{CoxGre94})
exhibit a dichotomy between
low (typically $d \leq 2$) and high dimensions (typically $d > 2$)
in their longtime behavior.
In high dimensions non-trivial equilibria
exist, while in low dimension such systems approach
laws which are concentrated on the ``traps'' of the
stochastic evolution, i.e.~on the configurations which
the system cannot ever leave with probability 1.

A special r{\^o}le is played by the \emph{critical dimension} $d=2$,
which is
characterized by the fact that the underlying (symmetrized)
migration random walk is recurrent, while its Green's function
$\sum_{k=1}^n{\bf P}\{X_k=0\}$ grows only logarithmical in $n$.
There (as in general for the recurrent setting)
the above processes converge weakly to a law
concentrated on mono-type configurations as time evolves from $0$ to infinity.
Somewhat surprisingly, as first explored for the voter model by
Cox and Griffeath in 1986 \cite{CoxGri86}, the order of magnitude of
the regions where the system looks mono-type
is not asymptotically deterministic  (unlike in the $d=1$ setting
where we get $\sqrt{t}$ as order of magnitude for the size of the mono-type regions).
In fact, the mono-type cluster containing the origin
has an area of the order
$t^{\alpha}$, as $t\to\infty$,  where the {\em random exponent} $\alpha$
takes values in $[0,1]$ and its
distribution can be specified as follows:
take a Fisher-Wright diffusion $(Z_t,t\geq 0)$, and define
\be{vl1}
   T
 :=
   \inf\big\{t\ge 0:\, Z_t\in \{0,1\}\big\},
\ee
then $\alpha\ed e^{-T}$ (see \cite{glw2} for details).
This phenomenon is called the \emph{diffusive clustering}.

Another interesting
question concerns the ``age'' of
a cluster.
%an ``monochrome'' cluster around the origin
%of area $t^\alpha$.
More precisely, in particle systems language, suppose the configuration at some large time $t$ contains
 a monochrome cluster around the origin
of area $t^\alpha$.
Then its age is, informally, the amount of time during which this
cluster has already persisted in the spatial volume of volume $t^\alpha$.
It turns out that this age is of the order $t^\beta$, for some
random  $\beta \in (\alpha,1)$.
To obtain more detailed results on cluster formation one needs to consider the
{\em time-space configuration} of the
process providing the information on which types populated a specific site in
space at a specific time.
This type of analysis for the
time-space configuration as a function of $\alpha$ and $\beta$
has been carried out by Fleischmann
and Greven \cite{FleGre94,FleGre96}
for interacting diffusions with components in $[0,1]$
on the hierarchical group
with a symmetric critically recurrent migration kernel,
by Klenke \cite{Klenke96} in the $[0,\infty)$-valued
(branching)
component case on the hierarchical group,
and by Winter \cite{Win02} on $\Z^2$.

The behavior of the above particle systems and their diffusion limit
is reflected in the behavior of their dual
processes, the spatial coalescent, and
the ``time-space'' dual processes, which we introduce here and which we call the  {\em spatial
coalescent with rebirth}. These dual processes generate the genealogies
of the population of Moran and Fisher-Wright systems.

In the present paper we systematically explore the longtime behavior
of the spatial coalescent with and without
rebirth in the geographic space $\Z^2$.
One of the interesting new features concerning the spatial coalescent with
rebirth is that it enables a description of
the whole genealogical structure (including ``fossils'')
rather than only that of the current population at a reference time
of the corresponding population models.
%In other words, it includes the ``fossils''.
The results we shall prove will replace and extend
the earlier ad hoc constructions via spatial or time-space moment dualities
used in previous work by various authors.
The full potential of this genealogical
viewpoint will become apparent in
future applications. For example, in
\cite{glw2} and \cite{GSW} we shall prove convergence theorems
 for the complete
genealogical structure of the coalescent with rebirth in order to
describe the genealogy in the interacting Moran models and the
interacting Fisher-Wright diffusions including ``fossils''.

We decided to devote a separate paper solely
to coalescent processes since we believe that
the coalescent process with rebirth constructed in
Subsection~\ref{Sub:CWR} is likely to appear in the scaling
limit for a whole class of
 similar mathematical population genealogy (particle
 and diffusion) models. In particular four points are important and different
 from previous work:
\begin{enumerate}
\item the universality of the scaling results
in the sense that the migration mechanism
belongs to a {\em large class} of random walks,
\item initial configuration may contain sites with
 countably {\em infinite} number of individuals,
\item
the concept of the coalescent with rebirth
allows for further applications to the study
 of the genealogies
for the underlying population
models via a {\em weighted $\mathbb{R}$-tree-valued} process
\item an {\em analytical} characterization of the coalescent with rebirth
as well as its construction via a look-down procedure.
\end{enumerate}

At the end of Subsection~\ref{Sub:descript}, after the outline, we comment in
detail on earlier work by Cox and Griffeath (\cite{CoxGri86}) and Bramson, Cox and Griffeath (\cite{BramCoxGri86}) who considered the spatial {\em instantaneous
coalescent} with a {\em simple random
walk migration} mechanism.

\section{Models}
\label{S:Models}

\begin{figure}
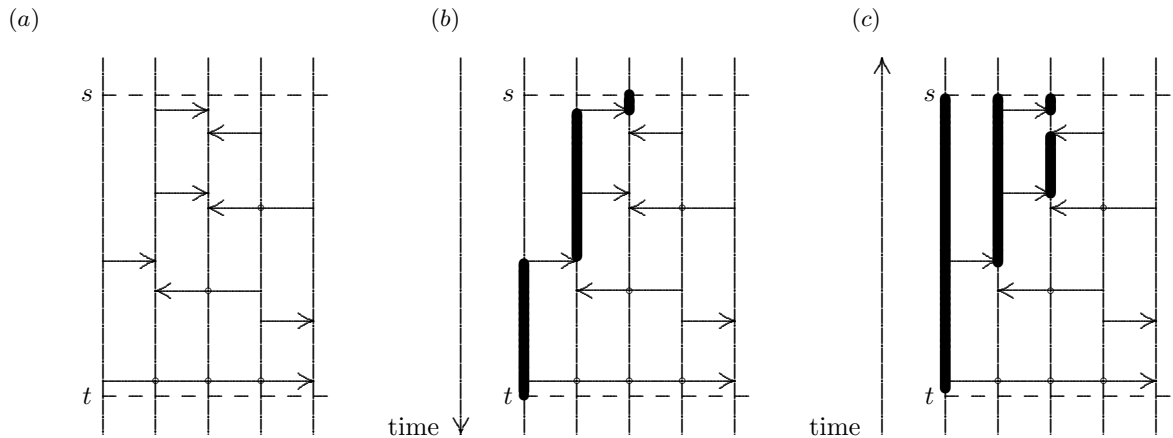

\beginpicture
\setcoordinatesystem units <.7cm,1cm>
\setplotarea x from -0.5 to 12, y from 0 to 5
\put{{\small $(a)$}} at -.5 5.5
\put{{\small $(b)$}} at 7.5 5.5
\put{{\small $(c)$}} at 15.5 5.5
\plot 1 0 1 5 /
\plot 2 0 2 5 /
\plot 3 0 3 5 /
\plot 4 0 4 5 /
\plot 5 0 5 5 /
\arrow <0.2cm> [0.375,1] from 1 0.7 to 5 0.7
\put{\tiny$\circ$} [cC] at 2 0.7
\put{\tiny$\circ$} [cC] at 3 0.7
\put{\tiny$\circ$} [cC] at 4 0.7
\arrow <0.2cm> [0.375,1] from  4 1.5 to 5 1.5
\arrow <0.2cm> [0.375,1] from  4 1.9 to 2 1.9
\put{\tiny$\circ$} [cC] at 3 1.9
\arrow <0.2cm> [0.375,1] from   1 2.3 to 2 2.3
%\arrow <0.2cm> [0.375,1] from 4 2.4 to 5 2.4
\arrow <0.2cm> [0.375,1] from  5 3.0 to 3 3.0
\put{\tiny$\circ$} [cC] at 4 3.0
\arrow <0.2cm> [0.375,1] from 2 3.2 to 3 3.2
\arrow <0.2cm> [0.375,1] from  4 4.0 to 3 4.0
\arrow <0.2cm> [0.375,1] from 2 4.3 to 3 4.3
%%% second figure comes %%%
\arrow <0.2cm> [0.375,1] from 7.8 5 to 7.8 0
\put{time} [cC] at 6.9 0.1
\thinlines
\plot 9 5 9 0 /
\thicklines
  \multiput {$\bullet$} at 9 2.25 *85 0 -.02  /
\plot 10 0 10 5 /
  \multiput {$\bullet$} at 10 2.35 *95 0 .02  /
\plot 11 0 11 5 /
  \multiput {$\bullet$} at 11 4.45 *8 0 -.02  /
\plot 12 0 12 5 /
%  \multiput {$\bullet$} at 12 4.45 *22 0 -.02  /
\plot 13 0 13 5 /
\arrow <0.2cm> [0.375,1] from 9 0.7 to 13 0.7
\put{\tiny$\circ$} [cC] at 10 0.7
\put{\tiny$\circ$} [cC] at 11 0.7
\put{\tiny$\circ$} [cC] at 12 0.7
\arrow <0.2cm> [0.375,1] from 12 1.5 to 13 1.5
\arrow <0.2cm> [0.375,1] from  12 1.9 to 10 1.9
\put{\tiny$\circ$} [cC] at 11 1.9
\arrow <0.2cm> [0.375,1] from 9 2.3 to 10 2.3
%\arrow <0.2cm> [0.375,1] from 12 2.4 to 13 2.4
\arrow <0.2cm> [0.375,1] from  13 3.0 to 11 3.0
\put{\tiny$\circ$} [cC] at 12 3.0
\arrow <0.2cm> [0.375,1] from 10 3.2 to 11 3.2
\arrow <0.2cm> [0.375,1] from  12 4.0 to 11 4.0
\arrow <0.2cm> [0.375,1] from 10 4.3 to 11 4.3
\put{$\bullet$} [cC] at 11 4.5
\put{$\bullet$} [cC] at 9 0.5
%%% third picture %%%
\arrow <0.2cm> [0.375,1] from 15.8 0 to 15.8 5
\put{time} [cC] at 14.9 0.1
\plot 17 0 17 5 /
\thicklines
  \multiput {$\bullet$} at 17 4.45 *193 0 -.02  /
\plot 18 0 18 5 /
  \multiput {$\bullet$} at 18 4.45 *109 0 -.02  /
\plot 19 0 19 5 /
  \multiput {$\bullet$} at 19 3.95 *38 0 -.02  /
  \multiput {$\bullet$} at 19 4.45 *8 0 -.02  /
\plot 20 0 20 5 /
%  \multiput {$\bullet$} at 20 4.45 *22 0 -.02  /
\plot 21 0 21 5 /
\arrow <0.2cm> [0.375,1] from 17 0.7 to 21 0.7
\put{\tiny$\circ$} [cC] at 18 0.7
\put{\tiny$\circ$} [cC] at 19 0.7
\put{\tiny$\circ$} [cC] at 20 0.7
\arrow <0.2cm> [0.375,1] from 20 1.5 to 21 1.5
\arrow <0.2cm> [0.375,1] from  20 1.9 to 18 1.9
\put{\tiny$\circ$} [cC] at 19 1.9
\arrow <0.2cm> [0.375,1] from 17 2.3 to 18 2.3
%\arrow <0.2cm> [0.375,1] from 20 2.4 to 21 2.4
\arrow <0.2cm> [0.375,1] from  21 3.0 to 19 3.0
\put{\tiny$\circ$} [cC] at 20 3.0
\arrow <0.2cm> [0.375,1] from 18 3.2 to 19 3.2
\arrow <0.2cm> [0.375,1] from  20 4.0 to 19 4.0
\arrow <0.2cm> [0.375,1] from 18 4.3 to 19 4.3
\setdashes
\plot 1 0.5 5.4 0.5 /
\plot 1 4.5 5.4 4.5 /
\put{$t$} [cC] at 0.7 0.5
\put{$s$} [cC] at 0.7 4.5
%\put{$t_1$} [cC] at 0.5 4.3
%\put{$t_2$} [cC] at 0.5 2.3
%\put{$u$} [cC] at 0.7  3.2
%
\plot 9 0.5 13.4 0.5 /
\plot 9 4.5 13.4 4.5 /
\put{$t$} [cC] at 8.7 0.5
\put{$s$} [cC] at 8.7 4.5
\plot 17 0.5 21.4 0.5 /
\plot 17 4.5 21.4 4.5 /
\put{$t$} [cC] at 16.7 0.5
\put{$s$} [cC] at 16.7 4.5
\endpicture
\caption{\label{fig:MM1}
(a) a realization of the resampling times, (b) the ancestral line of the individual $3$ which lives at time $s$
is indicated in bold, and (c)
the set of all descendants up to time $s$ of the individual~$1$
present at time $t$ is indicated in bold.
}
\end{figure}

The coalescent processes considered in the present paper are describing the genealogies of neutral population models involving resampling between any  two individuals
where two individuals are replaced by descendants of one of them
(the one at the end of the arrow in Figure~\ref{fig:MM1}).
In the time-reversed evolution these time points correspond to the times at which the ancestral lines of the
two individuals  {\em coalesce} to a common ancestral line (compare with Figure~\ref{fig:MM1}).

We shall define in this section the spatial coalescent, the spatial coalescent
with rebirth and finally the so-called look-down process, which allows for a
graphical representation of our coalescent processes which give straightforward
explicit constructions for a version of these processes.

\subsection{The spatial coalescent on $\Z^2$}
\label{Sub:reform}
Processes describing the dynamics of finitely
many moving and coalescing
particles appeared already in the 1980'ies
(see, for example, \cite{BraGri80,BramCoxGri86,CoxGri86} and
compare Liggett \cite{Lig85} for more detailed references).
For coalescents  representing genealogies of diffusion processes,
it is essential
to allow configurations
with {\em countably many particles per site}
on a countable geographic space. Moreover, while the above papers were only recording occupation numbers at various sites, we will provide a set-up which also exhibits the
partition structure.

The {\em spatial coalescent}
that we analyze in the current paper
was introduced on a class of Abelian groups in \cite{GreLimWin05}.
For the benefit of the reader we briefly recall in three steps the relevant notation,
appropriate topologies and
its construction.  We restrict the setting to $\Z^2$. \bi

\noindent{\em Step~1 (Migration)}\quad
Let $a(\boldsymbol{\cdot},\boldsymbol{\cdot})$ be an irreducible random
walk kernel which  has finite exponential moments, i.e.,
\be{vl2}
   a(x,y)=a(0,y-x),
\ee
for all $x,y\in\Z^2$, and
\be{vl3}
   \sum_{(z_1,z_2)\in\Z^2}e^{\lambda_1 z_1+\lambda_2 z_2}a(0,(z_1,z_2))<\infty,
\ee
for all $\lambda_1,\lambda_2\in\R$. We consider the continuous time random walk with
jump rate 1 and transition probability $a(\cdot, \cdot)$.\sm

We next present the standard way to
construct particle systems
that possibly start in configurations with
countably many particles at some or all sites and which involve migration as mechansim.
These particle systems are constructed as extensions of particle systems which
start in specific {\em locally finite} states and  the dynamics are such that they
guarantee local finiteness of the
particle process at all times $t>0$ (compare also with Remark~\ref{Rem:03}).
To construct such locally finite systems we follow an approach due to
Liggett and Spitzer \cite{LigSpi81}.

Fix a finite measure $\alpha$ on $\Z^2$ with
$\alpha(\{x\})>0$, for all $x\in\Z^2$, such that for a
constant $\Gamma$
\be{alpha}
   \sum_{y\in\Z^2}a(x,y)\alpha(\{y\})
 \leq
   \Gamma\cdot\alpha(\{x\}),
\ee
for all $x\in\Z^2$. Denote by ${\mathcal N}(\Z^2)$ the set of all
locally finite $\N_0$-valued measures on $\Z^2$.
Then
\be{LigSpi}
{\mathcal E}_\alpha
\equiv
{\mathcal E}
 :=
   \big\{\eta\in{\mathcal N}(\Z^2):\,\sum_{x\in \Z^2}
   \eta\{x\}\alpha(\{x\})<\infty\big\}
\ee
is the {\em Liggett-Spitzer} space
(corresponding to $\alpha$).

\begin{remark}[${\mathcal E}$ is a state space]\label{Rem:02}\rm
Let $\{(X_t^{i})_{t\geq 0}:\,{i}\in\CI\}$ be a countable collection of
independent random walks, and put for all $t\ge 0$,
$\eta_t:=\sum_{{i}\in\CI}\delta_{X_t^{i}}\in{\mathcal N}(\Z^2)$.
If
\be{eta0}
   \eta_0\in{\mathcal E},\qquad{\bf P}\mbox{-}\mbox{a.s.},
\ee
%is an element of ${\mathcal E}$
then an easy calculation shows that the process
$(e^{-\Gamma t}\sum_{{i}\in\CI}\alpha(\{X_t^{i}\}))_{t\ge 0}$
is a super-martingale.

In particular, under (\ref{eta0}), for all $t\ge 0$,
\be{vl7}
   {\bf P}\big\{\eta_t\in{\mathcal
     E}\big\}=1,
\ee
that is, $\CE$ is a state space for
$\{\eta_t, t\geq 0\}$. Note furthermore that
(\ref{vl7}) implies $\eta_t\in{\mathcal N}(\Z^2)$, for all $t\ge
0$, almost surely.
As topology on $\CN(\Z^2)$ choose the vague
topology, then $\CE_\alpha$ is a (not closed) subset of the polish space
$\CM(\Z^2)$, the locally finite measures on $\Z^2$.
\hfill$\qed$
\end{remark}\sm

To build in countably many individuals per site we shall make
use of the coalescence mechanism introduced next.
\bi

\noindent{\em Step~2 (Coalescence)}\quad
Recall that a {\em partition} of a set $\CI$ is a collection
$\{\pi_\lambda\}$ of pairwise disjoint subsets
of $\CI$
such that $\CI=\cup_\lambda \pi_\lambda$.  We refer to the
elements of a partition as {\em partition elements}.
Let us denote by
\be{e:PcI}
   \Pi^\CI
 :=
   \mbox{collection of all partitions of $\CI$.}
\ee

For all $\CI'\subseteq \CI$,
write $\rho_{\CI'}$ for the restriction
map from $\Pi^{\CI}$ to $\Pi^{\CI'}$ and hence
for any ${\mathcal P}\in\Pi^\CI$,
the {\em induced partition}
\be{e:ege}
   \rho_{\CI'} \diamond {\mathcal P}
 :=
   \big\{\rho_{\CI'}(\pi);\,\pi\in{\mathcal P}\big\}.
\ee

We say that a sequence $({\mathcal P}_n)_{n\in\N}$ {\em converges in
$\Pi^{\CI}$} if for all finite subsets $\CI'\subseteq \CI$,
the sequence $(\rho_{\CI'}\diamond {\mathcal P}_n)_{n\in\N}$ converges in
$\Pi^{\CI'}$ equipped with the discrete topology.
In particular,
a function
$f:\Pi^{\CI} \to \R$ that depends on $\Pi^\CI$
only through $\Pi^{\CI'}$,
for some finite subset $\CI'\subseteq \CI$, is continuous. Note that $\Pi^\CI$ equipped with this topology is a Polish space.

\begin{definition}[The ${\mathcal I}$-Kingman coalescent]
The ${\mathcal I}$-{\em Kingman coalescent},
or short the {\em Kingman-coalescent}, \label{Def:01}
\be{agre20}
   K:=(K_t)_{t\ge 0},
\ee
is the unique strong Markov process such that for all finite ${\mathcal I}'\subseteq{\mathcal I}$,
the restricted process
\be{restrpro}
   K^{{\mathcal I}'}
 :=
   \rho_{{\mathcal I}'}\diamond K
\ee
is a $\Pi^{{\mathcal I}'}$-valued Markov chain
which starts in some ${\mathcal P}\in\Pi^{{\mathcal I}'}$,
and given $K^{{\mathcal I}'}_t$, each pair of partition elements is merging
to form a single partition element after an exponential waiting time
with rate $\gamma_{\mathrm{King}}>0$.
\end{definition}\bi

\noindent{\em Step~3 (Migration and coalescence combined)}\quad
We next combine migration and coalescence. For that purpose, fix  a {\em site space} ${M}$ which in the present paper is $\Z^2$ unless stated otherwise.
Then
from any ${\mathcal P}\in \Pi^{{\mathcal I}}$
one can form a {\em marked} partition
\be{vl9}
   {\mathcal P}^M
 :=
   \big\{(\pi,L(\pi));\,\pi\in{\mathcal P}\big\},
\ee
by assigning to each partition element $\pi\in{\mathcal P}$,
its {\em location} $L(\pi)\in {M}$.
Put
\be{vl10}
   \Pi^{{\mathcal I},{M}}
 :=\mbox{set of all marked partitions}.
\ee

Note that $\Pi^{{\mathcal I},{M}}$ is a Polish space if we introduce
the topology as follows.
For all ${\mathcal I}'\subseteq {\mathcal I}$ and ${\mathcal P}\in\Pi^{{\mathcal I},{M}}$, we extend the restriction operator as
\be{vl10a}
   \rho_{{\mathcal I}'}\diamond{\mathcal P}
 :=
   \big\{(\rho_{\mathcal I}'(\pi),L(\pi));\,\pi\in{\mathcal P}\big\},
\ee
and say that a sequence $({\mathcal P}_n)_{n\in\N}$ converges in
$\Pi^{{\mathcal I},{M}}$ if and only if for all finite subsets
${\mathcal I}'\subseteq {\mathcal I}$,
the sequence $(\rho_{{\mathcal I}',{M}}\diamond {\mathcal P}_n)_{n\in\N}$
converges in $\Pi^{{\mathcal I}',{M}}$,
equipped with the discrete topology.\sm

We are now ready to define the spatial ${\mathcal I}$-coalescent.
\begin{definition}[The spatial ${\mathcal I}$-coalescent]
\label{Def:00}
The {\em spatial} ${\mathcal I}$-{\em coalescent} on $\Z^2$,
\be{0.4a}
   (C,L)
 :=
   (C_t,L_t)_{t \ge 0}
 =\Big(\big\{(\pi,L_t(\pi));\,\pi\in C_t\big\}\Big)_{t\ge 0},
\end{equation}
is a strong $\Pi^{{\mathcal I},\Z^2}$-valued
Markov process with c\`adl\`ag paths such that for all subsets ${\mathcal
  I}'\subseteq{\mathcal I}$ with
\be{Lip1}
   \sum\nolimits_{\pi\in C_0,\rho_{\CI'}(\pi)\not=\emptyset}
   %1_{\{\pi \cap \CI^\prime \neq \emptyset \}}
   \delta_{L_0(\pi)}
 \in
   {\mathcal E},
\end{equation}
the restricted process
is a $\Pi^{{\mathcal I}',\Z^2}$-valued strong Markov
    particle system which undergoes the following
    two {\rm independent} mechanisms:
\begin{itemize}
\item {\rm Migration }\;
The marks of the partition elements perform independent
continuous time random walks with rate 1 and transition kernel $a(\cdot,\cdot)$.
\item {\rm Coalescence }\;
Each pair of partition elements whose locations are equal merges into one
partition element independently after exponential waiting times with rate
$\gamma$.
\end{itemize}
\end{definition}\sm

\begin{remark}[Spatial coalescent is well-defined] \rm
\begin{itemize} \label{Rem:03}
\item[{}]
\item[(i)]
Note that for the marked ${\mathcal I}$-coalescent process above there is a natural coupling with
the migration random walks
$\{(X_t^{i})_{t\geq 0}:\,{i}\in{\mathcal I}\}$ such that
\be{vl14}
   \sum\nolimits_{\pi \in C_t}\delta_{L_t(\pi)}(B)
 \leq
   \sum\nolimits_{{i}\in{\mathcal I}}\delta_{X_t^{i}}(B),\qquad\mbox{a.s.},
\ee
for all $B\subseteq\Z^2$. Therefore, by Remark~\ref{Rem:02},
the spatial ${\mathcal I}$-coalescent is locally finite and in particular well-defined if
${\mathcal I}$ is already such
that (\ref{Lip1}) holds.
\item[(ii)] By Proposition 3.4 in \cite{GreLimWin05},
the spatial
${\mathcal I}$-coalescent is well-defined for all initial mark
configurations. Specifically, it is even well-defined if started in a
configuration which contains countable infinitely many partition
elements at each site in $\Z^2$.
In all cases, we have that $\sum_{\pi
  \in C_t} \delta_{L_t(\pi)}\in{\mathcal E}$, almost surely, for all $t>0$.
\hfill $\qed$
\end{itemize}
\end{remark}\sm

\begin{remark}[Consistency Property] \label{Rem:consistency}\rm
In all of our constructions of concrete realizations of coalescents below we use
the following important {\em consistency} property:
if $(C,L)$ is the $\CI$-coalescent and $\CI'\subseteq\CI$ then
$\rho_{{\mathcal I}'}\diamond(C,L)$ is the ${\mathcal I}'$-coalescent started in $\rho_{{\mathcal I}'}\diamond(C_0,L_0)$.
\hfill$\qed$\end{remark}\sm

\begin{remark}[Instantaneous coalescent; $\gamma=\infty$] \label{Rem:07}\rm
Note that if the finite parameter $\gamma$ is replaced by
$\infty$, the spatial coalescent  changes into
the system of {\em instantaneously} coalescing random walks
which can be obtained from particle systems as they are considered in
\cite{BramCoxGri86,CoxGri86}.
This observation will become important in Section~\ref{S:ASP}.
\hfill$\qed$
\end{remark}\sm

\subsection{The coalescent with rebirth}
\label{Sub:CWR}
In neutral population models
coalescent processes arise in the study of genealogical relationships
between individuals {\em currently} alive by looking in reversed time.
Each coalescent event corresponds to
a splitting of an ancestral line and a simultaneous death
of another ancestral line in a forward
population model.
However, if one considers genealogies which include also the
the individuals alive at earlier times (commonly referred to as ``fossils''),
then a richer object than the spatial coalescent is needed.
We call this new object the {\em coalescent with rebirth}.
The coalescent with rebirth accounts in the forward model
for the descendant lines which died
before the current time.
More precisely, whenever an individual dies and gets replaced by a descendent of another
individual in the forward model, in the time-reversed model
the coalescent dynamics with rebirth  generates a new individual at
the corresponding time.

\begin{figure}
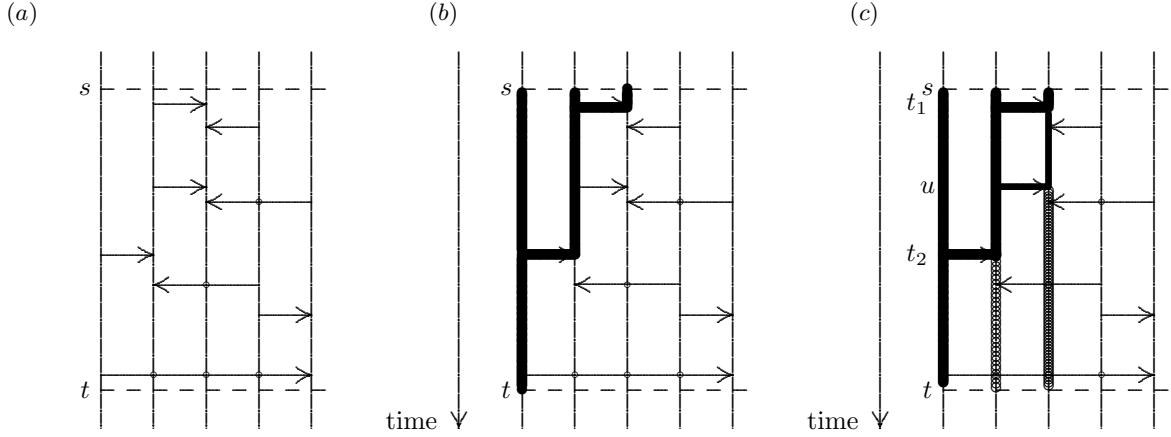

\beginpicture
\setcoordinatesystem units <.7cm,1cm>
\setplotarea x from -0.5 to 12, y from 0 to 5
\put{{\small $(a)$}} at -.5 5.5
\put{{\small $(b)$}} at 7.5 5.5
\put{{\small $(c)$}} at 15.5 5.5
\plot 1 0 1 5 /
\plot 2 0 2 5 /
\plot 3 0 3 5 /
\plot 4 0 4 5 /
\plot 5 0 5 5 /
\arrow <0.2cm> [0.375,1] from 1 0.7 to 5 0.7
\put{\tiny$\circ$} [cC] at 2 0.7
\put{\tiny$\circ$} [cC] at 3 0.7
\put{\tiny$\circ$} [cC] at 4 0.7
\arrow <0.2cm> [0.375,1] from  4 1.5 to 5 1.5
\arrow <0.2cm> [0.375,1] from  4 1.9 to 2 1.9
\put{\tiny$\circ$} [cC] at 3 1.9
\arrow <0.2cm> [0.375,1] from   1 2.3 to 2 2.3
%\arrow <0.2cm> [0.375,1] from 4 2.4 to 5 2.4
\arrow <0.2cm> [0.375,1] from  5 3.0 to 3 3.0
\put{\tiny$\circ$} [cC] at 4 3.0
\arrow <0.2cm> [0.375,1] from 2 3.2 to 3 3.2
\arrow <0.2cm> [0.375,1] from  4 4.0 to 3 4.0
\arrow <0.2cm> [0.375,1] from 2 4.3 to 3 4.3
%%% second figure comes %%%
\arrow <0.2cm> [0.375,1] from 7.8 5 to 7.8 0
\put{time} [cC] at 6.9 0.1
\thinlines
\plot 9 5 9 0 /
\thicklines
  \multiput {$\bullet$} at 9 2.25 *85 0 -.02  /
\plot 10 0 10 5 /
%  \multiput {$\bullet$} at 10 2.35 *40 0 .02  /
  \multiput {$\bullet$} at 9 2.35 *95 0 .02  /
\plot 11 0 11 5 /
  \multiput {$\bullet$} at 11 4.45 *8 0 -.02  /
  \multiput {$\bullet$} at 9 4.45 *8 0 -.02  /
  \multiput {$\bullet$} at 10 4.45 *108 0 -.02  /
\plot 12 0 12 5 /
%  \multiput {$\bullet$} at 12 4.45 *22 0 -.02  /
\plot 13 0 13 5 /
\arrow <0.2cm> [0.375,1] from 9 0.7 to 13 0.7
\put{\tiny$\circ$} [cC] at 10 0.7
\put{\tiny$\circ$} [cC] at 11 0.7
\put{\tiny$\circ$} [cC] at 12 0.7
\arrow <0.2cm> [0.375,1] from 12 1.5 to 13 1.5
\arrow <0.2cm> [0.375,1] from  12 1.9 to 10 1.9
\put{\tiny$\circ$} [cC] at 11 1.9
\arrow <0.2cm> [0.375,1] from 9 2.3 to 10 2.3
%\arrow <0.2cm> [0.375,1] from 12 2.4 to 13 2.4
\arrow <0.2cm> [0.375,1] from  13 3.0 to 11 3.0
\put{\tiny$\circ$} [cC] at 12 3.0
\arrow <0.2cm> [0.375,1] from 10 3.2 to 11 3.2
\arrow <0.2cm> [0.375,1] from  12 4.0 to 11 4.0
\arrow <0.2cm> [0.375,1] from 10 4.3 to 11 4.3
\put{$\bullet$} [cC] at 11 4.5
\put{$\bullet$} [cC] at 9 0.5
%%% third picture %%%
\arrow <0.2cm> [0.375,1] from 15.8 5 to 15.8 0
\put{time} [cC] at 14.9 0.1
\plot 17 0 17 5 /
\thicklines
  \multiput {$\bullet$} at 17 2.3 *50 .02 0 /
  \multiput {$\bullet$} at 18 4.25 *50 .02 0 /
  \multiput {$\bullet$} at 9 2.3 *50 .02 0 /
  \multiput {$\bullet$} at 10 4.25 *50 .02 0 /
  \multiput {$\bullet$} at 17 4.45 *193 0 -.02  /
\plot 18 0 18 5 /
  \multiput {$\bullet$} at 18 4.45 *109 0 -.02  /
\plot 19 0 19 5 /
  \multiput {{\tiny $\bullet$}} at 19 4.18 *48 0 -.02  /
  \multiput {{\tiny $\bullet$}} at 19 3.2 *50 -.02 0 /
  \multiput {$\circ$} at 19 3.15 *65 0 -.04  /
  \multiput {$\circ$} at 18 2.2 *28 0 -.06  /
  \multiput {$\bullet$} at 19 4.45 *8 0 -.02  /
\plot 20 0 20 5 /
%  \multiput {$\bullet$} at 20 4.45 *22 0 -.02  /
\plot 21 0 21 5 /
\arrow <0.2cm> [0.375,1] from 17 0.7 to 21 0.7
\put{\tiny$\circ$} [cC] at 18 0.7
\put{\tiny$\circ$} [cC] at 19 0.7
\put{\tiny$\circ$} [cC] at 20 0.7
\arrow <0.2cm> [0.375,1] from 20 1.5 to 21 1.5
\arrow <0.2cm> [0.375,1] from  20 1.9 to 18 1.9
\put{\tiny$\circ$} [cC] at 19 1.9
\arrow <0.2cm> [0.375,1] from 17 2.3 to 18 2.3
%\arrow <0.2cm> [0.375,1] from 20 2.4 to 21 2.4
\arrow <0.2cm> [0.375,1] from  21 3.0 to 19 3.0
\put{\tiny$\circ$} [cC] at 20 3.0
\arrow <0.2cm> [0.375,1] from 18 3.2 to 19 3.2
\arrow <0.2cm> [0.375,1] from  20 4.0 to 19 4.0
\arrow <0.2cm> [0.375,1] from 18 4.3 to 19 4.3
\setdashes
\plot 1 0.5 5.4 0.5 /
\plot 1 4.5 5.4 4.5 /
\put{$t$} [cC] at 0.7 0.5
\put{$s$} [cC] at 0.7 4.5
\put{$t_1$} [cC] at 16.5 4.3
\put{$t_2$} [cC] at 16.5 2.3
\put{$u$} [cC] at 16.7  3.2
\plot 9 0.5 13.4 0.5 /
\plot 9 4.5 13.4 4.5 /
\put{$t$} [cC] at 8.7 0.5
\put{$s$} [cC] at 8.7 4.5
\plot 17 0.5 21.4 0.5 /
\plot 17 4.5 21.4 4.5 /
\put{$t$} [cC] at 16.7 0.5
\put{$s$} [cC] at 16.7 4.5
\endpicture
\caption{\label{fig:MM2}
(a) a realization of the resampling times, (b) the
genealogy of the first three particles alive at time $s$ is drawn in bold, and (c) the enriched genealogy
after the fossils are included (here the different patterns correspond to different rebirth events).
}
\end{figure}

\begin{example}\rm
\label{Exp:01}
Assume ${\mathcal I}:=\{1,2,3\}$ and
consider the initial configuration
 $\{\{1\},\{2\},\{3\}\}$ at time $s$,
and the transitions (without rebirth)
at times $t_1$ and $t_2$, respectively,
\be{agre40}
   \big\{\{1\},\{2\},\{3\}\big\}
 \stackrel{\mbox{at time $t_1$}}{\longrightarrow}
   \big\{\{1\},\{2,3\}\big\}
 \stackrel{\mbox{at time $t_2$}}{\longrightarrow}
   \big\{\{1,2,3\}\big\}
\ee
(compare, for example, Figure~\ref{fig:MM2}).
Then the corresponding coalescent with rebirth would
start at time $s$ from  $\{\{(1,s)\},\{(2,s)\},\{(3,s)\}\}$,
 and at time $t_1$ the state would change to
$\{\{(1,s)\},\{(2,s) ,(3,s)\},\{(3,t_1)\}\}$.
In particular $(3,t_1)$ corresponds to the ``reborn'' individual $(3,s)$.
Next after time $t_2$ the new individual $(2,t_2)$ is born, etc.
All new born partition elements also undergo resampling and migration. Assume, for example, that in addition
to the above mentioned resampling events, there would be one at time $u\in(t_1,t_2)$ between
the ancestral lines of $2$ and $3$ then we would observe, for example,  the following transitions
\be{agre40a}
\begin{aligned}
   &\big\{\{(1,s)\},\{(2,s)\},\{(3,s)\}\big\}
  \\
  \stackrel{\mbox{at time $t_1$}}{\longrightarrow}
   &\big\{\{(1,s)\},\{(2,s),(3,s)\},\{(3,t_1)\}\big\}
  \\
  \stackrel{\mbox{at time $u$}}{\longrightarrow}
   &\big\{\{(1,s)\},\{(2,s),(3,s),(3,t_1)\},\{(3,u)\}\big\}
  \\
 \stackrel{\mbox{at time $t_2$}}{\longrightarrow}
   &\big\{\{(1,s),(2,s),(3,s),(3,t_1)\},\{(2,t_2)\},\{(3,u)\}\big\}.
\end{aligned}
\ee
Notice that the transition at time $u$ is not observable in the original coalescent.
%Moreover, since we replace the individual which is ``lost''.
\hfill$\qed$
\end{example}\sm

The goal of this subsection is to introduce the coalescent with rebirth first in the non-spatial setting
%in Subsubsection~\ref{SubSub:CR}
and then in the spatial setting.
%in Subsubsection~\ref{SubSub:spaco}.}
\bi

%\subsubsection{The Kingman-type coalescent with rebirth}
%\label{SubSub:CR}
\noindent {\em Step~1 (Coalescence with rebirth). }
As before let ${\mathcal I}$ be a countable set.

To define the state space of the coalescent with rebirth consider first a subset
%(note that here $\CI$ is only a superscript)
\be{agre3}
   \mathbb{S}(\CI)
 \subseteq
   \CI\times\R^+
\ee
such that for each $i\in{\mathcal I}$,
\be{S}
   \#\big\{t\in\R:\,(i,t)\in\mathbb{S}^{{\mathcal I}}\big\}
 <
   \infty.
\ee
We refer to the elements $(i,t)\in \mathbb{S}(\CI)$
as individuals, and call $t$ the {\em birth time} of the {\em individual}
$(i,t)$. Let $\mathrm{pr}_{\mathrm{index}}$ and $\mathrm{pr}_{\mathrm{time}}$ be the projection maps of individuals
to their {\em indices} and {\em birth times}, i.e.,
\be{proj}
   (i,t)
 =
   \big(\mathrm{pr}_{\mathrm{index}}(i,t),
   \mathrm{pr}_{\mathrm{time}}(i,t)\big),
\ee
for all $(i,t)\in{\mathcal I}\times\R$.

Recall from (\ref{e:PcI}) the collection $\Pi^S$ of all partitions of a set $S$.
Call $\CP$ a {\em sub-partition} of ${\mathbb{S}(\CI)}$, if
it is a partition of a subset of ${\mathbb{S}(\CI)}$, or equivalently;
a collection
$\{\pi_\lambda\}$ of pairwise disjoint subsets of
$\mathbb{S}(\CI)$.
With a slight abuse of notation
denote by
\be{subpart}
   \Pi^{\mathbb{S}(\CI)}
 :=
   \mbox{set of all sub-partitions of $\mathbb{S}(\CI)$.}
\ee

Notice that the coalescent was defined in a symmetric manner.
To define
the coalescence dynamics with rebirth we need to break the involved symmetry and declare which of the patches is getting ``lost'' and simultaneously reborn. For that recall that since ${\mathcal I}$ is countable we can fix
an order relation $\preceq$ such that
for all $i_0\in{\mathcal I}$, $\#\{i\preceq i_0\}<\infty$.
This extends to the {\em lexicographic order} relation on ${\mathcal I}\times\R$, that is,
for $(i,s),(j,t)\in{\mathcal I}\times\R$,
\be{vl19}
   (i,s)\preceq (j,t)
 \qquad
   \mbox{if and only if }
 \quad
   i\preceq j \mbox{ or } i=j \mbox{ and } s<t.
\ee

Given ${\mathcal P}\in\Pi^{\mathbb{S}(\CI)},$
define the label $\kappa(\pi)$ of a partition element $\pi\in{\mathcal P}$ as its smallest element with respect to $\preceq$, i.e.,
\be{label}
   \kappa(\pi)
 :=
   \min\big\{v;\, v\in \pi\big\}.
\ee
As illustrated in Example~\ref{Exp:01} the coalescent with rebirth dynamics relies on the rule that if two partition elements coalesce it is always the one with the bigger label that gets ``lost and reborn''. We therefore need our process to take values in the following
subset of $\Pi^{\mathbb{S}(\CI)}$:
\be{PiBBS}
   \widehat{\Pi}^{\mathbb{S}(\CI)}
:=
   \big\{{\mathcal P}\in\Pi^{\mathbb{S}(\CI)}:\,\exists\mbox{ bijection
   }\iota:{\mathcal I}\to{\mathcal P}\mbox{ s.t. }\forall i\in{\mathcal I}\,\exists t\in\R\mbox{ with }\kappa(\iota(i))=(i,t)\big\}.
\ee

We equip
$\R \times \wh \Pi^{\mathbb{S}}$ with a topology
that takes both the partition structure and the birth times into account.
Notice that in contrast to the original coalescent where the set of individuals is fixed, in the coalescent with rebirth the set of individuals increases as time increases.
Let therefore
for all ${\mathcal P}\in {\wh \Pi^{\mathbb{S}(\CI)}}$,
\be{e:basic}
   S({\mathcal P})
 :=
   \bigcup\nolimits_{\pi\in{\mathcal P}}\pi
\ee
be the {\em basic set} of ${\mathcal P}$.

Recall from (\ref{e:ege}) the restriction map, and abbreviate for
$\CI'\subseteq \CI$ and a subpartition $\CP\in\hat{\Pi}^{\mathbb{S}^{{\mathcal I}}}$,
\begin{equation}
\label{Erhorestr}
   \rho_{\CI'}\diamond \CP
 :=
   \rho_{(\CI'\times \R)\cap S(\CP)}\diamond\CP.
\end{equation}

\begin{example}\rm Take for example ${\mathcal I}:=\{1,2,3\}$ and ${\mathcal P}:=\big\{\{(1,0),(2,0),(3,0),(3,t_1)\},\{(2,t_2)\},\{(3,s)\}\big\}$. Then
\be{Erhores}
   \rho_{\{1,2\}}\diamond{\mathcal P}
 :=
   \big\{\{(1,0),(2,0)\},\{(2,t_2)\}\big\}.
\ee
\hfill$\qed$
\end{example}\sm

We now introduce a topology on the state space $\mathbb{S}(\CI)$
that accounts for the differences in both the indices and the birth
times. Loosely speaking, we say that a sequence
\be{grev1}
   \big({\mathcal P}_n\big)_{n\in\N} \mbox{ converges in }
   \hat{\Pi}^{\mathbb{S}(\CI)},
\ee
if and only if
for each finite subset $\CI'\subseteq\CI$, the projections to the index component of the restricted partitions $\rho_{\CI'}\diamond {\mathcal P}_n$ converge in the discrete topology and the corresponding birth times converge with respect to the Euclidian distance. More precisely, we consider the topology generated by a metric satisfying the properties (\ref{f:003}) through (\ref{f:005}).

We will need some further notation.
For a finite subset $\mathbb{S}'\subset\mathbb{S}(\CI)$,
denote the restriction map from
$\hat{\Pi}^{\mathbb{S}(\CI)}$ to $\hat{\Pi}^{\mathbb{S}'}$ by $\rho_{\mathbb{S}'}$.

Since the coalescent with rebirth is keeping track of the birth time of an individual we need
in addition (to obtain a time-homogeneous mechanism) to encode explicitly the time in the
state. That is, we finally choose
\be{agre8}
   \R\times\wh\Pi^{\mathbb{S}}
\ee
as the state space.
We also write $\rho_{\mathcal I}'\diamond (t,{\mathcal P}):=(t,\rho_{\mathcal I}'\diamond{\mathcal P})$. \sm

We are now ready to define the coalescent with rebirth.

\begin{definition}[Kingman-type coalescent with rebirth]\label{D:01}
Fix $t_0\in\R$. The Kingman-type coalescent with rebirth is a strong
$\R\times\hat{\Pi}^{\mathbb{S}(\CI)}$-valued Markov process
\be{grev1b}
   K^{\mathrm{birth}}
 =
   (K^{\mathrm{birth}}_u)_{u \geq t_0}, \quad t_0 \in \R,
\ee
whose initial condition $K^{\mathrm{birth}}_0:=(t_0,{\mathcal P}_0)$
satisfies
for all $i\in{\mathcal I}$,
\be{S}
   \#\big\{t\in\R:\,(i,t)\in {\mathcal P}_0\big\}
 <
   \infty,
\ee
and such that
for all finite subsets $\CI'\subseteq\CI$, the
  restricted process $\rho_{{\mathcal I}'}\diamond
  K^{\mathrm{birth}}$ is a $\R\times{\hat{\Pi}}^{\mathbb{S}'}$-valued
Markov chain which starts in $(t_0,\rho_{{\mathcal I}'}\diamond
  K^{\mathrm{birth}}_0)$ for some $K_0^{\mathrm{birth}}\in\hat{\Pi}^{\mathbb{S}(\CI)}$ such that
\begin{itemize}
\item the time coordinate grows at a deterministic speed one, and
\item given
the current
state $(t,{\mathcal P})\in\R\times\hat{\Pi}^{\mathbb{S}'}$
at time $t$,
each pair of partition elements $\pi_1,\pi_2\in{\mathcal P}$
merges into $\pi_1\cup\pi_2$ after an exponential waiting time with
rate $\gamma_{\mathrm{King}}>0$, and
at this time $t'>t$, instantaneously a new partition element
$\{(\mathrm{pr_{index}}(\kappa(\pi_1)\vee\kappa(\pi_2)),t')\}$ is born.
($\vee$ is the maximum taken in the sense of relation (\ref{vl19})).
\end{itemize}
\end{definition}\sm

\begin{proposition}[Existence and uniqueness in law]  \label{P:00}
\begin{itemize}
\item[{}]
\item[(a)] The Kingman-type coalescent with rebirth is a well-defined
pure jump process for every initial state with finitely many partition elements at time $t_0$.
\item[(b)] For every initial point in $\R\times\wh\Pi^{{\mathbb{S}(\CI)}}$ of the form
$(t_0,\{\{(i,t_0)\};\,i\in\N\})$, $t_0\in\R$, there exists a unique c\`adl\`ag process satisfying
the requirements of Definition~\ref{D:01}.
\end{itemize}
\end{proposition}\sm

\begin{proof}[{\bf Proof}] Proposition~\ref{P:00} is a special case of Proposition~\ref{P:02}. We therefore omit the proof at this point.
\end{proof}\sm

\begin{remark}\label{Rem:01}\rm
$(K^{\mathrm{birth}}_{u})_{u\in[a,b]}$
has the property that at each time $u>b$, for each $i\in\N$,
there is exactly
one (partition) element $\pi$ in $K^{\mathrm{birth}}_u$
with $\kappa(\pi)=(i,s)$, for some $s$.
Indeed the new individual $(i,s)$ will be born/introduced
at time $s$ only if a
partition element $\pi$ with label $\kappa(\pi)=(i,u)$, for some
$u<s$, coalesces at time $s$ with a partition element $\pi'$
such that $\kappa(\pi')<\kappa(\pi)$.
\hfill$\qed$\end{remark}\bi

%\subsubsection{The spatial coalescent with rebirth}
%\label{SubSub:spaco}
\noindent{\em Step~2 (Migration and coalescence with rebirth combined) }
In the case of the spatial coalescent with rebirth
all partition elements have in addition to an index and a birth-time
also a current location
that changes according to a random walk independently
over partition elements.
Fix again a countable index set $\CI\subseteq[0,\infty)$
and a {\em countable site space} ${M}$ which later will
be equal to $\Z^2$.
Then from any ${\mathcal P}\in \Pi^{\mathbb{S}(\CI)}$
one can form a {\em marked} partition
\be{vl9a}
   {\mathcal P}^M
 :=
   \big\{(\pi,L(\pi));\,\pi\in{\mathcal P}\big\},
\ee
by assigning to each partition element $\pi\in{\mathcal P}$,
its {\em location} $L(\pi)\in {M}$.
Denote the space of marked partitions in $\hat{\Pi}^{\mathbb{S}(\CI)}$ by
\be{vl10b}
   \hat{\Pi}^{\mathbb{S}(\CI),{M}}.
\ee

For all ${\mathcal I}'\subseteq {\mathcal I}$, recall from (\ref{Erhorestr}) the
restriction operator $\rho_{{\mathcal I}',{M}}$
from $\hat{\Pi}^{\mathbb{S}(\CI),{M}}$ to
$\hat{\Pi}^{\mathbb{S}({\mathcal I}'),{M}}$.
We say that a sequence
\be{e:coc}
   ({\mathcal P}_n)_{n\in\N}\quad
 \mbox{converges in }
   \hat{\Pi}^{\mathbb{S}(\CI),{M}}
\ee
if and only if the projections on the index component of
the restricted partitions $\rho_{\CI',M}\diamond {\mathcal P}_n$ converge
in the discrete topology
and their corresponding birth times and locations converge.
More precisely, we consider the topology generated by the metric (\ref{f:003}) through (\ref{f:005}).
%moreover the finitely many individuals with a particular index in
%$\CI'$ can be bijectively rearranged
%so that their corresponding birth-times and locations converge.

We are now prepared to define the spatial ${\mathcal
  I}$-coalescent with rebirth.
\begin{definition}[The spatial ${\mathcal I}$-coalescent with
  rebirth]\label{D:02}
%\begin{itemize}
%\item[{}] \hspace{3cm}
%\item[(a)]
The {\em spatial ${\mathcal I}$-coalescent with
  rebirth},
\be{0.4aa}
   (C^{\mathrm{birth}},L^{\mathrm{birth}})
 :=
   \big(t_0+t,(C^{\mathrm{birth}}_t,L^{\mathrm{birth}}_t)\big)_{t\ge 0},
\ee
is a strong $\R\times\hat{\Pi}^{\mathbb{S}(\CI),{M}}$-valued
  Markov process with c\`adl\`ag paths such that for all subsets ${\mathcal
  I}'\subseteq{\mathcal I}$ with
\be{Lip1a}
   \sum\nolimits_{\pi\in C_0,\rho_{{\mathcal I}'}(\pi)\not=\emptyset}
   \delta_{L^{\mathrm{birth}}_0(\pi)}
 \in
   {\mathcal E},
\end{equation}
and all initial birth times less than or equal to $t_0$,
the restricted process
$\rho_{{\mathcal I}'}\diamond (C^{\mathrm{birth}},L^{\mathrm{birth}})$
is a $\R\times\hat{\Pi}^{\mathbb{S}({\mathcal I}'),M}$-valued
    strong Markov particle system which undergoes the following
    three {\rm independent} transition mechanisms:
\begin{itemize}
\item {\rm Time growth }\; The time coordinate grows at deterministic rate one.
\item {\rm Migration }\;
The marks of the
partition elements perform independent
random walks.
\item {\rm Coalescence with rebirth}\;
Given the current state
$(t,\{(\pi,L^{\mathrm{birth}}(\pi));\,\pi\in{\mathcal P}\})\in\R\times\hat{\Pi}^{\mathbb{S}({\mathcal I}'),\Z^2}$,
each pair of partition elements $\pi_1,\pi_2\in{\mathcal P}$
merges into $\pi_1\cup\pi_2$ after an exponentially distributed waiting time
with hazard function given by the density
$\gamma\cdot\mathbf{1}\{L^{\mathrm{birth}}_{\boldsymbol{\cdot}}(\pi_1)=
L^{\mathrm{birth}}_{\boldsymbol{\cdot}}(\pi_2)\}$,
and at this random time $t'>t-t_0$, instantaneously the marked partition element
$\big(\{(\mathrm{pr_{index}}(\kappa(\pi_1)\vee\kappa(\pi_2)),t'{+t_0})\},
L^{\mathrm{birth}}_{t'-t_0}(\pi_1)\big)$ is created.
\end{itemize}
%\item[(b)] Consider the process defined above but allow
%rebirth only during the time interval $[a,b] \subset \R$.
%The process arising is denoted by
%\be{grev1cc}
%   \big(C^{\mathrm{birth}}[a,b],L^{\mathrm{birth}}[a,b]\big).
%\ee
%\end{itemize}
\end{definition}\sm

\begin{proposition}[The spatial ${\mathcal I}$-coalescent
  rebirth is well-defined]\label{P:02}
\begin{itemize}
\item[{}]
\item[(a)] The spatial ${\mathcal I}'$-coalescent with
is a well-defined particle system, for all ${\mathcal I}'$ satisfying (\ref{Lip1a}).
\item[(b)] For every initial point in
  $\R\times\wh\Pi^{\mathbb{S}(\CI),M}$ of the form
$(t_0,\{\{((i,t_0),L^{\mathrm{birth}}((i,t_0)))\};\,i\in\N\})$,
$t_0\in\R$, there exists a unique
c\`adl\`ag process satisfying
the requirements of Definition~\ref{D:02}.
\end{itemize}
\end{proposition}\sm

The proof will be given in the following subsection.\sm

\subsection{The look-down construction (Proof of Propositions~\ref{P:00} and~\ref{P:02})}
\label{Sub:look-down}

In this subsection we give the explicit construction of a version of the coalescent
and the coalescent with rebirth. For that purpose we will rely on the graphical representation
of the look-down process introduced first by Donnelly and Kurtz in \cite{DonKur96} and
generalized to the spatial setting in \cite{GreLimWin05}.
In the look-down construction we can link both the population model of locally
infinite population size in forward time and the coalescent starting with locally infinitely many patches in reversed time.

\begin{figure}
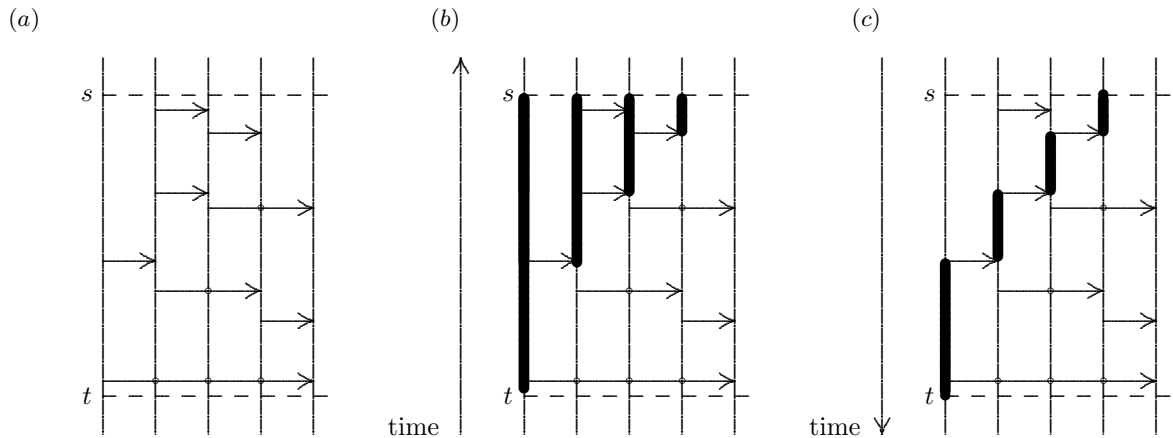

\beginpicture
\setcoordinatesystem units <.7cm,1cm>
\setplotarea x from -0.5 to 12, y from 0 to 5
\put{{\small $(a)$}} at -.5 5.5
\put{{\small $(b)$}} at 7.5 5.5
\put{{\small $(c)$}} at 15.5 5.5
\plot 1 0 1 5 /
\plot 2 0 2 5 /
\plot 3 0 3 5 /
\plot 4 0 4 5 /
\plot 5 0 5 5 /
\arrow <0.2cm> [0.375,1] from 1 0.7 to 5 0.7
\put{\tiny$\circ$} [cC] at 2 0.7
\put{\tiny$\circ$} [cC] at 3 0.7
\put{\tiny$\circ$} [cC] at 4 0.7
\arrow <0.2cm> [0.375,1] from 4 1.5 to 5 1.5
\arrow <0.2cm> [0.375,1] from 2 1.9 to 4 1.9
\put{\tiny$\circ$} [cC] at 3 1.9
\arrow <0.2cm> [0.375,1] from 1 2.3 to 2 2.3
%\arrow <0.2cm> [0.375,1] from 4 2.4 to 5 2.4
\arrow <0.2cm> [0.375,1] from 3 3.0 to 5 3.0
\put{\tiny$\circ$} [cC] at 4 3.0
\arrow <0.2cm> [0.375,1] from 2 3.2 to 3 3.2
\arrow <0.2cm> [0.375,1] from 3 4.0 to 4 4.0
\arrow <0.2cm> [0.375,1] from 2 4.3 to 3 4.3

\put{time} [cC] at 6.9 0.1

\plot 9 0 9 5 /
\thicklines
  \multiput {$\bullet$} at 9 4.45 *193 0 -.02  /
\plot 10 0 10 5 /
  \multiput {$\bullet$} at 9 4.45 *109 0 -.02  /
    \multiput {$\bullet$} at 10 4.45 *109 0 -.02  /
\plot 11 0 11 5 /
  \multiput {$\bullet$} at 9 4.45 *62 0 -.02  /
  \multiput {$\bullet$} at 11 4.45 *62 0 -.02  /
\plot 12 0 12 5 /
  \multiput {$\bullet$} at 12 4.45 *22 0 -.02  /
\plot 13 0 13 5 /
\arrow <0.2cm> [0.375,1] from 9 0.7 to 13 0.7
\put{\tiny$\circ$} [cC] at 10 0.7
\put{\tiny$\circ$} [cC] at 11 0.7
\put{\tiny$\circ$} [cC] at 12 0.7
\arrow <0.2cm> [0.375,1] from 12 1.5 to 13 1.5
\arrow <0.2cm> [0.375,1] from 10 1.9 to 12 1.9
\put{\tiny$\circ$} [cC] at 11 1.9
\arrow <0.2cm> [0.375,1] from 9 2.3 to 10 2.3
%\arrow <0.2cm> [0.375,1] from 20 2.4 to 21 2.4
\arrow <0.2cm> [0.375,1] from 11 3.0 to 13 3.0
\put{\tiny$\circ$} [cC] at 12 3.0
\arrow <0.2cm> [0.375,1] from 10 3.2 to 11 3.2
\arrow <0.2cm> [0.375,1] from 11 4.0 to 12 4.0
\arrow <0.2cm> [0.375,1] from 10 4.3 to 11 4.3

\arrow <0.2cm> [0.375,1] from 7.8 0 to 7.8 5
\arrow <0.2cm> [0.375,1] from 15.8 5 to 15.8 0
\put{time} [cC] at 14.9 0.1

\thinlines
\plot 17 5 17 0 /
\thicklines
  \multiput {$\bullet$} at 17 2.25 *85 0 -.02  /
\plot 18 0 18 5 /
  \multiput {$\bullet$} at 18 2.35 *41 0 .02  /
\plot 19 0 19 5 /
  \multiput {$\bullet$} at 19 3.95 *36 0 -.02  /
\plot 20 0 20 5 /
  \multiput {$\bullet$} at 20 4.45 *22 0 -.02  /
\plot 21 0 21 5 /
\arrow <0.2cm> [0.375,1] from 17 0.7 to 21 0.7
\put{\tiny$\circ$} [cC] at 18 0.7
\put{\tiny$\circ$} [cC] at 19 0.7
\put{\tiny$\circ$} [cC] at 20 0.7
\arrow <0.2cm> [0.375,1] from 20 1.5 to 21 1.5
\arrow <0.2cm> [0.375,1] from 18 1.9 to 20 1.9
\put{\tiny$\circ$} [cC] at 19 1.9
\arrow <0.2cm> [0.375,1] from 17 2.3 to 18 2.3
%\arrow <0.2cm> [0.375,1] from 12 2.4 to 13 2.4
\arrow <0.2cm> [0.375,1] from 19 3.0 to 21 3.0
\put{\tiny$\circ$} [cC] at 20 3.0
\arrow <0.2cm> [0.375,1] from 18 3.2 to 19 3.2
\arrow <0.2cm> [0.375,1] from 19 4.0 to 20 4.0
\arrow <0.2cm> [0.375,1] from 18 4.3 to 19 4.3

\put{$\bullet$} [cC] at 20 4.5
\put{$\bullet$} [cC] at 17 0.5

\setdashes
\plot 1 0.5 5.4 0.5 /
\plot 1 4.5 5.4 4.5 /
\put{$t$} [cC] at 0.7 0.5
\put{$s$} [cC] at 0.7 4.5

\plot 9 0.5 13.4 0.5 /
\plot 9 4.5 13.4 4.5 /
\put{$t$} [cC] at 8.7 0.5
\put{$s$} [cC] at 8.7 4.5

\plot 17 0.5 21.4 0.5 /
\plot 17 4.5 21.4 4.5 /
\put{$t$} [cC] at 16.7 0.5
\put{$s$} [cC] at 16.7 4.5
\endpicture
\caption{\label{fig:MM}
(a) a realization of the Poisson processes $M$, (b)
the set of all descendants up to time $s$ of the individual labeled $1$ at time $t$ is indicated in bold, (c)  the ancestral line of the individual $4$ alive at time $s$ is indicated in bold.
}
\end{figure}

In order to give the explicit construction based on a random graph we proceed as follows.
Fix a rate $\gamma>0$, and a non-empty countable set ${\mathcal I}$ referred to as the
{\em set of all individuals}.  Assume we are given
a collection
\be{agre21}
   \big\{(L^i_t)_{t \geq 0}, i\in
   {\mathcal I}\big\}
\ee
of the independent continuous time irreducible random walks on an Abelian group $G$. Then we can choose
a total order $\preceq$ on ${\mathcal I}$ such that for all $i\in{\mathcal I}$ and $x\in\Z^2$,
\be{orderreq}
  \#\big\{i'\preceq i:\,L^{i'}_0=x\big\}<\infty,\qquad\mbox{a.s.}
\ee

Let \be{MM}
   M
 :=
    \big\{M^{i,j}:\, i,j \in \CI; \;  i\preceq j\big\}
\ee
be a family of independent Poisson point processes on $\R^+$
with intensity measure $\gamma\,\mathrm{d}t$. The random collections in
(\ref{agre21}) and (\ref{MM}) are independent. This specifies our
probability space. Starting from $(L^i_0,0)$, $i \in \CI$, we follow the random walk
$(L^i_t)_{t \geq 0}$ and draw an arrow from $i$ to $j$ at time $t$ if
$t$ is a point of $M^{i,j}$ and
$L^i_t = L^j_t$. This defines a random graph embedded in $\R \times \CI$
with (random) marks in $\Z^2$, which is defined on our probability space.

For $s < t$ we say that $(i,s)$ and $(j,t)$ are connected by a {\em path} if in the
$\R \times \CI$ diagram we can move vertically without crossing the tip of
an arrow or horizontally along arrows from $(i,s)$ to $(j,t)$.
This means that in forward time we think of the points in $M^{i,j}$ as the times at which individual
$i$ is pushing out an individual $j$ from the population in order to replace it by
a new individual of its type.
We therefore call such a path a {\em line of descent} and $(j,t)$ a {\em descendent} of
$(i,s)$.

On the other hand we can reverse this path and say that the reverse
path is an {\em ancestral line} associating with $(j,t)$ its ancestor
$(i,s)$ at time $s$. In this path $(A^j_{s,t})_{t \geq 0}$ time now runs
backward, so that $A^j_{s,s}=j$. In reversed time we therefore
interpret the points in $M^{i,j}$ as the times at
which the ancestral lines of the
individuals $i$  and $j$ may {\em coalesce} to a common ancestral line.
For example, in Figure~\ref{fig:MM}(c) the ancestor back at time $t$ of the individual
which lives at time $s$ and corresponds to the fourth ancestral line is the individual which
corresponds to the first ancestral line.
If we define for each $j\in{\mathcal I}$ and $t\ge s\ge 0$,
\be{e:desc}
    \pi^s_j(t)
  :=
    \big\{i\in{\mathcal I}:\,A^{i}_{s,t}=j\big\},
\ee
we obtain that the partition element $\pi^s_j(t)$ as the set of all {\em descendants}
at time $s$ of the individual $j$ which lived at time $t$ in the past. For example, in
		   Figure~\ref{fig:MM}(b) the individual which lives at time $t$ and corresponds to the first ancestral line has $4$
		   descendants at time $s$ which correspond to the first four ancestral lines.

By condition (\ref{orderreq}) the forward construction is automatically
well-defined and hence the following key result holds:
\beL{L:01}{Ancestors are well-defined.}
For each $i\in{\mathcal I}$ and $s\ge 0$, there exists a unique function $(A^{i}_{s,t})_{t\ge 0}$
from $[s,\infty)$ into ${\mathcal I}$ with c\`adl\`ag paths such that $A^{i}_{s,s}=i$ and
\be{e:ances}
   A^{i}_{s,t-}\not =A^{i}_{s,t},\quad\mbox{if and only if}\quad t\in M^{A^{i}_{s,t},A^{i}_{s,t-}}\mbox{ and }L^{A^{i}_{s,t-}}_{t-}=L^{A^{i}_{s,t}}_{t-}.
\ee
\end{lemma}\sm

\begin{remark}[The look-down process and the spatial ${\mathcal I}$-coalescent]\rm

Construct the infinitely old population for the forward model in times
$(-\infty, s)$.
Recall from the look-down construction from Subsection~\ref{Sub:look-down} the notion
$A_{s,t}^i$ of the ancestor at time $t$ in the past of the individual $i$ which lives at time $s$ and the
notion $\pi^s_j(t)$ of the set of all descendants which are alive at time $s$ of the individual which lived in the past at time $t$.
\label{Rem:12}
Put $s=0$ and set
\be{e:partit}
   C_t
 :=
   \big\{\pi^0_j(t):\,j\in{\mathcal I},\pi_j(t)\not=\emptyset\big\}\in\Pi^{{\mathcal I}}.
\ee
Notice that if $\pi\in C_t$ then $A^{i}_{0,t}=A^{i'}_{0,t}$ for all $i,i'\in\pi$. Write therefore $A^{\pi}_{0,t}$ for
common ancestor of all individuals in $\pi$ at time $t$ in the past, and put
\be{e:loc}
   L_t(\pi)
 :=
   L_t^{A^{\pi}_{0,t}}.
\ee

Then the process $(C_t,L_t)_{t\ge 0}$ is the spatial ${\mathcal I}$-coalescent. Notice that the cadlag path property follows immediately from the choice of topology.
\hfill$\qed$
\end{remark}\sm

\begin{proof}[{\bf Proof of Proposition~\ref{P:02}}]
The proof of Assertion (a) is obvious and the proof of Part (b) will be given
with the look-down process we define next. \sm

(b) The {\em uniqueness} of the process is a direct consequence of the fact that all finite
sub-coalescents are uniquely determined and hence if the desired object exists
it must be unique.

In order to get
the {\em existence} of the process starting from a state containing countably
many individuals, we use once more the look-down construction.

Recall Lemma~\ref{L:01} the notion $A_{s,t}^i$ of the ancestor at time $t$ back
in the past of the individual $i$ which lives at time $s$.

We here let for each $t\ge 0$ and $j\in{\mathcal I}$,
$$\hat{\pi}_{(j,t)}:=\big\{i\in{\mathcal I};\,\exists\,s<t\mbox{ such that }A^j_{s,t}=i\big\},$$
and denote by
$$u^j_{t}(i):=\inf\{{s\le t};\,A^i_{s,t}=j\}$$
the {\em birth time} of the descendent $i\in\hat{\pi}_{(j,t)}$ of individual $j$ which lives at time $s$,
and put
$$C_t:=\Big\{\big\{(i,u_t^j(i));\,i\in\hat{\pi}_{(j,t)}\big\};\,j\in{\mathcal I}\Big\}.$$

For each $t\ge t_0$ and $j\in{\mathcal I}$, assign
\be{e:locreb}
   L_t\Big(\big\{\{(i,u_t^j(i));\,i\in\hat{\pi}_{(j,t)}\}\big\}\Big):=L^j_t.
\ee

Then $(t_0+s,(C_s,L_s))_{s\ge 0}$ is the spatial coalescent with rebirth.
It remains to show the {\em c\`adl\`ag path property}.

Notice first that the topology on $\hat{\Pi}^{\mathbb{S}^{{\mathcal I}}}$
can be metrized, for example, by the
metric $d^\mathbb{S}$ defined as follows: we let for each $N\in\N$,
\be{f:003}
   d^\mathbb{S}\big({\mathcal P},{\mathcal P}'\big)
 \ge
   2^{-N}
\ee
if and only if for ${\mathcal I}'_N\subseteq{\mathcal I}$ with $\#{\mathcal I}_N=N$,
\begin{itemize}
\item $\mathrm{pr}_{\mathrm{index}}\rho_{{\mathcal I}'_N}{\mathcal P}
 \not=
   \mathrm{pr}_{\mathrm{index}}\rho_{{\mathcal I}_N}{\mathcal
P}'$, or
\item $\mathrm{pr}_{\mathrm{index}}\rho_{{\mathcal I}'_N}{\mathcal
  P}=\mathrm{pr}_{\mathrm{index}}\rho_{{\mathcal I}'_N}{\mathcal
P}'$ and for any one to one map $\iota$ from $\mathrm{pr}_{\mathrm{index}}\rho_{{\mathcal I}'_N}{\mathcal
  P}$ onto $\mathrm{pr}_{\mathrm{index}}\rho_{{\mathcal I}'_N}{\mathcal
P}'$,
\be{f:005}
   \max\Big\{\big|\mathrm{pr}_{\mathrm{time}}(v)-
   \mathrm{pr}_{\mathrm{time}}(\iota(v))\big|:\,v\in S(\rho_{{\mathcal I}'_N,K}{\mathcal P})\Big\}
 \ge
   2^{-N}.
\ee
\end{itemize}
\sm

For all finite ${\mathcal I}\subseteq{\mathcal I}$, the restricted processes
$K^n:=\rho_{{\mathcal I}'}\diamond K^{\rm birth}$ is a pure jump process with c\`adl\`ag
paths. We will show that
\be{angrx}
  (K^n_s)_{s\ge 0}\nto(K_s^{\rm birth})_{s\ge 0},
\ee
in Skorohod topology, almost surely.

For that fix $T>0$. We will show that for all $(t_n)_{n\in\N}$ with
$(t_n)\downarrow t$, $K_{t_n}^n\nto K^{\rm birth}_t$ and for all $(t_n)_{n\in\N}$ with
$(t_n)\uparrow t$, $K_{t_n}^n\nto K^{\rm birth}_{t-}$, almost surely.
Indeed, if $(t_n)_{n\in\N}$ with
$(t_n)\downarrow t$ and $\varepsilon>0$ are given then there exists a
random $N=N(\varepsilon)$ such that for all $n\ge N(\varepsilon)$,
the set
$\cup_{1\le i\le j\le \lfloor
  1-\log_2{\varepsilon}\rfloor}M^{i,j}[t,t_n]$ is empty and therefore
$d^{\mathbb{S}} (K_{t_n},K_t^{\rm birth})<\varepsilon$ (recall
$d^{\mathbb{S}}$ from (\ref{f:003})). The other convergence
relation follows by a similar argument and the cadlag path property follows by the choice of the topology.
\end{proof}\sm

\section{Main results}
\label{Sub:descript}
We study the asymptotic behavior of the spatial coalescent with initial
configurations
concentrated on bounded regions
as the region and the time of observation both become large.
Our parameter tending to infinity will be $t$, size in the geographic space
will be measured on the scale $t^{\alpha/2}$ and the time at which we observe the
process is on the
scale $t^\beta$ with $\alpha$ and $\beta$ being the corresponding macroscopic
{\em space parameter} and  {\em time parameter}, respectively.
More precisely, set for
$\alpha\in(0,1]$ and  $t\ge 0$,
\be{rect}
   \Lambda^{\alpha,t}
 :=
   \big[-t^{\frac{\alpha}{2}},t^{\frac{\alpha}{2}}\big]^2\cap\Z^2,
\ee
to define the region where all the individuals will be placed initially
and then observe this process at time $t^\beta$. Note that we are interested
in $\Lambda^{\alpha,t} \uparrow \Z^2$ by letting $ t\to \infty$.

We consider three settings, (1) the spatial coalescent (without rebirth)
as process in the macroscopic time parameter $\beta$ for fixed
space parameter $\alpha$, (2)
the spatial coalescent as process in the macroscopic space
parameter $\alpha$ for fixed time parameter $\beta$,
and (3) the spatial coalescent with rebirth.
In all settings we state that certain functionals of
the spatial coalescent started from a configuration which contains
particles at each site of $\Lambda^{\alpha,t}$, and
which is observed at times $t^\beta$, for a $\beta\geq\alpha$,
converge to corresponding functionals of the Kingman coalescent
with or without rebirth.  \bi

\subsection{The spatial coalescent in the macroscopic time parameter}
\label{Sub:macrosc}

We are now in our setting (1) and we consider the various functionals
of the coalescent in different subsubsections containing each a theorem.

\subsubsection{The number of partition elements as a process indexed by the
  time parameter}
  \label{Sub:beta}
Recall from Definition~\ref{Def:00} the spatial coalescent $(C,L)$ on
$\Z^2$, and let $K$ be the Kingman coalescent.
Denote by
\be{Calphatrho}
   C^{\alpha,t,\rho},
\ee
$\alpha\in(0,1]$ and $t\ge 0$,
the spatial coalescent that starts in a Poisson configuration
with either intensity $\rho\in(0,\infty)$ or
with intensity ``$\rho = \infty$'', i.e.\ with initially countable
infinitely many particles at each site of $\Lambda^{\alpha,t}$.
We refer to these processes as to the {\em $\alpha$-coalescent}.

\begin{remark}\label{Rem:04}\rm
The case $\rho<\infty$ is used in the study of
the so-called interacting Moran models, while the case $\rho=\infty$ is
needed to analyze its diffusion limit, the so-called
Fisher-Wright diffusions, or in the setting of infinitely many types,
the so-called interacting Fleming-Viot processes.
See \cite{CoxGre91} and \cite{DGV95} for more on these processes.
\hfill$\qed$
\end{remark}\sm

The following result
states the convergence of the number of partition elements of the
$\alpha$-coalescent observed at time $t^\beta$ in the space
$D([\alpha',\infty),\N)$ of c\`adl\`ag functions on $[\alpha^\prime, \infty)$
with values in $\N$, equipped with the Skorohod topology, where $\alpha'>\alpha$.
Here and in the remainder of the paper, for any (marked) partition ${\mathcal P}$, we denote by
$\# {\mathcal P}$ the number of equivalence classes in ${\mathcal P}$.

\begin{theorem}[Number of partition elements as
  processes in $\beta$; $\rho<\infty$]\label{T:01}
Fix $0<\rho<\infty$ and consider the spatial coalescent and the Kingman  coalescent
for the same coalescence parameter $\gamma >0$. Then for all $\alpha'>\alpha>0$,
\be{agr5}
   \CL\big[(\# C^{\alpha,t,\rho}_{t^\beta})_{\beta\in[\alpha',\infty)}\big]
 \Tto
   \CL\big[(\# K_{\log(\beta/\alpha)})_{\beta \in [\alpha',\infty)}\big].
\ee

If $(C_0,L_0)\in\hat{\Pi}^{{\mathcal I},\Z^2}$ is such that
in addition to
(\ref{Lip1})
the following generalization holds:
\be{agr5b}
   \sup_{z \in \Z^2}{\bf E}\Big[
   \#\big\{\pi\in C_0,\;L_0(\pi)=z\big\}\Big]
 <
   \infty,
\ee
and
\be{agr5c}
 \# \{\pi \in C_0, L_0(\pi) \in \Lambda^{1,t}\} \la \infty
 \mbox{ in probability as } t \to \infty
 \ee
then (\ref{agr5}) holds.

\end{theorem}\sm

More generally a careful reader will note that
the r.h.s. in (\ref{agr5}) does not involve the parameter $\rho$
and hence the scaling limit does not depend on $\rho$.
Indeed the more general statement shows that there is very little dependence between
the initial state and the scaling limit.

\subsubsection{The number of partition elements as a process indexed by the
  time parameter; infinite intensity}
  \label{Sub:betarho}

We next turn to $\rho=\infty$.
This case arises if one studies the genealogies in a model corresponding to
interacting measure-valued Fleming-Viot diffusions. These models are limits of the
spatial Moran model as the number of individual
per site tends to $\infty$. The genealogy of the limiting model can be
represented by the spatial coalescent starting with countable many
particles at each site, see \cite{FleGre94}.

In this situation the total number of initial individuals (partition elements) does not come down
from infinity in positive time (compare \cite{ABHL06}) since partition elements can escape
into empty space. However we will show that the fraction
of the partition elements  which do escape quickly is small and
its relative frequency in the total population is in fact 0 and therefore
they can be neglected.

\begin{remark}\rm
The frequency of a certain property in the population is here defined by taking the ''n-smallest''
in the order individuals and counting how many of them have the property in
\label{Rem:13}
question. Then normalizing by $n$ and letting $n \to \infty$ gives the
frequency of the property. The limit exists using de Finetti if our property is
a function of the individual which generates an exchangeable array if we
observe occurrence or non-occurrence of the property.
\hfill$\qed$
\end{remark}\sm

To make our approach to the case $\rho = +\infty$  precise, assume we are given a realization of
$(C_s,L_s)_{s\ge 0}$ with $C_0:=\{\{{i}\};\,{i}\in{\mathcal I}\}$
which is constructed from collections of independent random walks for
migration $\{(L^i_s)_{s \geq 0}, i \in \CI\}$
and Poisson point processes $\{M^{i,j}; i<j\}$.
We construct now an increasing collection of sub-coalescents of this spatial coalescent,
namely we remove every individual in the original configuration
for which $L^i$ jumps   before a given time $\delta>0$.
Then we start the process in this new sub-configuration.
This gives the sub-coalescent of $(C^{\alpha,t,\infty}_s)_{s \geq 0}$ denoted by
\be{calphatrh2}
   (C^{\alpha,t,\infty,\delta}_s)_{s \geq 0}.
\ee
\sm

Note that $C^{\alpha,t,\infty,\delta}_s \uparrow  C^{\alpha,t,\infty}_s$ a.s.
$\delta \downarrow 0, s \geq 0$ in our topology.
The following result is the analogue of Theorem \ref{T:01} for
$\rho=\infty$.
\begin{theorem}[Number of partition elements  as
  processes in $\beta$; $\rho=\infty$]
Let $\rho=\infty$, and $\delta>0$.
Then for all $\alpha'>\alpha>0$, \label{T:02}
\be{agr5d}
   \CL\big[(\# C^{\alpha,t,\infty,\delta}_{t^\beta})_{\beta
   \in[\alpha',\infty)}\big]
 \Tto
   \CL\big[(\# K_{\log(\beta/\alpha)})_{\beta \in [\alpha',\infty)}\big].
\ee
\end{theorem}\sm

\begin{remark}\label{Rem:05} \rm Notice the following:
\begin{itemize}
\item[(i)] The right hand side of (\ref{agr5b}) does {\em not} depend
on $\delta>0$.
\item[(ii)] The frequency of the individuals in
  $C_0^{\alpha,t,\infty}$
that are not contained in
$C_0^{\alpha,t,\infty,\delta}$
tends to zero, as $\delta\to 0$.
Hence the theorem describes the behavior of
the coalescent's initial population of individuals
(partition elements) with exception of a subset of frequency $0$.
\hfill$\qed$
\end{itemize}
\end{remark}\sm

\begin{remark}\rm
Proving results for the system with the exception of a set of
frequency $0$ of initial individuals  important if one
\label{Rem:14}
anticipates describing the genealogy of the Fleming-Viot process by a weighted $\R$-tree since then
one gets convergence in the canonical Gromov-weak topology, as in \cite{GrePfaffWin}.
\hfill$\qed$
\end{remark}\sm

\subsubsection{The number of partition elements as a process indexed by the
  time parameter; refinement}
  \label{Sub:betarhorefine}
The next goal is to extend the results in Theorems \ref{T:01} and \ref{T:02}
to the case where $\alpha^\prime=\alpha$.
Let $\N$ be equipped with the discrete topology and denote by
\be{barN}
   \bar{\N}
 :=
   \N\cup\{\infty\}
\ee
its one point
compactification. This means that  a sequence
$(n_k)_{k\in\N}$ with values in $\bar{\N}$ converges in $\bar \N$ if either $(n_k)_{k\in\N}$ is a
convergent sequence in $\N$, or $(n_k)_{k\to\infty}$ diverges to
infinity.
%The two theorems above do not make sense for
%$\alpha^\prime = \alpha$ unless one
%introduces a new topology on the augmented set $\N\cup \infty$,
%for example generated by the following metric:
%\be{Emetricaug}
%\dm(n,m) = \sum_{i=\min\{n,m\}+1}^{\max\{n,m\}} \frac{1}{2^i}, \ d(n,\infty)= \sum_{i=n+1}^\infty \frac{1}{2^i}.
%\ee

Now we can consider the processes
\be{grev25}
   \big(\# C^{\alpha,t,\rho}_{t^\beta}\big)_{\beta\geq\alpha}, \,
   \big(\#
   C^{\alpha,t,\infty,\delta}_{t^\beta}\big)_{\beta\geq\alpha},
   \,\mbox{ and }
   \big(\# K_{\log{\beta/\alpha}}\big)_{\beta\geq\alpha}
\ee
in the Skorokhod space
$D([\alpha,\infty),\bar{\N})$.
For brevity, and in mind of future applications (see
\cite{glw2,glw3}),
we will consider only particular initial configurations.

\begin{theorem}[Convergence to the entrance law] Fix $\alpha>0$.
\begin{itemize}
\item[{}]
\item[(i)] Assume that the initial configuration is either a Poisson
  process with intensity $\rho$, or a
Bernoulli field with success probability $p\in(0,1]$,
for both choices we write $(C^{\alpha,t}_s, L^{\alpha,t}_s)_{s \geq 0}$
for the corresponding coalescent.
Then \label{T:03}
\be{agr5ent}
   \CL\big[(\# C^{\alpha,t}_{t^\beta})_{\beta\in[\alpha,\infty)}\big]
 \Tto
   \CL\big[(\# K_{\log(\beta/\alpha)})_{\beta\in[\alpha,\infty)}\big].
\ee
\item[(ii)] For each $\delta>0$,
\be{agr5bent}
   \CL\big[(\# C^{\alpha,t,\infty,\delta}_{t^\beta})_{
   \beta\in[\alpha,\infty)}\big]
 \Tto
   \CL\big[(\# K_{\log(\beta/\alpha)})_{\beta\in[\alpha,\infty)}\big].
\ee
\end{itemize}
\end{theorem}
%\begin{remark}\rm
%It turns out that we can prove that in fact
%not only the number of partition elements but partitions themselves
%converge in law modulo some random
%bijection between the basic set, ${\mathcal I}_{t^\alpha}$,
%of the $\CI^\alpha$-coalescent at time $t^\alpha$ and \label{Rem:11}
%the basic set $\N$ of the Kingman coalescent.
%\hfill$\qed$
%\end{remark}\sm\bi

\subsection{Spatial coalescent as a function of
 macroscopic spatial parameter $\alpha$}
  \label{SubSub:alpha}

We are now in our setting 2 and now the coalescent with rebirth occurs as limit object.
Fix $\alpha\in(0,1]$. Consider $C^{\alpha,t}$ the spatial coalescent on
$\Z^2$ but restricted to individuals initially in $\Lambda^{\alpha,t}$.
Let ${\mathcal I}^{\alpha,t}$ be the set of
individuals  initially placed in $\Lambda^{\alpha,t}$.
Then we can consider for every $t$ the collection of {\em sub-coalescents}
\be{grevcorr}
(\rho_{\CI^\alpha} \diamond C^{1,t})_{\alpha \in [0,1]}.
\ee
Notice that (equality in distribution)
\be{grevcorr2}
   C^{\alpha,t}
 \stackrel{\mathrm{d}}{=}
   \rho_{{\mathcal I}^\alpha}\diamond C^{1,t}
\ee
and,  of course, in (\ref{grevcorr}) the objects for different $\alpha$
all live on one probability space and they
are {\em coupled as sub-coalescents} of $C^{1,t}$.
Next we give a limiting object for partition element numbers.
Fix $0\le\alpha_l<\alpha_u\le1$.

Recall from Definition~\ref{D:01} the Kingman-type coalescent
$K^{\mathrm{birth}}$ with
rebirth and denote by
\begin{equation}\label{e:Kbab}K^{\mathrm{birth}}[a,b],\qquad -\infty\le a<b\le\infty,
\ee
the Kingman-type coalescent
with rebirth during the time interval $[a,b]$ only.

Next we give a limiting object for partition element numbers.
Fix $0\le\alpha_l<\alpha_u\le1$ and consider
the Kingman-type coalescent
$K^{\mathrm{birth}}[\log{\alpha_l},\log{\alpha_u}]$ with
rebirth during the time interval $[\log{\alpha_l},\log{\alpha_u}]$
only. We start the Kingman coalescent at time $\log \alpha$
(with $\CI=\N$) and we are interested in the latter, evaluated at time $0$.

\begin{remark}\rm
Recall the order relation (\ref{vl19})
that was used in the construction of a
particular realization of the process $K^{{\rm birth}}$.
Moreover, one could naturally order partition elements
within a partition according to their leading indices.
In the Definition~ (\ref{agr5e}) below, for reasons that will become
apparent later, we are introducing implicitly a ``reordering according to age''.
Note that, formally speaking, it is not a priori clear \label{Rem:15}
that the earliest born element of a partition element exists.
\hfill$\qed$
\end{remark}\sm

Recall
from (\ref{label}) the label $\kappa(\pi)$ of a partition element
$\pi$
and from (\ref{proj}) the projection maps which send the individual
$(i,t)$ to its
index and birth time.
For $\alpha\in [\alpha_l,\alpha_u]$, define:
\be{agr5e}
   N_\alpha
 :=
   \#\big\{\pi\in K^{\mathrm{birth}}_0[\log{\alpha_l},\log{\alpha_u}]:\,
  \mbox{ there exists } (s,i)\in \pi \mbox{ such that }
s \le\log(\alpha)\big\}.
\ee
and refer it to as the number of partition elements of
$K^{\mathrm{birth}}_0[\log{\alpha_l},\log{\alpha_u}]$
that are born before time $\log{\alpha}$.

The following result describes the asymptotic joint law of
the sizes of sub-coalescents
$(\# \rho_{{\mathcal I}^\alpha}C^{1,t})_{\alpha \in [\alpha_l,\alpha_u]}$
observed at time $t$, as $t \to \infty$.

\begin{theorem}[Convergence as processes in $\alpha$]
Fix $0\le\alpha_l<\alpha_u< 1$.
\begin{itemize}
\item[(i)] For all $\rho \in [0,\infty)$, \label{T:04}
\be{agr6ax}
   \CL\big[(\#\rho_{{\mathcal I}^{\alpha,t}}\diamond C^{1,t,\rho}_t\big)_{\alpha\in[\alpha_l,\alpha_u]}\big]
 \Tto
   \CL\big[(N_{\alpha})_{\alpha \in [\alpha_l,\alpha_u]}\big].
\ee
\item[(ii)] For all $\delta>0$,
\be{agr6axx}
   \CL\big[(\#\rho_{{\mathcal I}^{\alpha,t}}\diamond C^{1,t,\infty,\delta}_t\big)_{\alpha\in[\alpha_l,\alpha_u]}\big]
 \Tto
   \CL\big[(N_{\alpha})_{\alpha \in [\alpha_l,\alpha_u]}\big].
\ee
\end{itemize}
\end{theorem}

\begin{remark}\rm
\label{Rem:05b}
Notice that since
  $N_\alpha \to \infty$, as $\alpha \to 1$, and
  $\# C^{1,t,\rho}_t \to\infty$, as $t \to \infty$, the result holds
  also for $\alpha_u=1$. However, in order to rigorously 
include  $\alpha_u=1$ in the statement we would again have to consider the
  one point compactification of $\bar{\N}$ and apply similar
  techniques as in the proof of Theorem~\ref{T:03}.
\hfill$\qed$
\end{remark}\sm

\subsection{Rescaling the spatial coalescent with rebirth}
  \label{ss.spatialrebirth}

We are now in the setting 3.
Recall from Definition~\ref{D:02} the spatial
coalescent $\big(C^{\textrm{birth}}, L^{\textrm{birth}}\big)$
with rebirth.
Fix $\alpha\in(0,1)$ and $t>1$.
At time $t$ we observe the spatial coalescent with
rebirth which started
in Poisson configuration on $\Z^2$ with intensity
$\rho\in(0,\infty]$  at time $0$. In particular, even if $\rho<\infty$, the total number of initial partition elements is infinity.
Note that if $\rho<\infty$, all the partition
elements which ``die'' due to coalescence get replaced.
Hence during $[0,t]$ the configuration of locations on $\Z^2$
of partition elements remains Poisson.

Observe in the spatial coalescent with rebirth at a late time $t$ the partition
elements which are observed in a box $\Lambda^{\alpha,t}$ at some times $s_1,...,s_m \leq t$, where $m\in\N$, which forms a sub-coalescent. How many partition elements has this sub-coalescent currently in the limit as $t\to\infty$?
This question is also of interest since this sub-coalescent arises as a dual object in resampling models
if one considers the configuration in macroscopic time-space windows, and is explained in
Remark \ref{Rem:06} below.
In view of the previous scaling results we look at the system in times
$t^u$, with $u \in (0,1)$ the macroscopic time parameter and then let
$t \to \infty$. Since the times $t^u, t^{u^\prime}$ for $u^\prime \neq u$
separate, we cannot use a continuous macroscopic time parameter in our
analysis. We have to discretize.

For $m\in\N$, fix parameters
$\alpha < u_1<u_2<\ldots<u_m < 1$.
We are now
interested in the asymptotic behavior, as $t\to\infty$, of the
number of those partition elements observed in the population
at time $t$
which were located in $\Lambda^{\alpha,t}$ at an $m$-tuple of time points of the form $t^{u}$,
where $u \in \{u_1,\ldots,u_m\}$.
For $u\in(\alpha,1]$, we therefore put
\be{grev25a}
\begin{aligned}
   &N^{\alpha,t,\vec{u},\rho}_u
  \\
 &:=
   \#\big\{ \pi \in C^{\textrm{birth}}_t:\,
   \exists i \in \{1,\ldots,m\} \mbox{ s.t. } u_i \leq u, \,
   \mathrm{pr_{time}}(\kappa(\pi)) \leq t^{u_i}, \,
   L^{\textrm{birth}}_{t^{u_i}} (\pi) \in \Lambda^{\alpha,t}\big\}.
\end{aligned}
\ee
%\be{grev25a}
%N^{\alpha,t,\rho}_u = \# \{ \pi \in C^{\textrm{birth}}_t
%: pr_{\textrm{time}}(x(\pi) \leq t^u , \quad
%L^{\textrm{birth}}_{t^u} (\pi) \in \Lambda^{\alpha,t}\}.
%\ee

The dependence on $u_1,\ldots,u_m$ in the above definition is recorded by the
third superscript $\vec{u}$.
We chose to form the vector $\vec{u}=(u_1,\ldots,u_m)$ out of notational convenience.
Similarly we will write below
$\log{\vec{u}}$ for the vector $(\log{u_1},\ldots,\log{u_m})$.
We keep the dependence on $\vec{u}$
in mind, yet we omit it sometimes from the
notation by setting $N^{\alpha,t,\rho}\equiv N^{\alpha,t,\vec{u},\rho}$.
The parameter $u$ (in the subscript) will play the r{\^o}le of the new time index running
in $(\alpha, 1)$.

As before, if we consider $\rho=\infty$, we
need to observe
$\delta$-thinnings of our spatial coalescent with rebirth. We will
denote the corresponding functionals by
\be{agre18}
N^{\alpha,t,\infty,\delta}\equiv N^{\alpha,t,\vec{u},\infty,\delta}.
\ee

In order to study the asymptotic behavior of a suitable rescaling of
$N^{\alpha,t,\rho}_u$, we introduce a limit object, which we call
the {\em family of merging coalescents},
a collection of coalescents which start at specified times to
interact by coalescence, and which we denote by
\[
(K^{\mathrm{mer},\log{\vec{u}}}_s)_{s\ge\log{u_1}}.
\]
Here $1<\log{u_1}<\ldots<\log{u_m}<\infty$ is a given
sequence of merging times, at which the ``inter-coalescing'' of partitions belonging to
two or more different coalescents is enabled, as described precisely below.
Note that the coalescent structure is of Kingman-type, that is,
only pairs of partition elements (for which the coalescence is enabled)
coalesce at a constant rate.
%For that consider
%a collection
%\be{grev25b}
%   \big\{K^{\log{\alpha_i}};\,i\in\{1,\ldots,m\}\big\}
%\ee
%of independent Kingman coalescents starting with countably many
%partition elements.
The process $(K^{\mathrm{mer},\log{\vec{u}}}_s)_{s\ge\log{u_1}}$
is $\Pi^\N$-valued and evolves informally as follows.

We consider $m$ copies of the Kingman coalescent
$\{K^i,\,i\in\{1,\ldots,m\}\}$, where the $i^{\mathrm{th}}$ copy is
initially started in the configuration $\{\{km+i-1\}:\,k\in\N\}$ and runs
from time $\log{u_1}$ until time
$\log{u_i}$ independently from the others, but after time $\log{u_i}$ its
partition elements coalesce mutually as well as with
the partition elements of $\{K^j,\,j\in\{1,\ldots,i-1\}\}$.

A realization is constructed as follows.
The family of merging coalescents process
starts at time $\log{u_1}$ in
$\{\{n\};\,n\in\N\}$, and given a time $s\geq \log{u_1}$ two partition
elements $\pi_1,\pi_2\in
K^{\mathrm{mer},\log{\vec{u}}}_s$ with
\be{e:0012}
   \kappa(\pi_l)
 =
   mk_l+n_l, \mbox{ for some $k_l\in\Z$ and $n_l\in\{0,1,\ldots,m-1\}$, $l=1,2$},
\ee
in which case we write for $n_l:=[\kappa(\pi_l)]_{\mathrm{mod} m}$, coalesce at rate
$1\{n_1=n_2\}+1\{n_1\not=n_2,s\ge\log{u_{n_1}}\vee\log{u_{n_2}}\}$.
Upon coalescing the new partition inherits, as usual, the smaller label
where we define that for $\kappa(\pi_l)$, $l=1,2$, of the form (\ref{e:0012})
\be{grev25a2}
   \kappa(\pi_1)\leq\kappa(\pi_2) \mbox{ iff } n_1 \leq n_2 \mbox{ or } n_1=n_2, \, k_1\leq k_2.
\ee

%Now let each of these evolve independently up to time $u$ from where a
%coalescence with partition elements for $v < u$ takes place. This object
%we call
%\be{grev25c}
%K^{\mathrm{mer}} = (K^{\mathrm{mer}}_t)_{t \in \R},
%\ee
%the Kingman coalescent with immigration.
Note that $(K^{\mathrm{mer},\log{\vec{u}}}_s)_{s\ge\log{u_1}}$ can be
coupled with the coalescent with rebirth on $\Z^2$, where in the latter at each time
countably many individuals are reborn (immigrate back into the system), so that both can be constructed in such a way that the number of
partition elements in $(K^{\mathrm{mer},\log{\vec{u}}}_s)_{s\ge\log{u_1}}$ is
almost surely smaller than the number of partition elements in the coalescent with rebirth.
In particular, $(K^{\mathrm{mer},\log{\vec{u}}}_s)_{s\ge\log{u_1}}$ is well-defined,
and its number of partition elements is finite at all times $s> \log{u_1}$, almost surely.

%Let
%\be{grev25d}
%(K^{\mathrm{mer}} ([a,b])
%\ee
%be the process with immigration taking place onlygrev25
 %from $u \in [a,b]$.
%Now define
%\be{grev25e}
% (\wt K^{\mathrm{mer}}_s = K^{\mathrm{mer}}_{\log(\frac{s}{\alpha})}
%\ee
%and then for $0 < \alpha <\beta<1$
%\be{grev25f}
% (\wt K^{\mathrm{mer}}_0 ([\log \alpha, \log \beta])
%\ee
%is the quantity of interest for us.

Put (note $0<\log(u_1/\alpha)<\ldots<\log(u_m/\alpha)<\log(1/\alpha)$):
\be{Nimmalpha}
%N^{\mathrm{mer}}_{u_i}
%\equiv
N^{\mathrm{mer},\vec{u}/\alpha}_{u_i}
:=
   \#\big\{\pi\in K^{\mathrm{mer},\log(\vec{u}/\alpha)}_{\log(1/\alpha)}:\,
   [(\kappa(\pi))]_{\mbox{mod}(m)} \le i-1\big\}.
\ee

\begin{theorem}[Asymptotics of coalescent with rebirth]\label{T:05}
Fix  $0<\alpha<\beta<1$.
\begin{itemize}
\item[(a)] If $\rho<\infty$, then for all $m\in\N$ and $\alpha\le u_1,\ldots,u_m\le\beta$,
\be{agr6y}
   \CL\Big[\big(N^{\alpha,t,\rho}_{{u_1}},\ldots,
   N^{\alpha,t,\rho}_{{u_m}}\big)\Big]
 \Tto
   \CL\Big[\big(N^{\mathrm{mer},\vec{u}/\alpha}_{u_1},\ldots,N^{\mathrm{mer},\vec{u}/\alpha}_{u_m}\big)\Big].
\ee
\item[(b)] If $\rho =\infty$,
then for all $\delta>0$, $m\in\N$ and $\alpha\le u_1,\ldots,u_m\le\beta$,
\be{agr6y2}
   \CL\Big[\big(N^{\alpha,t,\infty,\delta}_{{u_1}},\ldots,
   N^{\alpha,t,\infty,\delta}_{{u_m}}\big)\Big]
 \Tto
   \CL\Big[\big(N^{\mathrm{mer},\vec{u}/\alpha}_{u_1},\ldots,N^{\mathrm{mer},\vec{u}/\alpha}_{u_m}\big)\Big].
\ee
\end{itemize}
\end{theorem}\sm

\begin{remark}[Space-time cluster formation]\rm
\label{Rem:06}\rm
As already indicated,
the spatial coalescent with rebirth describes the space-time
genealogy of the interacting Moran models.
To make this more precise, let us fix some large $t$
and introduce the reversed time $s^\leftarrow(s)\equiv s_t^\leftarrow(s):=t-s$.
Then, provided that the
original configuration of particles is Poisson,
and that the particles evolve according to the
enriched interacting Moran models
in forward time (where the types that die due to resampling are
kept as fossils), then their paths observed in reversed time
evolve according to the spatial coalescent with rebirth.
Moreover, a resampling event that occurs at time $t-s$ corresponds to
a unique rebirth event occurring at time $s_t^\leftarrow(s)$.

In this way,
Theorem~\ref{T:05} plays an important r{\^o}le in the study of the
space-time cluster formation of the interacting Moran models on
$\Z^2$. 
Namely, assume that at the initial time $0$
each individual (particle) carries its own type.
The following questions arise naturally in this context: 
if we fix a time $t>0$ and a large window $W$ of
observation, how far back in time do we have to look so
that most of the population present in $W$ at time $t$ has a single ancestor
and hence carries a single type (color)?;
how does this information
change if
the population is sampled at several time instances from the same window $W$?

The natural time scale for answering these questions is the
logarithmic time scale.
Fix $\alpha\in(0,1]$ and choose, out of convenience, the 
{\em $\alpha$-box} $\Lambda^{t,\alpha}$ as the window of observation.
%We are interested in the {\em age} of this $\alpha$-cluster.
Fix $u\in(0,1]$, and observe the subpopulation,
located in the $\alpha$-box at time $t$,
during the time interval $[t-t^{u},t]$.
It turns out that for large $t$ and
$u\le\alpha$, with overwhelming probability we find in this subpopulation a certain non-trivial ($\geq 2$) number of types
at time $t-t^{u}$, and a non-trivial
number ($\geq 2$) of these types are still
visible in the $\alpha$-box at time $t$.
However, if $u>\alpha$,
the event that only one of the types observed at time $t-t^{u}$
is visible in the $\alpha$-box at time $t$,
has positive probability, for $t$ large.
This is equivalent to saying
that with positive probability, all the
individuals observed in the
$\alpha$-box at time $t$ have a common ancestor
among the particles observed at time $t-t^{u}$, for $t$ large.
If this happens, we say that the {\em age (on logarithmic scale)} of our chosen
subpopulation is at most
$u$, since one type is carried by
a substantial fraction of the subpopulation,
and its original carrier can therefore be considered as the ancestor.

More generally, fix $m\geq 1$ and $u_1,\ldots,u_m$ such that
$0<u_1<\ldots <u_m\leq 1$.
During the time interval $[t-t^{u_m},t]$, consider 
the joint evolution
of $m$ different Moran model subpopulations 
where the ``$0$th'' subpopulation consists of
particles present in the $\alpha$-box at 
time $t$, and for $i=1,\ldots,m-1$, the 
$i$th subpopulation consists of particles
present in the $\alpha$-box at time $t-t^{u_i}$. 
By reversing time %in the same way as above,
all the interesting information about their joint genealogy is
expressed precisely in terms of the quantities 
$(N^{\alpha,t,(u_1,...,u_m),\rho}_{u_1},...,N^{\alpha,t,(u_1,...,u_m),\rho}_{u_m})$
as defined in (\ref{grev25a}). 
For example, the event on which the latter vector takes value $(1,...,1)$
is precisely the event that all the individuals (in $m$ subpopulations
combined) have a common ancestor at time $t-t^{u_m}$.
%all partition elements
%we have the same $\alpha$-cluster
%at times $t,t-t^{u_1}, \cdots, t-t^{u_{m-1}}$ persisting with age $u_m$ (in %logarithmic scale).
\hfill$\qed$
\end{remark}\sm

\begin{proof}[{\bf Outline}] The rest of the paper is organized as follows:
In Section~\ref{S:prelim} we recall and extend some basic facts on
coalescents on $\Z^2$, and in Sections~\ref{S:ASP} and~\ref{S:ADPST}
we provide
the asymptotic analysis of coalescents which allows us to prove
Theorems~\ref{T:01} and~\ref{T:02} in Section~\ref{S:ltsas}, and
Theorems~\ref{T:04} and~\ref{T:05} in Section~\ref{S:ConvSpatScale}.
Section \ref{S:proofmobo} contains
the proof of a moment estimate on the number of partition elements.
\end{proof}\sm

\begin{proof}[{\bf Result and Problem History}]
Here we give some information concerning the history of the problems
treated in this paper.
In the setting of instantaneous coalescence for simple random walks on
$\Z^2$, i.e., two partition
elements coalesce immediately when the hit the same site,
Lemma \ref{L3.1mar}
was proved in \cite{CoxGri86},
and Proposition~\ref{P:tight}
in \cite{BramCoxGri86}.
Propositions \ref{P3.1} and \ref{P3.1a} and Lemma
\ref{L3.1b}, are to the best of our knowledge
novel in the setting of any spatial coalescent model on $\Z^2$.
Due to the applications we have in mind (using duality with the
IMM and IFWD) in the subsequent papers, we are primarily
interested in the spatial (and delayed)
coalescents, and therefore the results are phrased and proved
in the current setting. However, it is important to note that the
arguments, and therefore statements,
in Section~\ref{Sub:descript} remain to hold in the
setting of \cite{BramCoxGri86} and \cite{CoxGri86}.
\end{proof}\sm

\section{Preliminaries}
\label{S:prelim}
In this section we present several basic techniques on coalescents
and present the key properties of random walks which we will need for our
subsequent arguments. We first state in Subsection~\ref{Sub:convent}
some notational conventions which will be used throughout the rest of
the paper. In
Subsection~\ref{Sub:Erdos} we recall a famous result by Erd\"os-Taylor
which gives the asymptotics of the hitting time of a planar random
walk. In Subsection~\ref{Sub:asympexch} we state the asymptotic
exchangeability for the spatial coalescent on $\Z^2$. In
Subsection~\ref{Sub:monoton} we recall some consequences of monotonicity
properties.

\subsection{Notational conventions}
\label{Sub:convent}
In the rest of
the paper we often use the following convention concerning
notation.
\begin{itemize}
\item For functions $g,h:[0,\infty)\to\R$, we write
$g(t)=O(h(t))$ or $g(t)=o(h(t))$ if and only if
$\limsup_{t\to \infty}\tfrac{g(t)}{h(t)}<\infty$ or
$\lim_{t\to \infty}\tfrac{g(t)}{h(t)}=0$, respectively.
\item
For a set $A$,  we denote by
$A^c$ its complement (with respect to the natural superset, determined
by the context).
\item Recall $\bar{\N}$ from (\ref{barN}), and denote for a finite or countable
set $A$ by
$\# A\in\bar{\N}$ the number of elements in $A$.
\item
If $a,b\in\R$, let $a\wedge b$ denote the minimal, and $a \vee b$
the maximal element of $\{a,b\}$.
\item
Poisson($\rho$) random variable (or distribution) has intensity (rate, expectation) $\rho$.
\item
For a partition $\CP$, recall that $\# \CP$ denotes the number of partition
elements of $\CP$.
\item
If $\CP$ is a partition then we write $i\sim_{\CP} j$ if $i$ and $j$ belong to the same partition element of $\CP$.
If $(\CP_t, \, t\geq 0)$ is a partition-valued process then
$i \sim_{\CP_t} j$ will be sometimes abbreviated as $i \sim^t j$.
\end{itemize}

\subsection{Erd\"os-Taylor formula}
\label{Sub:Erdos}
Recall
a well-known result by Erd\"os and Taylor \cite{ErdoesTaylor1960} for planar random walks
with finite variance:  if $\tau$ is the first hitting time of the
origin of a two-dimensional random walk, then
\be{GA2}
   \lim_{t\to\infty}
   {\bf P}^{xt^{\alpha/2}}\big\{\tau>t^{\beta}\big\}
 =
   \frac{\alpha}{\beta}\wedge 1,
\ee
for all $\alpha,\beta\in[0,1]$, and all $x\in\R^2\setminus\{(0,0)\}$ (see, for example,
Proposition 1 in \cite{CoxGri86}). In particular, the right hand side
of (\ref{GA2}) does not depend on $x\in\R^2\setminus\{(0,0)\}$.
Due to this peculiar (specific to $d=2$) property,
the behavior
of the spatial coalescent
started in $\Lambda^{\alpha,t}$
and observed at time $t^\beta$,
asymptotically as $t \to \infty$,
depends  only on the
logarithmic scales $\alpha$ and $\beta$, while all the finer distinctions
are washed out.

For $c\in (0,\infty)$, define
\be{EI}
   I_\alpha(c,t)
 :=
   \big[A^-_\alpha (t), A^+_\alpha(t)\big]
 :=
   \big[(c\log{t})^{-1}\cdot t^{\frac{\alpha}{2}},(c\log{t})\cdot
   t^{\frac{\alpha}{2}}\big].
\ee

Say that a set of locations (marks) $\{x_1,\ldots,x_n\}$
{\em is contained in} $I_\alpha(c,t)$ if and only if
$\parallel x_i-x_j\parallel\in I_\alpha(c,t)$,
for all $1\leq i< j\leq n$.

From (\ref{GA2}) one sees immediately that if
$\{x_1,\ldots,x_n\}$
is contained in $I_\alpha(c,t)$, then for the corresponding random
walks $\{(X_t^j)_{t\geq 0},\,j=1,\ldots,n\}$
with $X_0^j:=x_j$,
\be{Efirstest}
   {\bf P}\big\{X_s^i\neq X_s^j,\;\forall i\neq j,\forall s\in
   [0,g(t)]\big\}
 \tto
   1,
\ee
whenever $g(t)$ is a function satisfying
$g(t)=O(t^{\alpha+\varepsilon})$ for all $\varepsilon>0$.

\subsection{Asymptotic exchangeability}
\label{Sub:asympexch}
In this subsection we perform some
preliminary calculations implying ``asymptotic
exchangeability'' that will be useful in the sequel.
The main result is Proposition~\ref{P:Casyexch} below.

Let $\alpha\in(0,1]$, and set
\be{ag10}
   g_\alpha(t)
 :=
   t^\alpha\log^3 t.
\ee

\begin{remark} In fact, any function $g_\alpha(t)$ with
$t^\alpha\log^2{t}=o(g_\alpha(t))$
%satisfying
%$\lim_{t\to \infty} t^{\alpha/2} \log{t}/\sqrt{g_\alpha(t)} =0$
could be used instead of $t^\alpha\log^3{t}$.
\hfill$\qed$
\end{remark}\sm

For $k\in\N$, let $\zeta$ be a permutation on $\{1,\ldots,k\}$.
Given $\{x_1(t),\ldots,x_k(t)\}\subset\Z^2$, we
denote by $(C^{\alpha,t}_s,L^{\alpha,t}_s)_{s\geq 0}$
the spatial coalescent that starts from
$C^{\alpha,t}_0=\{\{1\},\ldots,\{k\}\}$,
$L^{\alpha,t}_0(\{i\})=x_i(t)$, $i=1,\ldots,k$,
and by $(C^{\alpha,t,\zeta}_s,L^{\alpha,t,\zeta}_s)_{s\geq 0}$
 the spatial coalescent that starts from
$C^{\alpha,t,\zeta}_0=C^{\alpha,t}_0$,
$L^{\alpha,t}_0(\{i\})=x_{\zeta_i}(t)$, $i=1,\ldots,k$.

\begin{proposition}[Asymptotic exchangeability for the spatial coalescent]
\label{P:Casyexch}
Fix $\alpha\in(0,1]$ and $k\in\N$, and
assume that $\{x_1(t),\ldots,x_k(t)\}
\subset\Z^2$ is
contained in $I_\alpha(c,t)$.
If the spatial coalescent $(C^{\alpha,t}_{s},L^{\alpha,t}_s)_{s\geq
  0}$
starts in the marked partition
$\{(\{1\},x_1(t)),\ldots,(\{k\},x_k(t))\}$,
then for all $M\in{\mathcal B}(D([0,\infty),\Pi^{{\mathcal I}}))$,
\be{vl23}
   \lim_{t\to \infty}
   \Big|{\bf P}\big\{(C^{\alpha,t}_s)_{s\geq g_\alpha(t)}\in M\big\}
   -{\bf P}\big\{(C^{\alpha,t,\zeta}_s)_{s\geq g_\alpha(t)}\in M\big\}\Big|
 =
   0.
\ee
\end{proposition}\sm

We prepare the proof by stating the corresponding result for the
underlying random walks.
\begin{lemma}[Asymptotic exchangeability for random walks] Fix
  $\alpha\in(0,1]$, $k\in\N$, and let $\zeta$ be a permutation on
  $\{1,\ldots,k\}$. Let
  $(Y_s)_{s\ge 0}$ be the $2k=k\times 2$ dimensional random walk
\be{EYwalk}
   Y_s
 :=
   \big(X_{1,s}^1,X_{2,s}^1,X_{1,s}^2,X_{2,s}^2,\ldots,X_{1,s}^k,X_{2,s}^k\big),
\ee
where, for each $i\in\{1,\ldots,k\}$,
$(X_s^i)_{s\ge 0}=(X_{1,s}^i,X_{2,s}^i)_{s\geq 0}$
is the two dimensional random walk with transition kernel $a(x,y)$,
and the $k$ random walks are taken to be independent. Moreover, let
\be{EYwalkPi}
   (Y^\zeta_s)_{s\ge 0}
 :=
   \big(X_{1,s}^{\zeta_1},X_{2,s}^{\zeta_1},X_{1,s}^{\zeta_2},
   X_{2,s}^{\zeta_2},\ldots,X_{1,s}^{\zeta_k},X_{2,s}^{\zeta_k}\big).
\ee
\label{Layexch}
Then for all $M\in{\mathcal B}(D([0,\infty),\Z^{2k}))$,
\be{grev27}
   \Big|{\bf P}\big\{(Y_s)_{s\geq g_\alpha(t)}\in M\big\}
   -{\bf P}\big\{(Y^\zeta_s)_{s\geq g_\alpha(t)}\in M\big\}\Big|
 \tto
   0.
\ee
\end{lemma}\sm

\begin{proof}[{\bf Proof}]
The proof relies on
a consequence of the local central limit theorem for
continuous time random walks that we recall next:
if $(Z_s)_{s\ge 0}$ is a random walk in $\Z^d$
(here no moment assumption is needed),
then there exists a finite constant $c$ (see, for example for
our setting \cite{lawlim})
that depends on the dimension and the transition mechanism only,
such that for all $y\in \Z^2$,
\be{Ediffest}
   \sum_{z\in\Z^2}\big|{\bf P}(Z_s=z|Z_0=0)-{\bf P}(Z_s=z|Z_0=y)\big|
 \leq
   \frac{c\parallel y\parallel}{s^{1/2}}.
\ee
We will apply the above difference estimate (\ref{Ediffest})
to $(Y_s)_{s\ge 0}$ and $(Y^\zeta_s)_{s\ge 0}$.
%a particular $2k=k\times 2$ dimensional random walk
%\be{EYwalk}
%Y_s=(X_{1,s}^1,X_{2,s}^1,X_{1,s}^2,X_{2,s}^2,\ldots,
%X_{1,s}^k,X_{2,s}^k),
%\ee
%where, for each $i\in\{1,\ldots,k\}$, $(X_s^i=(X_{1,s}^i,X_{2,s}^i),s\geq 0)$ is the familiar two dimensional random walk with transition kernel $a(x,y)$,
%and the $k$ random walks are taken to be independent.
%Let $M_{1,2}$ denote the event that $X_{\boldsymbol{\cdot}}^1$ and $X_{\boldsymbol{\cdot}}^2$
%simultaneously visit a site of $\Z^{2k}$ before any other pair of random walks above.
%Clearly, $M_{1,2}$
%can also be expressed in terms of $Y$ as the event
%on which $Y$ visits a particular $2k-4$-dimensional subspace before
%${{2k}\choose 2}-1$ other $2k-4$-dimensional subspaces.

Let $M\in{\mathcal B}(D([0,\infty),\Z^{2k}))$.
For each $2k$-tuple $(z^1_1,z^1_2,z^2_1,z^2_2,\ldots,
z^k_1,z^k_2)\in \Z^{2k}$,
set
\be{Eqdef}
   q\big(z^1_1,z^1_2,z^2_1,z^2_2,\ldots,z^k_1,z^k_2\big)
 :=
   {\bf P}\big((Y_s)_{s\geq 0}\in M|Y_0=(z^1_1,z^1_2,z^2_1,z^2_2,\ldots,
   z^k_1,z^k_2)\big).
\ee
%and for each $u\in [0,1)$ define
%\be{EAdef}
%A(u):=\{(z^1_1,z^1_2,z^2_1,z^2_2,\ldots,z^k_1,z^k_2):
%q(z^1_1,z^1_2,z^2_1,z^2_2,\ldots,z^k_1,z^k_2)>u\}.
%\ee

Denote by $B(r)$ the ball in $\R^2$ of radius $r$ centered at $0$.
Suppose
$x^1,\ldots,x^k\in \Z^2\cap B(c \log(t)\, t^{\alpha/2})$,
and let $X^1,\ldots,X^k$ be $k$ independent random walks
with transition kernel $a(\boldsymbol{\cdot},\boldsymbol{\cdot})$
started at locations $x^1,x^2,\ldots,x^k$, respectively.
Let $Y$ be the walk formed as in (\ref{EYwalk}) but using the
walks $X^1,X^2,\ldots,X^k$ as input.
For a
permutation $\zeta$ of $\{1,2,\ldots,k\}$, let $Y^\zeta$ be the walk
formed as in (\ref{EYwalkPi})
%Y_s=(X_{1,s}^{\zeta_1},X_{2,s}^{\zeta_1},X_{1,s}^{\zeta_2},X_{2,s}^{\zeta_2},\ldots,
%X_{1,s}^{\zeta_k},X_{2,s}^{\zeta_k}),
%\ee
%formed
using
$X^{\zeta_1},X^{\zeta_2},\ldots,X^{\zeta_k}$ as input, instead.
Then clearly $Y$ and $Y^\zeta$ have the same transition mechanism,
and the difference $Y_0-Y^\zeta_0$
of their starting locations is a vector with norm bounded by
$O(t^{\alpha/2}\log{t})$.
Therefore, by (\ref{Ediffest}), for all $u\in[0,1]$,
\be{grev26}
%\begin{aligned}
   \Big|{\bf P}\big\{q(Y_{g_\alpha(t)})\ge u\big\} -
   {\bf P}\big\{q(Y^\zeta_{g_\alpha(t)})\ge u\big\}\Big|
 \leq
   \frac{O(t^{\alpha/2}\log{t})}{t^{\alpha/2}\log{t}^{3/2}}
 \tto
   0.
%\end{aligned}
\ee
%which converges to $0$ as $t\to\infty$.
%This is equivalent to saying that
That is,
the $[0,1]$-valued random variables $q(Y_{g_\alpha(t)})$
and $q(Y^\zeta_{g_\alpha(t)})$
are asymptotically equal in distribution.
In particular,
\be{Eimplim}
\begin{aligned}
   &\Big|{\bf E}\big[q(Y_{g_\alpha(t)})\big]-
   {\bf E}\big[q(Y^\zeta_{g_\alpha(t)})\big]\Big|
  \\
 &=
   \Big|{\bf P}\big\{(Y_s)_{s\geq g_\alpha(t)}\in M\big\}-
   {\bf P}\big\{(Y^\zeta_s)_{s\geq g_\alpha(t)}\in M\big\}\Big|
  \\
 &\tto
   0,
\end{aligned}
\ee
and we are done.
\end{proof}\sm

%More importantly we obtain

%\beC{Casyexch}{asymptotic exchangeability for spatial coalescent}
%Fix $ \alpha \in [0,1]$ and
%$k \in \N$ and
%assume that
%$x_1(t), \ldots, x_k(t)$ are
%contained in $I_\alpha(c,t)$.
%If the spatial coalescent $(C^\alpha_{t},t\geq 0)$ starts
%%in the marked partition
%$\{(\{1\}, \ldots, \{k\}, (x_1 (t), \ldots, x_k(t)) \}$, and if
%$Y^{c}$, $Y^{c,\zeta}$ (``c'' stands for coalescent)
%are defined in (\ref{EYwalk}), (\ref{EYwalkPi})
%with input processes
%$X_s^i := L_s^{\{i\sim^s\}}$, $s\geq0$, $i\in \{1,\ldots,k\}$, then
%\be{grev28}
%\lim_{t\to \infty}
%|{\bf P}((Y^{c}_s)_{s\geq g_\alpha(t)} \in M) - {\bf P}((Y^{c,\zeta}_s)_{s\geq g_\alpha(t)}\in M )|
%= 0.
%\ee
%\end{corollary}

\begin{proof}[{\bf Proof of Proposition~\ref{P:Casyexch}}]
Let
the $2k$-dimensional processes $Y^{sc}$ and $Y^{sc,\zeta}$
(``sc'' stands for semi-coalescent)
be formed as in (\ref{EYwalk}) and (\ref{EYwalkPi}),
however the input random processes $X^1,\ldots,X^k$
are changed so that
$X^i$s are independent continuous-time random walks
with kernel $a(\boldsymbol{\cdot},\boldsymbol{\cdot})$
until time $g_\alpha(t)$,
and after time $g_\alpha(t)$ their joint evolution is
the evolution of the location process
of the spatial coalescent with initial
configuration $(X_{g_\alpha(t)}^1,\ldots,
X_{g_\alpha(t)}^k)$. Moreover, let $Y^{c}$ and $Y^{c,\zeta}$
(``c'' stands for coalescent) be the
$2k$-dimensional processes whose joint evolution is
the evolution of the location process
of the spatial coalescent with initial
configuration $(X_{g_\alpha(t)}^1,\ldots,
X_{g_\alpha(t)}^k)$.

It is obvious how to construct couplings
$(Y^{sc}, Y^c)$ and $(Y^{sc,\zeta},Y^{c,\zeta})$, so
 that on the event $\{$no coalescence up to time $g_\alpha(t)\}$
%=$\{\tau_{N-1}^\alpha(t)>g_\alpha(t)\} $
the two processes, the coalescent and the corresponding
semi-coalescent,
in both couplings above agree for all times.
Hence,
\be{grev29}
\begin{aligned}
   \Big|{\bf P}\big\{&(Y^{c}_s)_{s\geq g_\alpha(t)}\in M\big\} -
    {\bf P}\big\{(Y^{\zeta}_s)_{s\geq g_\alpha(t)}\in M\big\}\Big|
  \\
 &\le
   2{\bf P}\big\{\mbox{coalescence occurs before time
    }g_\alpha(t)\big\}
  \\
 &\quad+
   \Big|{\bf P}\big\{(Y^{sc}_s)_{s\geq g_\alpha(t)}\in M\big\} -
   {\bf P}\big\{(Y^{sc,\zeta}_s)_{s\geq g_\alpha(t)} \in M\big\}\Big|.
\end{aligned}
\ee

The claim follows immediately from the previous observations
and from the fact
\be{grev30}
   {\bf P}\big\{\mbox{coalescence occurs before time
   }g_\alpha(t)\big\}
 \tto
   0,
\ee
which is a direct consequence of (\ref{Efirstest}).
\end{proof}

\subsection{Monotonicity and consequences}
\label{Sub:monoton}
Recall the set of
marked partitions $\Pi^{{\mathcal I},\Z^2}$ from (\ref{vl10}).
It is convenient to introduce a partial order ``$\leq$'' on
$\Pi^{\mathcal I,\Z^2}$. Let for $\CP_1,\CP_2 \in \Pi^{\mathcal I,\Z^2}$,
\be{Epartord}
   \CP_1\leq \CP_2
\ee
iff for each $g\in \Z^2$
the number of partition elements in $\CP_1$ with
mark $g$ is bounded above by
the number of partition elements in $\CP_2$ with
mark $g$.
For brevity reasons,
we will often omit from (\ref{Epartord}) the dependence
on the location processes when evident from the context,
so we will write
\be{vl15}
   C_1\leq C_2
\ee
to mean $(C_1,L_1)\leq (C_2,L_2)$.

\begin{remark}\label{Rem:00}\rm
Note that if $\CP_1\leq \CP_2$,
%due to the nature of migration
%and coalescence mechanisms,
one can easily construct a coupling
$((C_s^1,L_s^1),(C_s^2,L_s^2))_{s\geq 0}$ of the spatial coalescents where
$(C^j_0,L^j_0)={\mathcal P}_j$,
$j=1,2$,
such that $C_s^1\leq C_s^2$, for all $s\geq 0$, almost surely.
\hfill$\qed$
\end{remark}\sm

%and its partial order from (\ref{Epartord}).
Suppose that $f:\Pi^{{\mathcal I},\Z^2}\to\R$ is
non-decreasing,
%a function
%$f:{\mathcal P}\to \R$ which is {\em non-decreasing} in the sense that
%\be{vl27}
%\CP_1 \leq \CP_2 \Rightarrow f(\CP_1) \leq f(\CP_2)
%\ee
and let $g:[0,\infty)\to(0,\infty)$.
For $a,b\in [-\infty,\infty]$ consider asymptotic behavior(s)
of the type
\begin{equation}
\label{Easybeh}
\limsup(\liminf)_{t\to\infty}\frac{f(C_t,L_t)}{g(t)} =a, \, \
\limsup(\liminf)_{t\to\infty}\frac{E(f(C_t,L_t))}{g(t)} = b.
\end{equation}

An important
observation is the next easy consequence of monotonicity and
Remark~\ref{Rem:00}. Namely, if any of the four types of
asymptotic behavior (\ref{Easybeh}) holds for
both spatial coalescents $(C_t^j, t\geq 0)$, $j=1,3$, and if
\be{vl28}
   {\bf P}\big\{C_0^1 \leq C_0^2 \leq C_0^3\big\} =1,
\ee
then the same asymptotic behavior holds for the spatial coalescent $(C_t^2,
t\geq 0)$.

Moreover, let $\CA\subseteq \R$, and suppose we are given
three coalescent families
\be{vl29}
  \big\{(C_s^{j,\alpha})_{s\geq 0};\,\alpha \in \CA,j\in\{1,2,3\}\big\},
\ee
with initial states such that
\be{Esandwich}
   {\bf P}\big\{C_0^{1,\alpha}\leq C_0^{2,\alpha}\leq C_0^{3,\alpha},\,
   \forall\alpha\in\CA\big\}
 =
   1.
\ee

In addition, assume that %there exists a process $(C_s)_{s\ge 0}$ with
c\`adl\`ag path such that
\be{Esandwich1}
   \lim_{\alpha\to\alpha_0}(C_s^{1,\alpha})_{s\geq 0}
 =
   \lim_{\alpha\to\alpha_0}(C_s^{3,\alpha})_{s\geq 0},
\end{equation}
where the above convergence is weak convergence on $D([0,\infty),\Pi^{{\mathcal
    I}})$equipped with the Skorokhod topology.

\begin{lemma}
\label{Lsandwich}
If (\ref{Esandwich}) and (\ref{Esandwich1}) hold, then
$(C_s^{2,\alpha})_{s\geq 0}$ also converges in law as $\alpha \to \alpha_0$,
and
\be{grev31}
   \lim_{\alpha\to\alpha_0}(C_s^{2,\alpha})_{s\geq 0}
 =
   \lim_{\alpha\to\alpha_0}(C_s^{1,\alpha})_{s\geq 0}.
\ee
\end{lemma}

The next result will, together with the above consequences of monotonicity,
eventually be used
for deducing various asymptotics for the spatial
coalescent started from
infinite configurations, given the results for the spatial coalescents started
from finite configurations.

Let $(K_s)_{s\geq 0}$ be  the Kingman coalescent.
\begin{lemma}
\label{LPoisson}
For each $\delta>0$ there exists $\rho=\rho(\delta)\in (0,\infty)$
such that
\be{grev6c}
   {\bf P}\big\{\# K_\delta\geq n\big\}
 \leq
   {\bf P}\big\{X_\rho \geq n\big\},
\ee
where $X_\rho\stackrel{\mathrm{d}}{=}1+\mathrm{Poisson}(\rho)$.
That is,
${\bf P}\{X_{\rho}=k\}=e^{-\rho}\rho^{(k-1)}/(k-1)!$, for all $k\geq 1$.
\end{lemma}\sm

\begin{remark}\label{Rem:08}\rm
The shift by one unit is necessary here
since ${\bf P}(K_\delta \geq 1) =1$.
\hfill$\qed$
\end{remark}\sm

\begin{proof}[{\bf Proof}]
Let $\{\Upsilon_n;\,n\geq 1\}$ be the family of independent
exponential random variables where $\Upsilon_n$ has rate $n(n+1)/2$.
Then by construction of Kingman's coalescent (see, for example, \cite{Kingman82,Aldsurvey}),
\be{EtailK}
%\begin{aligned}
   {\bf P}\{\# K_\delta>n\}
 =
   {\bf P}\{\sum_{k\geq n}\Upsilon_k > \delta\}
%  \\
% &=
%   {\bf P}\{e^{\theta\sum_{k\geq n}\Upsilon_k}> e^{\theta\delta}\}
%  \\
 \leq
   e^{-\theta\delta}{\bf E}\big[e^{\theta\sum_{k\geq n}\Upsilon_k}\big],
%\end{aligned}
\ee
for all $\theta\in\R$. Assume that $\theta<n(n+1)/2$, and consequently
that ${\bf E}\big[e^{\theta\sum_{k\geq
 n}\Upsilon_k}\big]<\infty$.

Since
\be{EtailK2}
\begin{aligned}
   {\bf E}\big[e^{\theta \sum_{k\geq n}\Upsilon_k}\big]
 &=
   \prod_{k=n}^\infty\frac{\frac{(k+1)k}{2}}{\frac{(k+1)k}{2}-\theta}
  \\
 &=
   \exp \left [\sum_{k=n}^\infty\log\big(
   1+\tfrac{\theta}{\frac{(k+1)k}{2}-\theta}\big)\right]
  \\
 &\leq
   \exp \left [\sum_{k=n}^\infty
   \left(\frac{\theta}{\frac{(k+1)k}{2}-\theta}+O(\frac{\theta^2}{(\frac{(k+1)k}{2}-\theta)^2}\big)
   \right)\right],
\end{aligned}
\ee
by (\ref{EtailK})
\be{Etail3}
   {\bf P}\big\{\# K_\delta>n\big\}
 \leq
   \exp\big[-\delta\theta+\sum_{k=n}^\infty\frac{\theta}{\frac{(k+1)k}{2}-
   \theta}\big].
\ee

Plugging in, for example, $\theta=n\log^2{n}$ gives
\be{EtailK1}
   {\bf P}\big\{\# K_\delta>n\big\}
 =
   O(e^{-\delta n (\log{n})^2/2}),
\ee
which is of a smaller order than
\be{grev6c2}
   {\bf P}\big\{\mathrm{Poisson}(\rho)+1>n\big\}
 \asymp
   C(\rho)e^{-n (\log{n}+O(1))},
\ee
for all large $n$, where $O(1)$ indicates a term that stays bounded as $n \to \infty$.
Since the sum of independent Poisson random
variables is another Poisson random variable, we
can choose $\rho$ appropriately large
so that
${\bf P}\big\{\# K_\delta>n\big\}\leq {\bf
  P}\big\{\mathrm{Poisson}(\rho)>n-1\big\}$, for
all $n\geq 1$.
\end{proof}

\section{Asymptotics for sparse particles}
\label{S:ASP}
Fix throughout this section $\alpha\in(0,1]$.
Our goal in this section is to analyze the behavior of a finite coalescent
with particles spaced at distance $t^{\alpha/2}$ and observed
at time $t^\beta$, $\beta>\alpha$, as $t \to \infty$.

Recall the instantaneous coalescent that
corresponds to the spatial coalescent with
resampling rate $\gamma=\infty$.
In our setting $\gamma\in (0,\infty)$ is fixed.
Nevertheless, we still can rely on
the ``loss of the spatial
structure'' property of the coalescent on time scales $t^\beta$
for the instantaneous coalescent with
partition elements situated  initially at mutual distances of order
$t^{\alpha/2}$ that was exploited in \cite{CoxGri86}.

%For fixed $\alpha \in (0,1]$ and
% finite (large) t consider the torus
%\be{vl29a}
%[-t^{\alpha/2}, t^{\alpha/2}]^2 \cap \Z^2.
%\ee
%Of course, $t^{\alpha/2}$ may not be integral but we will use it
%in limits of sums regardlessly.

Recall $\Lambda^{\alpha,t}$ from (\ref{rect}).
We denote by
\be{vl29a2}
   (C^\alpha_{s})_{s\geq 0},\qquad \mbox{ and }\;(IC^\alpha_{s})_{s\geq 0}
\ee
the spatial coalescent and the instantaneous coalescent
starting from initial configuration
$C_0^\alpha$ with marks contained in $\Lambda^{\alpha,t}$.
Notice that $t$ is suppressed from the notation,
but this should not cause confusion.
% and the symbols $\alpha$ and $t$ are
%reserved for the purpose
%of parametrizing the size of the torus and in this way
%all the spatial coalescents we consider.
%We denote by
%\be{vl29b}
%(IC^\alpha_{s},s\geq 0)
%\ee
%any analogous instantaneous coalescent.

There are classical results on $(IC^\alpha_{s})_{s\geq 0}$
with initially $N$ individuals spread out in $\Lambda^{\alpha,t}$, and
observed at time $t^\beta$, where $\beta>\alpha$, which we wish to
recall first.
%This has been investigated by Bramson, Cox and
%Griffeath \cite{BramCoxGri86} for the case of \emph{instantaneous} coalescence.
Let $c>0$ and recall $I_\alpha(c,t)$ from (\ref{EI}).
The following result was proved in a beautiful paper by
Cox and Griffeath \cite{CoxGri86}
under the additional assumption that the underlying random walks are simple
random walks: for fixed $N \in \N$,  the
initial
locations
$\{x_1(t),\ldots,x_N(t)\}$ contained in
$I_\alpha(c,t)$ and for each $\beta>\alpha$,
\be{grev8}
   \CL\big[\# IC^\alpha_{t^\beta}\big]
 \Tto
   \CL\big[\# K_{\log\frac{\beta}{\alpha}}\big].
\ee

We next consider the spatial (delayed) coalescent,
and show the stronger form of weak convergence
in two ways:
(i) in the sense of path-valued random variables
where $\beta$ is the ``time''-parameter, and (ii)
accounting for the partition structure.
Note that the weak convergence is done in the sense of the
discrete topology.
\begin{proposition}[Finite sparse coalescents: large time
    scales]
\label{P3.1}
Fix $N\in\N$, $\alpha\in(0,1]$ and $c>0$, and assume that
$\{x_1(t),\ldots,x_N(t)\}\subset\Z^2$ is
contained in $I_\alpha(c,t)$.
Let the spatial coalescent $(C^\alpha_{s})_{s\geq 0}$ start in
$\{(\{1\},x_1(t)),\ldots,(\{N\},x_N(t))\}$.
Then
\be{GR15x}
   \CL\big[(C^\alpha_{t^\beta})_{\beta \in [\alpha,\infty)}\big]
 \Tto
   \CL\big[(K^N_{\log (\frac{\beta}{\alpha})})_{\beta
     \in[\alpha,\infty)}\big],
\ee
where $(K^N_t)_{t \geq 0}$ is the Kingman coalescent started in
$\{\{1\},\ldots,\{N\}\}$.
\end{proposition}\sm

The proof of this result is given in the next two subsections.

\subsection{Convergence of marginal distributions}
\label{Sub:convfdd}
A key element of the proof is the following fact which we state
for future reference.
\begin{lemma}[Lemma~1 from \cite{CoxGri86}]
\label{L1CoxGri}
Fix $\alpha_0>0$, and $c>0$.
Let $\{(X_s^i)_{s\geq 0};\,i=1,\ldots,4\}$ be a family of
independent  random walks with $X_0^i=x_i$, for $i=1,\ldots,4$.
Then uniformly in $\alpha\in [\alpha_0,\infty)$ and
$\{x_1,\ldots,x_4\}\subset\Z^2$ contained in $I_\alpha(c,t)$, we have
\be{Eunifpos}
   \int_{t^\alpha}^\infty\mathrm{d}s\,
   {\bf P}\big(\{X^1_s = X^2_s \}\cap
   \{X^1_s, X^3_s,X^4_s \mbox{ not contained in }
   I_1(4c,s)\}\big)
 \tto
   0.
\ee
\end{lemma}\sm

\begin{remark}\label{Rem:09}\rm
In the setting of \cite{CoxGri86} the walks are simple
symmetric walks, but the proof of the corresponding lemma is more general,
depending solely on the uniform bound
\be{vl29d}
  {\bf P}\{X_s = y\}
 \leq
   \frac{c}{s},
\ee
for all $y\in \Z^2$, and $s \geq 0$. It therefore applies to our
setting (see \cite{spitzer}).
\hfill$\qed$\end{remark}\sm

The first major step to prove Proposition~\ref{P3.1} is to show:

\begin{lemma}[Finite sparse coalescents: convergence of marginals]
\label{L3.1mar}
Fix $N\in\N$, $\alpha\in(0,1]$ and $c>0$, and assume that
$\{x_1(t),\ldots,x_N(t)\}\subset\Z^2$ is
contained in $I_\alpha(c,t)$.
Let the spatial coalescent $(C^\alpha_{s})_{s\geq 0}$ start in
$\{(\{1\},x_1(t)),\ldots,(\{N\},x_N(t))\}$.
Then for all $\beta>\alpha$,
\be{GR15}
   \CL\big[\# C^\alpha_{t^\beta}\big]
 \Tto
   \CL\big[\# K^N_{\log{(\frac{\beta}{\alpha})}}\big],
\ee
where $(K^N_s)_{s\geq 0}$ is the Kingman coalescent started in
$\{\{1\},\ldots,\{N\}\}$.
\end{lemma}\sm

\begin{proof}[{\bf Proof}]
The argument makes use of an obvious coupling of
$(C^\alpha_{\cdot},L^\alpha_{\boldsymbol{\cdot}})$ and
$(IC^\alpha_{\cdot},IL^\alpha_{\boldsymbol{\cdot}})$ where
$IC^\alpha_{0}:=C^\alpha_{0}$. We proceed by induction on $N\in\N$.

We start with $N=2$. Put
\be{tauregul}
   \tau_1'(t)
 :=
   \inf\big\{s>0:\,\# IC^\alpha_s=1\big\},
\ee
and set $C^\alpha_s:=IC^\alpha_s$, for all $s\in[0,\tau_1'(t)]$.
Define $(C^\alpha_s)_{s>\tau^\prime_1(t)}$ in a standard way,
using additional (independent)
randomness. Let then
\be{tau}
   \tau_1(t)
 :=
   \tau_1^{\alpha,t}
 :=
   \inf\big\{s>0:\,\# C^\alpha_{s}=1\big\},
\ee
so that $\tau_1'(t)$ and $\tau_1(t)$
are the coalescence times of the two particles in $IC^\alpha$, and
$C^\alpha$, respectively.
Then clearly
\be{Edomin}
   \tau'_1(t)
 \leq
   \tau_1(t)
 \leq
   \tau_1'(t)+\sum_{i=0}^G \tau_i^0
\ee
where $G$ has shifted geometric distribution with
success probability $\gamma/(2+\gamma)$, i.e.,
${\bf P}\{G\ge m\}=(2/(2+\gamma))^{m-1}$, for all $m\geq 1$,
$\tau_i^0$, $i\ge 1$,  is
distributed as the length of the (almost surely finite) excursion
away from $0$ for the underlying migration walk, and where the family
$\{\tau_1'(t),G,\{\tau_i^0,i\geq 0\}\}$ is an independent family of
random variables.

The result of Cox and Griffeath discussed above is based on
the Erd\"os-Taylor asymptotics (\ref{GA2}) and stronger estimates
of a similar type.
In particular,
we rewrite (\ref{GA2}) in the current setting,
where $\beta>\alpha$ and the random walk is twice as fast as the
simple one, as
\be{ErdTay}
   {\bf P}\big\{\tau_1'(t)>t^{\beta}/2\big\}
 \tto
   \frac{\alpha}{\beta}.
\ee

Note that
(\ref{ErdTay}) can be restated as
the following convergence in distribution: for all $u\ge 0$,
\be{Econexp}
   \lim_{t\to\infty}
   {\bf P}\Big\{\log\big(\frac{\log\tau_1'(t)}{\alpha\log t}\big)<u\Big\}
 =
   1-e^{-u}.
\ee

We would like to show the same convergence holds with $\tau_1(t)$ in
place of $\tau_1'(t)$.
Due to $(\ref{Edomin})$ it suffices to show that, as $t \to \infty$,
with overwhelming probability,
\be{Edomin1}
\sum_{i=0}^G \tau_i^0  \leq \tau_1'(t),
\ee
since then $\log(\tau_1'(t)+\sum_{i=0}^G \tau_i^0)\leq \log{\tau_1'(t)} +
\log 2$, and $\log 2/\log t$ becomes negligible in the limit.
Since
%$G$ is geometric with positive success probability
%we have that
$\sum_{i=0}^G\tau_i^0<\infty$, almost surely,
%is a finite random variable while
and $\tau_1'(t)\to \infty$, as $t\to \infty$, in
probability,
%being of the order
%\be{vl29e}
%t^{\alpha \exp{U}}, \mbox{ where } \CL (U) = \mbox{ Exponential (rate $1$)},
%\ee
(\ref{Edomin1}) trivially follows, and we have
%by (\ref{Edomin}) and (\ref{Edomin1})
\be{Econexp1}
   \lim_{t\to\infty}
   {\bf P}\Big\{\log\big(\frac{\log\tau_1(t)}{\alpha\log t}\big)<u\Big\}
 =
   1-e^{-u},
\ee
for all $u\ge 0$.

Now note that for $N>2$
and for $\beta\geq \alpha$, using analogous coupling
of $(C^\alpha_{\boldsymbol{\cdot}},L^\alpha_{\boldsymbol{\cdot}})$ and
$(IC^\alpha_{\boldsymbol{\cdot}},IL^\alpha_{\boldsymbol{\cdot}})$ up to the first coalescence
time $\tau_{N-1}'(t)$ in $IC^\alpha$,
\be{Enocoal}
   \lim_{t\to \infty}{\bf P}\big\{\# C^\alpha_{t^\beta} = N\big\}
 =
   \lim_{t\to \infty}
   {\bf P}\big\{\# IC^\alpha_{t^\beta}=N\big\}
 =
   \Big(\frac{\alpha}{\beta}\Big)^{N \choose 2}.
\ee
where the second limit above was evaluated in Proposition~2
of \cite{CoxGri86}.
Moreover, if
\be{vl29f}
   \tau_{N-1}(t)
 :=
   \inf\big\{s>0:\, \# C^\alpha_{s}=N-1\big\},
\ee
due to the fact that
$|\log\tau_{N-1}'(t)-\log\tau_{N-1}(t)|\to 0$, as $t\to\infty$,
almost surely (argue as for (\ref{Edomin1}) above),
the induction step in the proof of  \cite{CoxGri86} Theorem~3
can be carried out verbatim.
The details are tedious, so we omit them, and state instead that
\be{vl30}
   p_{N,k}(\alpha/\beta)
 :=
   \lim_{t\to \infty}{\bf P}\big\{\# C^\alpha_{t^\beta}=k\big\}
\ee
satisfies the recursion of \cite{CoxGri86} Theorem~3,
\be{ErecCG}
   p_{N+1,k}(\frac{1}{s})
 =
   {{N+1} \choose 2}\int_1^s\mathrm{d}y\,y^{-{{N+1} \choose 2}-1}
   p_{N,k}(y/s),
\ee
for all $s\ge 1$ and $1\le k\le N+1$.

Since the initial conditions (\ref{Econexp1}), (\ref{Enocoal})
to the recursion are identical to those in Theorem~3 in
\cite{CoxGri86},
as argued above, the solution is the same, and so we
have verified that
for each $\beta>\alpha$ and each $k\in\{1,\ldots,N\}$,
\be{vl31}
   \lim_{t\to \infty}{\bf P}\big\{\# C^\alpha_{t^\beta}= k\big\}
 =
   \lim_{t\to \infty}{\bf P}\big\{\# IC^\alpha_{t^\beta}=k\big\}
 =
   {\bf P}\big\{\# K^N_{\log (\frac{\beta}{\alpha})}=k\big\},
\ee
where the last identity was again obtained in \cite{CoxGri86}.
\end{proof}\sm

\subsection{Convergence in path space}
\label{Sub:convps}
In order to show path convergence
of $(\# C^\alpha_{t^\beta})_{\beta\geq \alpha}$
to
$(\# K_{\log{\beta/\alpha}})_{\beta\geq \alpha}$
one defines a sequence of random times
$\{\tau^{\alpha}_k(t);\,1\leq k\leq N\}$,
where for each $k\ge 1$,
\be{Elap}
   \tau^{\alpha}_k(t)
 :=
   \inf\big\{s\ge 0:\;\#C^{\alpha}_s\le k\big\},
\ee
where as usual $\tau^{\alpha}_k(t)=\infty$ if $\inf_{s\ge
  0}\#C^{\alpha}_s>k$. That is,
$\tau^\alpha_N=0$, and
$\tau^{\alpha}_{N-1}(t)$ is the first coalescence time,
(also denoted by $\tau_{N-1}(t)$ in the proof of Lemma \ref{L3.1mar}),
$\tau^{\alpha}_{N-2}(t)$ is the second coalescence time,
etc.
It is not difficult to see that the arguments of the proof of
Theorem~3 in \cite{CoxGri86}
extend to showing that, with  probability one
$\# C_{\tau_k^\alpha(t)}=k$, for all $k=N-1,\ldots,1$
(see also Lemma \ref{L1CoxGri}),
and that with respect to convergence in probability,
\be{Eorderofmag}
   \lim_{t\to\infty}\frac{\tau_k^\alpha(t)}{\tau_{k-1}^\alpha(t)}
 =
   0,
\ee
for each $k\geq 2$. (Note here that the remaining $k$ partition elements
are spread out). Moreover, the
following joint convergence in distribution holds
\be{Ejointcon}
\begin{aligned}
   \Big(&\log\big(\tfrac{\log(\tau^{\alpha}_{N-1}(t))}{\alpha\log t}\big),
   \log\big(\tfrac{\log(\tau^{\alpha}_{N-2}(t)- \tau^{\alpha}_{N-1}(t))}
   {\log(\tau^{\alpha}_{N-1}(t))}\big),\ldots,
   \log\big(\tfrac{\log(\tau^{\alpha}_1(t)-\tau^{\alpha}_2(t))}
   {\log(\tau^{\alpha}_2(t))}\big)\Big)
  \\
 \Tto&
   \big(U_{N-1},U_{N-2},..,U_1\big),
\end{aligned}
\ee
where $\{U_i;\,i=1,\ldots,N-1\}$ is a family of independent random
variables such that for all $i=1,\ldots,N-1$,
$U_i$ has the rate ${i+1 \choose 2}$ exponential distribution.
Now (\ref{Eorderofmag}) and (\ref{Ejointcon})
imply the convergence of random vectors
\be{Ejointcon1}
\begin{aligned}
   \Big(&\log\big(\tfrac{\log(\tau^{\alpha}_{N-1}(t))}{\alpha\log t}\big),
   \log\big(\tfrac{\log(\tau^{\alpha}_{N-2}(t))}{\alpha\log t}\big)-
   \log\big(\tfrac{\log(\tau^{\alpha}_{N-1}(t))}{\alpha\log t}\big),\ldots,
   \log\big(\tfrac{\log(\tau^{\alpha}_1(t))}{\alpha\log t}\big)-
   \log\big(\tfrac{\log(\tau^{\alpha}_2(t))}{\alpha\log t}\big)
   \Big)
  \\
 \Tto&
   \big(U_{N-1},U_{N-2},\ldots,U_1\big).
\end{aligned}
\ee

Since
\be{ECnice}
   \# C^\alpha_{t^\beta}
 =
   1_{\big[\frac{\log(\tau^\alpha_1(t))}{\alpha\log t},\infty\big)}
   \big(\tfrac{\beta}{\alpha}\big)
 +
   \sum_{k=2}^N k 1_{\big[\frac{\log(\tau^\alpha_k(t))}{\alpha\log t},
   \frac{\log(\tau^\alpha_{k-1}(t))}{\alpha\log t}\big)}
   \big(\tfrac{\beta}{\alpha}\big)
\ee
and with $\bar U_k = exp (U_N+ \cdots+U_k)$,
\be{EKnice}
   \# K^N_{\log(\beta/\alpha)}
 =
   1_{\big[\bar U_1,\infty\big)}\big(\tfrac{\beta}{\alpha}\big)
 +
   \sum_{k=2}^N k 1_{\big[\bar U_k,\bar U_{k-1}\big)}\big(\tfrac{\beta}{\alpha}\big),
\ee
it immediately follows that the process
$(\# C^\alpha_{t^\beta})_{\beta\geq \alpha}$
converges in the sense of Skorokhod topology to the process
$(\# K^N_{\log\frac{\beta}{\alpha}})_{\beta\geq \alpha}$, as $t\to\infty$.

In order to upgrade the above convergence
to the one on the level of partitions, as stated in Proposition~\ref{P3.1},
we need to make sure that for any fixed $N$ and any choice of initial locations
$x_1,\ldots,x_N$ contained in $I_\alpha(c,t)$, asymptotically as $t\to
\infty$, any two current partitions  elements coalesce equally likely and independently of
the coalescent time. That is,
\be{Eunicolap}
   {\bf P}\big(i,j\mbox{ coalesce at time }\tau^\alpha_{N-1}(t)|
   \tau^\alpha_{N-1}(t)\big)
 \tto
   {N \choose 2}^{-1}.
\ee

Assume without loss of generality that $i<j$.
Fix $\beta >\alpha$, and let (with $C^\alpha_0:=\{\{1\},\ldots,\{N\}\}$
and $g_\alpha$ as in (\ref{ag10}))
\be{EMbetaij}
%\begin{aligned}
   M^\beta_{i,j}(t)
 :=
   \bigcup_{s\in[0,t^\beta-g_\alpha(t)]}\big\{
   C^\alpha_{s-}=C^\alpha_0,\,C^\alpha_s=\{i,j\}\cup C^\alpha_{s-}\setminus\{\{i\},\{j\}\}
   \big\},
\ee
and put
\be{ENbeta}
   N^\beta(t)
 :=
   \bigcup_{1\leq i < j \leq N}M^\beta_{i,j}(t).
\ee
Note that the events $\{M_{i,j}^\beta(t);\,1\le i<j<\infty\}$
are disjoint.

Recall from the proof of
Proposition~\ref{P:Casyexch} that $(Y_s^c)_{s\ge 0}$ denotes the
$2N$-dimensional process (i.e. $\Z^{2N}$-valued), whose joint evolution is the evolution of
the location processes of $(C^\alpha_s,L_s^\alpha)_{s\ge 0}$ but started at time 0 in
$(X^1_{g_\alpha (t)}, \cdots, X^N_{g_\alpha (t)})$ where the latter
are $N$-independent $a(\cdot,\cdot)$-random walks on $\Z^2$.
We consider the path of $Y^c$ after time $g_\alpha(t)$ up to time $t$
and ask whether the coalescent with these paths in the time interval
$[0,t^\beta-g_\alpha(t)]$ would have a first coalescence event, we write
$(Y^c_s)_{s \geq g_\alpha(t)} \in M^\beta(t)$ for this event.

Then
\be{Easym1}
   \Big|{\bf P}\big\{(Y^{c}_s)_{s\geq g_\alpha(t)}\in N^\beta(t)\big\}
   -{\bf P}\big\{\tau_{N-1}^\alpha(t)\leq t^\beta\big\}\Big|
 \leq
   {\bf P}\big\{\tau_{N-1}^\alpha(t) \leq g_\alpha(t)\big\},
\ee
and similarly for each $i<j$,
\be{Easym2}
   \Big|{\bf P}\big\{(Y^{c}_s)_{s\geq g_\alpha(t)}\in
   M^\beta_{i,j}(t)\big\}
   -{\bf P}\big\{\tau_{N-1}^\alpha(t)\leq t^\beta,
   i\sim^{\tau_{N-1}^\alpha(t)} j\big\}\Big|
 \leq
   {\bf P}\big\{\tau_{N-1}^\alpha(t)\leq g_\alpha(t)\big\}.
\ee

Proposition~\ref{P:Casyexch} together with Lemma \ref{L3.1mar} and
(\ref{Easym2}) imply
\be{vl29g}
   \Big|{\bf P}\big\{\tau_{N-1}^\alpha(t)\leq t^\beta,
   1\sim^{\tau_{N-1}^\alpha(t)} 2\big\}-
   {\bf P}\big\{\tau_{N-1}^\alpha(t) \leq t^\beta,
   i\sim^{\tau_{N-1}^\alpha(t)} j\big\}\Big|
 \tto 0,
\ee
and again due to (\ref{EMbetaij}), (\ref{ENbeta}),
and  (\ref{Easym1}),
\be{vl29h}
   \Big|{\bf P}\big\{\tau_{N-1}^\alpha(t) \leq t^\beta,
   1\sim^{\tau_{N-1}^\alpha(t)} 2\big\}
   -{{N \choose 2}}^{-1}{\bf P}\big\{\tau_{N-1}^\alpha(t)
   \leq t^\beta\big\}\Big|
 \tto
   0,
\ee
which proves (\ref{Eunicolap}).

Due to the asymptotic exchangeability given by
Proposition~\ref{P:Casyexch} and uniform estimates
(\ref{Eunifpos})
on locations of
partition elements at each coalescence time, it is easy to extend
(for example by induction)
(\ref{Eunicolap}) to an analogous statement
at any future
coalescence time.
This indeed confirms that the limiting object $K^N$ is
the Kingman coalescent, since the right hand sides of
(\ref{Ejointcon1}) and (\ref{Eunicolap})
characterizes its law completely.

\section{Asymptotics for dense particles at small times}
\label{S:ADPST}

This section concentrates on the behavior of the system
 for fixed $\alpha\in[0,\infty)$
at times of order only slightly larger than the area of the
rectangle on which the initial configuration is supported.
More precisely, we set

\be{e:015}
   \Lambda(r)
 :=
   [-r,r]^2\cap\Z^2.
\ee
and study the corresponding restricted spatial coalescent.

\subsection{Coupled spatial coalescents and moment bound}
\label{Sub:CSC}

Here and at many other occasions it is useful to couple coalescents
starting in different but comparable initial configurations.
We next describe a formal setting that will be used in
Sections \ref{S:ADPST}, \ref{S:ltsas} and~\ref{S:ConvSpatScale}.

Let
\be{Fz}
   F:=\big\{F_z,\,z\in\Z^2\big\}
\ee
be a family of
$\bar{N}$-valued valued random variables.
We think of $F_{z}$ as the number of partition elements (particles)
present at site $z\in\Z^2$ in the coalescent at time $0$.
In symbols,
\be{vl16}
   F_z
 :=
   \#\big\{\pi\in C_0:\,L_0(\pi)=z\big\}.
\ee

Typically we will choose the collection $F$ such that
$\sum_{z\in\Z^2}F_z\delta_z\in{\mathcal E}$, almost surely.
In addition, for the applications we have, we often assume $F$ to be a
family of independent random variables with the same Poisson
(rate $\rho\in (0,\infty)$) distribution. \sm

Fix a countable (or finite) set ${\mathcal I}$, and recall for all
${\mathcal I}'\subseteq{\mathcal I}$ satisfying (\ref{Lip1}) the
restricted process $(C_s^{{\mathcal I}'},L_s^{{\mathcal I}'})_{s\ge 0}$.

Sometimes we are interested in restricting $(C_s,L_s)_{s\ge 0}$ to
geographical information. That is, for
$A\subseteq\Z^2$, let
\be{IA}
  \CI_A
 :=
   \big\{i\in{\mathcal I}:\,L_0(\{j:\,j\sim_{C_0}i\})\in A\big\}.
\ee
In this particular case, we write
\be{CA}
   C(A)
 :=
  C^{\CI_A}.
\ee

In particular, if $F^A$ gives the number of partition elements of the
restricted coalescent $C(A)$, then
\be{vl18}
   F_z^A
 :=
   \left\{\begin{array}{cl}
   F_z,& z\in A
  \\
   0, & z\not \in A.\end{array}\right.
\ee

Moreover,
if $A\subset B\subset \Z^2$
then $C_0^{\CI_A}\leq C_0^{\CI_B}$ and due to the comment following
(\ref{Epartord}), the two coalescents
$(C^{\CI_A}, L^{\CI_A})$
and
$(C^{\CI_B}, L^{\CI_B})$
can be coupled so that at any point in time and space,
the number of partition elements in
$(C^{\CI_B}, L^{\CI_B})$
dominates from above the number of partition elements in
$(C^{\CI_A}, L^{\CI_A})$.
\sm

Assume we are given the coupled spatial coalescents from
above and recall
$\{F_z;\,z\in\Z^2\}$ from (\ref{Fz}). Assume that
\be{Econdcoup}
  {\bf E}\big[F_z\big]>0,\quad\mbox{and }{\bf Var}\big[F_z\big]<\infty,
\ee
for all $z\in\Z^2$.

Our goal is to show next
that the sparse initial configurations necessary for the results of
the previous section arise if the coalescent is started in the torus
$\Lambda^{\alpha,t}=\Lambda(t^{\alpha/2})$, and
observed at time $t^{\alpha^\prime}$
for $\alpha^\prime > \alpha$
and $\alpha^\prime$ approaching $\alpha$.

We will rely on the following tightness result for $C^{\alpha,t}_{t^\beta}$
started in $\Lambda^{\alpha,t}$, whose somewhat technical proof is
given in Section~\ref{S:proofmobo}.
Denote by
\be{e:055y}
  \big\{C^{\alpha,t};\,\alpha\in(0,1]\big\},
\ee
the collection of coupled coalescent processes constructed in
(\ref{grevcorr}), (\ref{grevcorr2}),
where we abbreviate
\be{e:055z}
   \big(C^{\alpha,t}_s, L^{\alpha,t}_s\big)_{s \geq 0}
 :=
   \big(C^{\CI_{\Lambda^{\alpha,t}}}_s,
   L^{\CI_{\Lambda^{\alpha,t}}}_s\big)_{s\geq 0}.
\ee
\begin{proposition}[Uniformly bounded expectation on logarithmic scale]
\label{P:tight}
There are finite constants $M$ and $t_0$ such that for all $t\ge t_0$,
satisfying $\alpha\in (0,\infty)$, and $\beta\in(\alpha,3\alpha/2)$,
\be{e:013}
   {\bf E}\big[\# C^{\alpha,t}_{t^\beta}\big]
 \le
   M\,\big\{\frac{\alpha}{2(\beta-\alpha)}\vee
   \frac{{\bf E}\big[\# C^{\alpha,t}_2\big]}{t^{\alpha}}\vee 1\big\}.
\ee
\end{proposition}\sm

\begin{remark}\label{Rem:10}\rm
The $C^{\alpha,t}_2$ in (\ref{e:013}) denotes the coalescent
partition evaluated at time $2$,
any finite positive time could be taken instead of $2$ here, and the
two constants $t_0$ and $M$ would change accordingly.
Our special choice of the time point $2$ is convenient from the
perspective of the time discretization used in the proof of
Proposition~\ref{P:tight} (compare with (\ref{e:016})).
\hfill$\qed$
\end{remark}\sm

\subsection{Consequences of the expectation bound: Tightness}
\label{Sub:ConTig}
Recall notation
(\ref{e:055y})-(\ref{e:055z}) and in
addition assume that
\be{Erhofin}
   \rho
 :=
   \limsup_{t\to \infty}\sup_{z\in\Lambda^{1,t}}{\bf
     E}\big[F_z\big]
 <
   \infty.
\ee

The next result states that as $t \to \infty$ the coalescents in
(\ref{e:055y}) remain finite and localized in certain boxes.
\begin{proposition}[The asymptotically infinite spatial case: small time
  scales]
\label{P3.1a}
Consider the coalescent restricted
to $\Lambda^{\alpha,t}$.
Let $t_0$ be as specified in Proposition~\ref{P:tight}. Then the
following holds.
\begin{itemize}
\item[(a)]
For each fixed $\alpha^\prime>\alpha$,
there exists a sequence $(a_N)_{N\in\N}\uparrow 1$
such that for all $N\in\N$,
\be{agr1}
   \inf_{t\geq t_0}{\bf P}\big\{\# C^{\alpha,t}_{t^{\alpha^\prime}}
   \leq  N\big\}
 \geq
   a_N,
\ee
and ($\sim t$ denoting the equivalence relation w.r.t. time $t$ partition)
\be{vli1}
   \liminf_{t\to\infty}{\bf P}\big\{\max_{{i}}
   \parallel L_{t^{\alpha'}}^{\alpha,t}( \{{i}\sim^{t^{\alpha'}}\})\parallel
   \leq  t^{\alpha'/2}\log t\big\}
 \geq
   a_N.
\ee
\item[(b)] For each fixed $\alpha'>\alpha$ and each $N\in\N, L^{\alpha,t}_{t^{\alpha'}}$,
the set of all marks at time $t^{\alpha'}$ and
$I_\alpha (1,t)$ as in (\ref{EI}) we have:
\be{agr2}
   {\bf P}\big(\{L^{\alpha,t}_{t^{\alpha^\prime}} \mbox{ is
contained in }
   I_{\alpha'}(1,t)\} \big | \# C^{\alpha,t}_{t^{\alpha^\prime}}\leq N\big)
 \tto
   1.
\ee
\item[(c)] For each $N\in\N$,
\be{agr3}
%{\underline \liml}_{t \to \infty}
%{_{\D{\underline \lim}\atop t \to \infty}}
%\underline \liml_{t \to \infty}
   \lim_{\alpha^\prime \downarrow \alpha}
   \liminf_{t \to \infty}{\bf P}\big\{\# C^{\alpha,t}_{t^{\alpha^\prime}}
    \geq N\big\}
 =
   1. % \quad \mbox{as } \alpha^\prime \downarrow \alpha.
\ee
\end{itemize}
\end{proposition}\sm

\begin{proof}[{\bf Proof}]
Assertion (\ref{agr1}) is now an immediate consequence of
Proposition~\ref{P:tight}
and the Markov inequality.

Assertion (\ref{vli1}) follows from a large deviation estimate.
It will be convenient here and below to set
\be{Ealphastar}
   \alpha^\ast
 =
   \alpha^\ast(\alpha,\alpha')
 :=(\alpha+\alpha')/2.
\ee

Let $\{(X_s^i)_{s\geq 0};\,i\geq 1\}$
be an infinite collection of independent random walks with
kernel $a(\boldsymbol{\cdot},\boldsymbol{\cdot})$ such that the
initial locations $\{X_0^i;\,i\geq 1\}$
are distributed as the location process $L_0^{\alpha,t}$ of the
coalescent restricted to the box
$\Lambda^{1,t}$.
Take $\varepsilon <(\alpha'-\alpha)/2$ so that
$\alpha^\ast + \varepsilon < \alpha'$.
Since (\ref{Erhofin}) holds, we have that
 $\# C^{\alpha,t}_0$ is bounded by $2\rho t^{\alpha}$ with
overwhelming probability. Due
to a large deviation estimate (for example,
(\ref{e:026x}) is more than needed here)
\be{lim1}
   \lim_{t\to\infty}{\bf P}\big\{\max
   \parallel X_{t^{\alpha^\ast}}^i\parallel
   > t^{(\alpha^\ast + \varepsilon)/2}
   : i \in \{1,\cdots, \lfloor 2 \rho t^\alpha\rfloor \}\big\} = 0,
\ee
and hence
\be{lim2}
   \lim_{t\to\infty}{\bf P}\big\{
   \max\parallel
   \{X_{t^{\alpha^\ast}}^i\parallel> t^{(\alpha^\ast + \varepsilon)/2}
    : i \in \{1,\cdots, \# C^\alpha,t_0\} \big\}
 =
   0.
\ee
Therefore
\be{Ebootbeg}
   \lim_{t\to\infty}{\bf P}\big\{\max_{{i}}\parallel
   L_{t^{\alpha'}}^{\alpha,t}( \{{i}\sim^{t^{\alpha'}}\})\parallel
   >t^{(\alpha' + \varepsilon)/2}\big\}
 =
   0.
\ee

In order to get (\ref{vli1}) from (\ref{Ebootbeg}) we use (\ref{agr1})
together with the fact that during the remaining time
$t^{\alpha'} - t^{\alpha^\ast}$
none of the
finitely many partition classes reaches distance larger than $t^{\alpha'/2}\log t$,
with overwhelming probability.

In order to prove (\ref{agr2})
%we use the random walk asymptotics derived in the
%work of Cox and Griffeath \cite{CoxGri86} and adapt the arguments to
%delayed coalescing processes.
fix  $\alpha'>\alpha$.

Fix $N\geq 1$, and note that (\ref{agr1}) implies
the uniform lower bound $\bar{p}$ on the probability of
$\{\# C^{\alpha,t}_{t^{\alpha'}}\leq N\}$.
So (\ref{agr2}) will follow provided we show that for any $\varepsilon
>0$
we have
\be{ElowP}
   {\bf P}\big\{L^{\alpha,t}_{t^{\alpha^\prime}}
   \mbox{ is contained in } I_{\alpha'}(1,t)\big\}
 \geq
   1-\varepsilon \bar{p}.
\ee

Again due to part (a), it is possible to pick $M_\varepsilon$ so that
$C^{\alpha,t}_{t^{\alpha^\ast}}$ contains
at most $M_\varepsilon$ equivalence classes, with probability higher than
$1-\bar{p}\varepsilon/3$,
and such that any pair of them is at mutual distance
smaller than $2 t^{\alpha^\ast/2}\log t$ with probability higher than
$1-\bar{p}\varepsilon/3$.
During the remaining time interval $(t^{\alpha^\ast},t^{\alpha'}]$ of length
$t^{\alpha'} - t^{\alpha^\ast}$, which is of order
$t^{\alpha'}$, each pair of non-coalescing walks
(out of at most ${M_\varepsilon \choose 2}$ many pairs) achieves, with
overwhelming probability, a mutual distance of order $N(0,1)\times$
$2 \sigma^2
t^{\alpha'/2}$, which is with overwhelming probability in the interval
$I_\alpha(1,t)$.
The set of distances between pairs of elements of $C^{\alpha,t}_{t^{\alpha'}}$
is a subset of the set of distances between
the pairs of above random walks.
Therefore, one can choose $t$ large enough so that
\be{grev10}
\begin{aligned}
   &{\bf P}\big\{L^{\alpha,t}_{t^{\alpha^\prime}} \mbox{ is
   contained in } I_\alpha(1,t)\}\big| \# C^{\alpha,t}_{t^{\alpha^\ast}}
   \leq M_\varepsilon,\max_{{i}}\parallel
   L_{t^{\alpha'}}^{\alpha,t}(\{{i}\sim^{t^{\alpha'}}\})\parallel
   \leq  2 t^{\alpha^{\ast/2}} / log t  \big\}
  \\
 &\geq
   1-\bar{p}\varepsilon/3,
\end{aligned}
\ee
so (\ref{ElowP}), and therefore (\ref{agr2}) holds.

It still remains to prove (\ref{agr3}).
%We have to prove that
%$C^{\alpha,t}_{t^{\alpha^\prime}}$
%with $\alpha^\prime = \alpha - \ve$
%diverges stochastically as $\ve \downarrow 0$ and that it does so in a monotone way.
Fix $\alpha'>\alpha>0$. Note
that for any $N$ particles started at locations $x_1,\ldots,x_N$
contained in $I_\alpha(1,t)$ we have by convergence of the first
component in (\ref{Enocoal}) that
\be{grev11}
   {\bf P}\big\{\mbox{no coalescence by time }t^{\alpha'}\big\}
 \tto
   (\alpha/\alpha')^{{N \choose 2}}.
\ee

For fixed $N$ first choose
large $t$ so that it is possible to find $N$ particles from the initial
configuration at time $0$
with locations contained in $I_\alpha(1,t)$, and then note that as
$\alpha'\downarrow \alpha$ the right hand side above converges to $1$.
\end{proof}\sm

\section{Large time-space scale asymptotics of coalescent}
\label{S:ltsas}

In this section we combine the results of Sections~\ref{S:prelim},~\ref{S:ASP} and~\ref{S:ADPST} to prove Theorems~\ref{T:01} through~\ref{T:03}.

\subsection{Proof of Theorem \ref{T:01}}
\label{Sub:prooft}
Fix $1\ge\alpha'>\alpha>0$ and $\varepsilon\in(0,\alpha'-\alpha)$.

By
(\ref{agr3}), for all $N\in\N$ there exists
$\alpha^\ast\in(\alpha,\alpha+\varepsilon)\subset(\alpha,\alpha')$
and $t_1=t_1(N)$  such
that for all $t\ge t_1$,
\be{Elarprob}
   \mathbf{P}\big\{\# C^{\alpha,t}_{t^{\alpha^\ast}}\geq N\big\}
 \ge
   1-\varepsilon.
\ee

%It can easily be seen that $\alpha^\ast(N)\to \alpha$ as $N \to\infty $.
From now on assume that $t\geq t_0$  where $t_0$ is specified as in
Proposition~\ref{P:tight}.
Proposition~\ref{P3.1a} implies that with
probability tending to $1$ as $t\to \infty$,
the configuration $C^{\alpha,t}_{t^{\alpha^\ast}}$ has
finitely many
particles in locations contained in $I_{\alpha^\ast}(1,t)$.

Put
\be{grev13}
   n^{\alpha^\ast,t}
 :=
   \# C^{\alpha,t}_{t^{\alpha^\ast}}.
\ee

Then Proposition~\ref{P3.1} joint with Proposition~\ref{P3.1a} (a)
and (b), yield
\be{Eresfirs}
   d_{\mathrm{Pr}}\Big(
   \CL\big[(\# C^{\alpha,t}_{t^\beta})_{\beta\in[\alpha^\ast,\infty)}\big],
   \CL\big[(\# K^{n^{\alpha^\ast,t}}_{\log(\beta/\alpha^\ast)})_{\beta\in
   [\alpha^\ast,\infty)}\big]\Big)
 \tto
   0,
\ee
where $d_{\mathrm{Pr}}$ is the Prohorov metric  which is known to
metrize the weak topology (see, for example,
\cite{EthierKurtz86}).  Moreover, for a random variable $n$ and $s\ge 0$,
$\# K^n_s$ is a random variable which, given $n=k$, is distributed as
the Kingman coalescent started in $\{\{1\},\{2\},\ldots,\{k\}\}$ and
evaluated at time $s$.
%In particular,
%\be{Eresfir}
%   d_{\mathrm{Pr}}\Big(
%   \CL\big[(\# C^{\alpha,t}_{t^\beta})_{\beta\in[\alpha',\infty)}\big],
%   \CL\big[(\# K^{n^{\alpha',t}}_{\log(\beta/\alpha^\ast)})_{\beta\in
%   [\alpha',\infty)}\big]\Big)
% \tto
%   0,
%\ee

Recall that we denote by $ K_{\boldsymbol{\cdot}}$ the Kingman
coalescent started
from the trivial infinite partition $\{\{i\}:\,i\in \N\}$.
Easy properties of the Kingman coalescent guarantee that
for all $\delta>0$,
\be{EKineasy1}
   \big(\# K^n_s\big)_{s\geq\delta}
 \Tno
   \big(\# K_s\big)_{s\geq\delta},
\ee
and
\be{EKineasy2}
   \big(\# K^\infty_{s+u}\big)_{s\geq \delta}
 \TuO
   \big(\# K^\infty_s\big)_{s\geq \delta}.
\ee

Note that Proposition~\ref{P3.1a}(a) insures that the family
$\{n^{\alpha^\ast,t};\,t\geq t_0\}$ is tight.
Choose $(t_m)\to\infty$ and $n^{\alpha^\ast}$ such that
$n^{\alpha^\ast,t_m}\to n^{\alpha^\ast}$, as $m\to\infty$.
Then $n^{\alpha^\ast}$ is a finite
random variable and
\be{grev14}
   \big(\# C^{\alpha,t_m}_{t_m^\beta}\big)_{\beta\in[\alpha',\infty)}
 \Tmo
   \big(\# K^{n_{\alpha^\ast}}_{\log(\beta/\alpha^\ast)}\big)_{\beta\in
   [\alpha',\infty)}.
\ee

The left hand side of (\ref{grev14}) does
not depend on $\varepsilon$.
By (\ref{Elarprob}), (\ref{EKineasy1}) and (\ref{EKineasy2})
we have, after letting $\varepsilon\to 0$,
\be{grev15}
   \big(\# C^{\alpha,t_m}_{t_m^\beta}\big)_{\beta\in[\alpha',\infty)}
 \Tmo
   \big(\# K_{\log(\beta/\alpha)}\big)_{\beta\in
   [\alpha',\infty)}.
\ee

Since one obtains the same limit regardless of the choice of
the subsequence $(t_m)$, the statement of the theorem follows.

\subsection{Proof of Theorem \ref{T:02}}
\label{Sub:infty}
Recall from (\ref{Calphatrho}) that
$\{(C^{\alpha,t,\rho}_{s},L^{\alpha,t,\rho}_{s})_{s\ge 0};\,\alpha\in
(0,1]\}$ denotes the family of
spatial coalescents on $\Lambda^{\alpha,t}$ corresponding to the
parameter $\rho\in(0,\infty]$.

Recall the initial states $\{F^A_z;\,z\in A\}$ from (\ref{vl18}).
In this section we assume that
$\{F^\rho_z;\, z\in\Lambda^{1,t},\rho\ge 1\}$
is for fixed $\rho$ a family of independent
identically distributed random
variables with Poisson($\rho$) distribution.
In fact, due to thinning and superposition properties of
the Poisson process on the line
we can consider a coupling
such of the families for different $\rho$ that if $\rho_1\le\rho_2$ then
\be{grev17}
  F_{z}^{\rho_1}\leq F_{z}^{\rho_2},
\ee
for all $z\in\Lambda^{1,t}$.

Due to this coupling and the monotonicity properties collected in
Subsection~\ref{Sub:monoton},
\be{grev18}
   \big(C^{\alpha,t,\rho}_{s},L^{\alpha,t,\rho}_{s}\big)_{s\ge 0}
 \rhoto
   \big(C^{\alpha,t,\infty}_{s},L^{\alpha,t,\infty}_{s}\big)_{s\ge 0},
\ee
here convergence is meant in the sense of convergence
defined (\ref{vl10a}).

The goal of this subsection is to show that
the results obtained in Subsection~\ref{Sub:ConTig}
hold in the limit $\rho\to\infty$.

Fix $\delta>0$ and recall from (\ref{calphatrh2}) the spatial coalescent
$(C^{\alpha,t,\infty,\delta}_s)_{s\ge 0}$,
thinned out by those particles which were attempted to jump in the time period
$[0,\delta]$.

%We denote thus constructed thinned spatial coalescent
%using the additional superscript $\delta$, for example
%\be{GR15b}
%C^{\alpha,t,\infty,\delta}_s,
%\ee
%denotes its partition structure at time $s\geq 0$.

\begin{lemma}[The limit of infinite density]\label{L3.1b}
%Let $(C^{\alpha,t,\infty,\delta}_s,L^{\alpha,t,\infty,\delta}_s)_{s\geq 0}$
%be the thinned spatial coalescent started with
%countably many partition elements as described above.
For each $\delta>0$ fixed,
\be{agr4}
   \lim_{N\to\infty}\liminf_{t \to \infty}
   {\bf P}\big\{\# C^{\alpha,t,\infty,\delta}_{t^{\alpha^\prime}}\leq
   N\big\}
 =
   1.
\ee
\end{lemma}\sm

\begin{proof}[{\bf Proof}]
Recall from  Lemma~\ref{LPoisson} that the number of partition
elements of a Kingman coalescent can be dominated by a Poisson
variable with suitably large parameter $\rho_0$. By monotonicity we can
construct a coupling
\be{agr4b}
   \Big(\big(C^{\alpha,t,\infty,\delta}_s,L^{\alpha,t,\infty,\delta}_s\big),
   \big(C^{\mathrm{Poisson}(\rho_0)+1}_s,L^{\mathrm{Poisson}(\rho_0)+1}_s\big)
   \Big)_{s\geq\delta},
\ee
where
$(C^{\mathrm{Poisson}(\rho_0)+1}_{s},
L^{\mathrm{Poisson}(\rho_0)+1}_{s})_{s\ge 0}$ is started from
the initial configuration where
$\{F_{z};\,z\in\Lambda^{1,t}\}$ is a family of independent random
variables which equal in distribution one plus a rate $\rho_0$ Poisson
distributed random variable such that
% for $\rho$ chosen as in Lemma \ref{LPoisson},
%so that
$C^{\alpha,t,\infty,\delta}_s\leq C^{\mathrm{Poisson}(\rho_0)+1}_s$,
for all $s\geq \delta$, almost surely.
The statement now follows from Proposition~\ref{P3.1a}(a) applied to
$(C^{\mathrm{Poisson}(\rho_0)+1}_s)_{s\ge 0}$.
\end{proof}\sm

In addition, notice that
${\bf P}\{\# C^{\alpha,t,\infty,\delta}_\delta\geq 1\}=1$,
so if $(C_s^1,L_s^1)_{s\geq \delta}$ is the family of
spatial coalescents started with $1$ particle at each site
of $\Lambda(t^{1/2})$,
we have
\be{agr4c}
   C_\delta^1
 \leq
   C_\delta^{\alpha,t,\infty,\delta}
 \leq
   C_\delta^{\mathrm{Poisson}(\rho)+1}.
\ee

The extension of Theorem \ref{T:01} (in Proposition \ref{P:02})
proved in Subsection \ref{Sub:prooft} clearly applies to both the left-most and the right-most
family of coalescents. Therefore, by Lemma~\ref{Lsandwich},  for fixed
$\alpha'>\alpha>0$,
\be{agr5bhelp}
   \big(\# C^{\alpha,t,\infty,\delta}_{t^\beta}\big)_{\beta\in
   [\alpha',\infty)}
 \Tto
   \big(\# K_{\log(\beta/\alpha)}\big)_{\beta\in[\alpha',\infty)},
\ee
and Theorem \ref{T:02} follows.

\subsection{Proof of Theorem \ref{T:03}}
\label{Sub:EAPA}
%Recall the metric (\ref{Emetricaug}).

The proof of Theorem \ref{T:03} makes use of a convergence result
stated in Theorem~1 in \cite{donnelly},
which applies in a much more general setting than ours.
For the benefit of the reader, we will rephrase it in our setting.\sm

\begin{lemma}[Donnelly, 1991]
\label{L:Do1}
Suppose $\{(B_{s}^N)_{s\ge 0};\,N\geq 1\}$ is a family of
$D([0,\infty),\N)$-valued
random variables which satisfy the
following three assumptions:
\begin{itemize}
\item[{(A1)}]
For all $N\in\N$, $l\ge n\in \N$, $s\ge\alpha$ and $y\ge 1$,
\be{grev40}
   \mathbf{P}\big(\inf_{u\in [\alpha,s]} B_u^N\le y|
   B_\alpha^N=l\big)
 \leq
   \mathbf{P}\big(\inf_{u\in [\alpha,s]} B_u^N\le y|
   B_\alpha^N=n\big).
\ee
\item[{(A2)}]
For all $n\in \N$,
\be{grev41}
   \CL\big[(B_u^N)_{u\geq \alpha}|B_\alpha^N=n\big]
 \Nto
   \CL\big[(\#K_{\log(u/\alpha)}^n)_{u\ge\alpha}\big],
\ee
\item[{(A3)}]
Suppose we have a sequence $(n_M)\to\infty$, such that for each $u>\alpha$,
\be{grev42}
   \lim_{M\to \infty}\lim_{N\to \infty}\mathbf{P}\big(
   B_u^N \leq M|B_\alpha^N=n_M\big)=1.
\ee
\end{itemize}

Then
\be{don_conc}
   \CL\big[(B_u^N)_{u\geq \alpha}|B_\alpha^N=n_N\big]
 \Nto
   \CL\big[(\# K_{\log(u/\alpha)})_{u\geq \alpha}\big].
\ee
\end{lemma}\sm

\begin{proof}[{\bf Proof of Theorem \ref{T:03}}]
(i)
Take a subsequence $(t_N)\uparrow\infty$, and let
\be{sN}
   s_N
 :=
   \#\Lambda^{\alpha,t_N}.
\ee
We consider first a {\em special} case.
Draw $\mathrm{Bin}_N$ according to the Binomial distribution with parameters
$s_N$ and $p\in(0,1]$ or the Poisson distribution with parameter
$s_N\cdot\rho$. Given
$\mathrm{Bin}_N=k$, place $k$ particles uniformly
without and with replacement  at $k$ positions in
$\Lambda^{\alpha,t_N}$.
Notice that the random configurations obtained this way will equal in law to
$C^{\alpha,t_N}_{0}$ under the assumption that $\{F_{z};\,z\in\Lambda^{\alpha,t_N}\}$ are
independent and identically distributed random variables with the
Bernoulli (parameter $p$) or with the Poisson($\rho$) distribution, respectively.

Put for all $u\ge\alpha$,
\be{grev43}
   B_u^N
 :=
   \# C^{\alpha,t_N}_{(t_N)^u}.
\ee
By Lemma~\ref{L3.1mar}, given that $\mathrm{Bin}_N =k$,
with probability tending to $1$, $B_\alpha^N=k$.
The advantage of the above construction(s) is that the assumption (A1)
is automatically satisfied provided we keep the same algorithm
for ``positioning
the $k$ particles in
$\Lambda^{\alpha,t_N}$'', for all $k\in\N$, i.e.,
provided that for each $k\in\N$ and all $l< k$, the
first $l$ points in
$C^{\alpha,t_N}_{0}$ given $\mathrm{Bin}_N=l$ match those
in $C^{\alpha,t_N}_{0}$ given $\mathrm{Bin}_N=k$.
The assumptions (A2) and (A3) (provided that $n_N=O(s_N)$)
are implied by Lemma~\ref{L3.1mar} and (\ref{agr1}), respectively.
Therefore, (\ref{don_conc}) holds in the (special) Binomial case for any $p\in(0,1]$ and any
sequence $n_N\leq s_N$ going to $\infty$.
Similarly, (\ref{don_conc}) holds in the (special) Poisson case for any $\rho\in(0,\infty)$ and
any sequence $n_N=O(s_N)$.\\
In particular, if $p=1$ and $n_N=s_N$ (almost surely) then
\be{grev44}
   \big(\# C^{\alpha,t_N}_{{t_N}^\beta}\big)_{\beta\geq \alpha}
 \Nto
   \big(\# K_{\log(\beta/\alpha)}\big)_{\beta\geq \alpha}.
\ee
Since the limit is uniform in the choice of the subsequence $t_N\to
\infty$, we conclude the statement of the theorem in this case.

The general Bernoulli($p$) case can be dealt with similarly as the
general Poisson($\rho$) case, as we explain next.
Fix $\rho\in(0,\infty)$ and note that
(\ref{don_conc}) holds both with $n_N:= \lfloor \rho s_N/2\rfloor$ and with
$n_N:= \lfloor 2\rho s_N\rfloor$.
Since with probability tending to $1$, the Poisson
($\rho s_N$) distributed random variable $\mathrm{Bin}_N$ satisfies
\be{grev45}
   \lfloor\rho s_N/2\rfloor
 \leq
   \mathrm{Bin}_N
 \leq
   \lfloor 2 \rho s_N\rfloor,
\ee
we can apply Lemma~\ref{Lsandwich} to conclude the needed statement as
done before. \sm

(ii)
Note that due to part (i), the family of processes
$(t_N$ is a sequence with $t_N \uparrow \infty$ as $N \to \infty$)
\be{grev46}
   \big(\# C^{\alpha,t_N}_{t_N^\beta}\big)_{\beta \in [\alpha,\infty)},
\ee
where the family $\{F_{z};\,z\in\Lambda^{1,t}\}$
is drawn from the ``Poisson($\rho$)+1'' distribution
%as in Lemma~\ref{LPoisson},
is tight in $D([\alpha,\infty),\bar{\N})$
since we can sandwich it between from below the case where we start with exactly $1$
particle per site (Bernoulli with $p=1$) and from above with the independent sum of two spatial coalescent processes  one started in
$\mathrm{Poisson}(\rho)$ and the other one with exactly $1$
particle per site (Bernoulli with $p=1$). Here we use monotonicity in
$\beta$ for every $N$.
Moreover, the process
$(\# K_{\log(\beta/\alpha)})_{\beta\in[\alpha,\infty)}$ is
the only possible
(subsequential) limit
due to Theorem \ref{T:01}.
Therefore, applying monotonicity and using (\ref{agr4c}) as in
the proof of Theorem~\ref{T:02} implies the statement.
\end{proof}\sm

\section{Convergence on the spatial scale (Proof of Theorems~\ref{T:04}
and~\ref{T:05})}
\label{S:ConvSpatScale}

In this section we prove results which involve the coalescent with
rebirth using the results established in  Sections~\ref{S:ASP}
and~\ref{S:ADPST}.

\subsection{Proof of Theorem \ref{T:04}}
\label{Sub:tskoro}
Consider the family $\{\rho_{{\mathcal I}^\alpha}\diamond
C^{1,t};\,\alpha\in(0,1]\}$ from
(\ref{grevcorr}).  Fix $\rho\in(0,\infty)$, and assume
(\ref{Erhofin}).
By Theorem~\ref{T:01} (with $\beta=1$) and (\ref{grevcorr2})
we already know that, for a fixed $\alpha\in(0,1]$,
\be{agr6}
   \# \rho_{{\mathcal I}^\alpha}\diamond C^{1,t}_t
 \Tto
   \# K_{-\log{\alpha}}.
\ee
Our first and key goal is to extend (\ref{agr6}) to the f.d.d.~convergence
of $(\# \rho_{\CI^\alpha} \diamond C^{1,t}_t)_{\alpha\in[\alpha_l,\alpha_u]}$
to $(N_\alpha)_{\alpha\in[\alpha_l,\alpha_u]}$,
where $\alpha_l,\alpha_u\in (0,1)$,
stated below in Proposition~\ref{P.Cofdd1}.
In particular, here $N_\alpha$ is the number of partition elements of
$K_0[\log{\alpha_l},\log{\alpha_u}]$ born before time
$\log{\alpha}$, as defined in (\ref{agr5e}).
As a second (small) step we derive at the end the
convergence on path space as stated in Theorem~\ref{T:04}.
\begin{proposition}[Partition number f.d.d.~convergence]\label{P.Cofdd1}
\begin{itemize}
\item[{}]
\item[(a)] Fix $\rho \in [0,\infty)$.
\begin{itemize}
\item[(i)] For all $m\in\N$ and $\alpha_l\le\alpha_1<\ldots<\alpha_m\le\alpha_u$,
\be{agr6a}
  \big(\# \rho_{{\mathcal I}^{\alpha_1}}\diamond C^{1,t}_t,\ldots,\# \rho_{{\mathcal I}^{\alpha_m}}\diamond
  C^{1,t}_t\big)
 \Tto
  \big(N_{\alpha_1},\ldots,N_{\alpha_m}\big).
\ee
\item[(ii)] For any $\varepsilon>0$, the family
\be{agr7}
   \Gamma
 :=
   \Big\{\CL\big[\#\rho_{{\mathcal I}^{\alpha}}\diamond C^{1,t}_t\big];\,\alpha\in
   [0,1-\ve],t>0\Big\}
\ee
is tight.
\end{itemize}
\item[(b)] The statements of (a) remain valid if $\rho=\infty$ and,
for a fixed
$\delta>0$, $C^{1,t}_t$ is replaced by
$C^{1,t,\infty,\delta}_t$.
\end{itemize}
\end{proposition}\sm

\begin{remark}
Note that the generalization of the proposition
in terms of the corresponding convergence of the partition structure
%as in Theorem \ref{T:03b} using the set-up of the end of subsection \ref{SubSub:alpha}.
could be formulated in the setting of finite (and sparse) initial configurations considered in
the proof of (a.i) below,
and proved by applying
the technique of Section \ref{Sub:convps} (see also Lemma 7.3 in \cite{GreLimWin05} or
Proposition 14 in \cite{limstu06}).
\hfill $\square$ \end{remark}\sm

Before giving the argument we present a key tool.
For $m\in\N$, fix parameters $0<\alpha_l\le\alpha_1<\alpha_2<\ldots<
\alpha_m\le\alpha_u$.
To make the argument more transparent we introduce the coalescent
with {\em rebirth at
finitely many prescribed times} $\{\log{\alpha_1},\ldots,\log{\alpha_m}\}$
only and call this process:
\be{agr6b}
   \big(K_s^{\vec{\alpha}}\big)_{s\geq \log \alpha_\ell}
 :=
   \Big(K_s\big[\{\log{\alpha_1},\ldots,\log{\alpha_m}\}\big]\Big)_{s\geq \log \alpha_\ell}.
\ee
In words, the process $(K_s^{\vec{\alpha}})_{s\geq \log \alpha_\ell}$
behaves as follows: it starts in the configuration
$\{\{(n,\log{\alpha_1})\};\,n\in\N\}$
%, the time coordinate grows at unit speed,
and during each interval of the
form $[\log{\alpha_{(k-1)}},\log{\alpha_k})$, the partition-valued
component behaves like a
coalescing process without rebirth, while at times of the form
$\log{\alpha_k}$ the partition elements that were ``lost'' during the time
interval $[\log{\alpha_{(k-1)}},\log{\alpha_k})$, (i.e., their
 label $(n,\log{\alpha_{k-1}})$
 does not label any partition element of
$K_{\log{\alpha_k}-}^{\vec{\alpha}}$)
get ``reintroduced'' at time $\log{\alpha_k}$ as partition elements
$\{(n,\log{\alpha_k})\}$, $k=1,\ldots,m$.
Compare Figure~{4} for an illustration.

To define the new process formally
replace in $(K^{\mathbb{S}}_s)_{s\ge \log \alpha_\ell}$ every $(i,t)$ by
$(i,\log{\alpha_k})$ if $t\in[\log{\alpha_{(k-1)}},\log{\alpha_k})$,
$k=1,\ldots,m$, where $\alpha_0:=0$.

We are particularly interested in the state of
$K^{\vec{\alpha}}$ at time $0$.
We shall show that it
%$(K_s[\{\log{\alpha_1},\ldots,\log{\alpha_m}\}])_{s\ge \log \alpha_\ell}$
agrees with the coalescent with
(continuous) rebirth $(K_s[\log{\alpha_l},\log{\alpha_u}])_{s\ge \log \alpha_\ell}$
with respect to the functional of interest. This observation will then imply the
statement once we have handled the case of finitely many rebirth
times introduced above.

Denote by $\tilde{N}_k^{\vec{\alpha}}$
the total number of partition elements of
$K_0^{\vec{\alpha}}$
with birth time {\em equal to or smaller than} $\log{\alpha_k}$.
Note that $1\leq \tilde{N}_1^{\vec{\alpha}}\leq \tilde{N}_2^{\vec{\alpha}}\leq\ldots\leq
\tilde{N}_m^{\vec{\alpha}}$, almost surely.
By construction, it is not difficult to verify
that the following key identity holds,
\be{agr6b2}
 \tilde{N}^{\vec{\alpha}}
  =
   \big(N_{\alpha_1},\ldots,N_{\alpha_m}\big).
\ee

%%%% Figure 1 %%%%%%%%
\psset{xunit= 1.05cm,yunit= 1cm}
%\psset{xunit= 0.4725cm,yunit= 0.45cm}
\pspicture(-2,-10.1)(15,-.2)
\psset{linewidth=.5pt}
%draw lines
\psline(0,-2.8)(0,-8.0)
\psline{<-}(0,-0.8)(0,-1.3)
\psline(0,-8)(3,-8)
\psline(6,-8)(9,-8)
\psline{->}(9.3,-8)(10,-8)
\uput[270](10,-8){$s$}
%draw branches
\psset{linewidth=.2pt}
\pnode(0,-7){A}
\pnode(0,-6){B}
\pnode(3,-7){C}
\pnode(3,-8){D}
\pnode(3,-6){E}
\pnode(3,-5){F}
\pnode(6,-4){G}
\pnode(6,-3){H}
\pnode(6,-6){I}
\pnode(6,-7){J}
\pnode(2,-6.5){AB}
\pnode(4.5,-7.5){CD}
\pnode(5,-5.5){EF}
\pnode(8.5,-4.1){GH}
\pnode(8,-6.4){IJ}
\twocoaln{A}{B}{AB}
\twocoaln{C}{D}{CD}
\twocoaln{E}{F}{EF}
\twocoaln{G}{H}{GH}
\twocoaln{I}{J}{IJ}
%%% Anita's extension %%%
\pnode(0,-5){K}
\pnode(0,-4){L}
\pnode(0,-3){M}
\pnode(3,-4){O}
\pnode(3,-3){P}
\pnode(6,-5){R}
\pnode(9,-6){T}
\psline(6,-5)(9,-6)
\pnode(2,-5.6){KL}
\twocoaln{L}{K}{KL}
\psline(2,-5.6)(3,-6)
\psline(2.5,-5.8)(0,-3)
\psline[linearc=.6](3,-3)(5,-3.2)(6,-4)
%%%
\psline(2,-6.5)(3,-7)
\psline(4.5,-7.5)(6,-7)
\psline(5,-5.5)(6,-6)
\psline(8.5,-4.1)(9,-4)
\psline(8,-6.4)(9,-7)
\psline[linearc=.6](3,-4)(5,-4.2)(6,-5)
%\psline[linearc=.6](2,-5.7)(.5,-5.7)(3,-6) \psset{linewidth=.5pt}
%draw circles
\pscircle*(0,-8){.1}
\pscircle*(0,-7){.1}
\pscircle*(0,-6){.1}
\pscircle*(0,-5){.1}
\pscircle*(0,-4){.1}
\pscircle*(0,-3){.1}
%\pscircle*(0,-2){.1} \pscircle*(0,-1){.1} \pscircle*(0,0){.1}
\pscircle*(3,-8){.1}
\pscircle*(3,-7){.1}
\pscircle*(3,-6){.1}
\pscircle*[linecolor=gray](3,-5){.1}
\pscircle*[linecolor=gray](3,-4){.1}
\pscircle*[linecolor=gray](3,-3){.1}
%\pscircle*[linecolor=gray](3,-2){.1}
%\pscircle*[linecolor=gray](3,-1){.1}
%\pscircle*[linecolor=gray](3,0){.1}
\pscircle*(6,-7){.1}
%\pscircle*[linecolor=gray](6,-8){.1}
\pscircle*(6,-6){.1}
\pscircle*[linecolor=gray](6,-5){.1}
\pscircle*[linecolor=gray](6,-4){.1}
\pscircle[fillstyle=solid,fillcolor=white](6,-3){.1}
%\pscircle[fillstyle=solid,fillcolor=white](6,-2){.1}
%\pscircle[fillstyle=solid,fillcolor=white](6,-1){.1}
%\pscircle[fillstyle=solid,fillcolor=white](6,0){.1}
\pscircle*(9,-7){.1}
%\pscircle*[linecolor=gray](9,-8){.1}
\pscircle[fillstyle=solid,fillcolor=white](9,-8){.1}
\pscircle[fillstyle=solid,fillcolor=white](6,-8){.1}
\pscircle*[linecolor=gray](9,-6){.1}
\pscircle*[linecolor=gray](9,-4){.1}
%\pscircle[fillstyle=solid,fillcolor=white](9,-3){.1}
%\pscircle[fillstyle=solid,fillcolor=white](9,-2){.1}
%\pscircle[fillstyle=solid,fillcolor=white](9,-5){.1}
%\pscircle[fillstyle=solid,fillcolor=white](9,0){.1}
%dot dot dots
\psset{dotscale=.5}
\psdots*(4.4,-8.0)(4.6,-8.0)(4.8,-8.0)
\psdots*(0,-2.3)(0,-2.0)(0,-1.7)(3,-2.3)(3,-2.0)(3,-1.7)(6,-2.3)(6,-2.0)(6,-1.7)(9,-2.3)(9,-2.0)(9,-1.7)
%interval labels
\psset{linewidth=.3pt}
%\psline(6,-8.2)(6,-8.5)(6.9,-8.5)
%\uput[0](6.9,-8.5){$\log\frac{1}{\alpha_{3}}$}
%\psline(8.1,-8.5)(9,-8.5)(9,-8.2) \psline(3,-8.2)(3,-9)(5.4,-9)
\uput[0](3.0,-8.5){$\log{\alpha _{2}}$}
\uput[0](6.0,-8.5){$\log{\alpha _{3}}$}
\uput[0](9.0,-8.5){$0$}
%\uput[0](5.4,-9){$\log\frac{1}{\alpha _{2}}$}
%\psline(6.6,-9)(9,-9)(9,-8.6) \psline(0,-8.2)(0,-9.5)(3.9,-9.5)
\uput[0](0.0,-8.5){$\log{\alpha _{1}}$}
\label{FigD1}
%\uput[0](3.9,-9.5){$\log\frac{1}{\alpha _{1}}$}
%\psline(5.1,-9.5)(9,-9.5)(9,-9.1)
\endpspicture
\vspace{-1cm}

%\centerline{Figure 2}
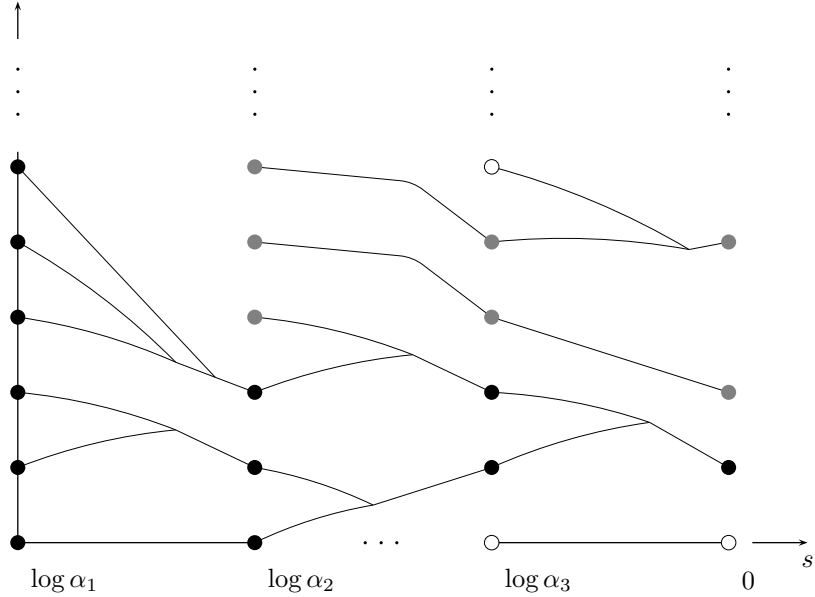
\noindent {\em Figure~{4} illustrates the evolution of the
process $K_{\boldsymbol{\cdot}}^{\vec{\alpha}}[\{\log{\alpha_1},\ldots,
\log{\alpha_3}\}]$ to the {\em left}
of time 0, where
partition elements
with birth time $\log{\alpha_1}$ are colored black, those
with birth time $\log{\alpha_2}$ are colored gray, etc.
For this realization
we see that $N_1^{\vec{\alpha}}\ge 1$,
$N_2^{\vec{\alpha}}\ge 3$, and $N_3^{\vec{\alpha}}\geq 4$.}
\bi

The scenario of Figure~{4} corresponds in the spatial set-up to the following.
For $\alpha\in (0,1]$ we refer to the set $\Lambda^{\alpha,t}$ as the
{\em $\alpha$-box}. We observe the coalescents corresponding to the
$\alpha_1$, \ldots, $\alpha_m$-boxes at times $t^{\alpha_1}$, \ldots,
$t^{\alpha_m}$
and finally at time $t$, and apply here our results from
Sections~\ref{S:ASP} and~\ref{S:ADPST}.
The structure of the arising coalescents is depicted in Figure~{5}.

\begin{proof}[{\bf Proof of Proposition~\ref{P.Cofdd1}}]

\noindent
(a.i) Theorem 6 of \cite{CoxGri86} gives some information for the case of the
{\em sparse} particles and for {\em instantaneous} coalescence
in terms of the convergence in the sense of finite dimensional distributions.

The gap between the instantaneous coalescent f.d.d.~convergence
case of \cite{CoxGri86} and our delayed coalescent
 path space convergence case is bridged as in Section \ref{S:ASP}.
It would be tedious to write out (again) all the details, yet we encourage
the reader to verify the steps of the argument outlined below.

{\em Step 1} (Sparse individuals). We first treat finitely many sparse
particles as initial state, where we can use some techniques from \cite{CoxGri86}.
As above, fix $\alpha_1,\ldots,\alpha_m$, where $0<\alpha_l\leq \alpha_1< \alpha_2<\ldots \alpha_{m-1}< \alpha_m\leq \alpha_u$.
Initially consider {\em finitely many particles} (independent of $t$)
in each of the boxes $\Lambda^{\alpha_i,t}$, $i=1,\ldots,m$,
such that,
in analogy to the statement of Proposition \ref{P3.1},
the initial positions of particles in the box $\Lambda^{\alpha_1,t}$ are contained in
$I_{\alpha_1}(c,t)$ and
moreover that,
for each $i=2,\ldots,m$:
\be{agre30}
\mbox{ the positions of all particles initially in  }
\Lambda^{\alpha_i,t}\setminus \Lambda^{\alpha_{i-1},t},
\mbox{ is in } I_{\alpha_i}(c,t).
\ee
For concreteness, assume that there are initially :
\be{agre31}
n_1 \mbox{ particles in  } \Lambda^{\alpha_1,t}, \mbox{ and }
n_i \mbox{  particles in  }
\Lambda^{\alpha_i,t}\setminus \Lambda^{\alpha_{i-1},t}, i= 2,\ldots,m.
\ee
We write for this spatial coalescent
\be{agre32}
(C^{\vec{n},t}_t, L^{\vec{n},t}_t)_{t \geq 0}.
\ee

%The Figure 3 below explains how the random vector $\tilde{N}^{\vec{\alpha}}:=(\tilde{N}^{\vec{\alpha}}_1,\ldots,\tilde{N}^{\vec{\alpha}}_m)$
%arises as a limiting object from the spatial coalescent.

%%%% Figure 3 %%%%%%%%
\psset{xunit= 1.05cm,yunit= 1cm}
%\psset{xunit= 0.4725cm,yunit= 0.45cm}
\pspicture(0,-10.6)(15,0.9)
\psset{linewidth=.5pt}
%draw vertical lines
\psline{[-}(0,0)(0,-2) \psline{[-}(0,-2)(0,-4)
\psline{[-]}(0,-4)(0,-6) \psline{-]}(0,-6)(0,-8)
\psline{-]}(0,-8)(0,-10) \psline(2,0)(2,-10) \psline(5,0)(5,-10)
\psline(9,0)(9,-10) \psline(11,0)(11,-10) \psline(13,0)(13,-10)
%label the lines
\uput[180](0,-6){$\alpha_{_1}$-box}
\uput[180](0,-8){$\alpha_{_2}$-box}
\uput[180](0,-10){$\alpha_{_3}$-box} \uput[270](0,-10){$0$}
\uput[270](2,-10){$t^{\alpha_{_1}}$}
\uput[270](5,-10){$t^{\alpha_{_2}}$}
\uput[270](9,-10){$t^{\alpha_{_3}}$}
\uput[270](11,-10){$t^{\alpha_{_k}}$} \uput[270](13,-10){$t$}
%draw branches
\psset{linewidth=.2pt} \pnode(0,-3.1){A} \pnode(0,-3.4){B}
\pnode(0,-5.4){C}
\pnode(0,-5.7){D}
\pnode(0,-6.3){E}
\pnode(0,-7.7){F}
\pnode(0,-8.3){G}
\pnode(1.1,-2.8){AB}
\pnode(0.9,-5.6){CD}
\pnode(1.8,-6){CDE}
\pnode(1.5,-8){FG}
\pnode(4,-7){CDEFG}
\twocoaln{A}{B}{AB}
\twocoaln{C}{D}{CD}
\twocoaln{CD}{E}{CDE}
\twocoaln{F}{G}{FG}
\twocoaln{CDE}{FG}{CDEFG} \psset{linewidth=.5pt}
%draw circles at time 0
\pscircle*(0,-4.3){.1} \pscircle*(0,-5.1){.1}
\pscircle*(0,-5.4){.1} \pscircle*(0,-5.7){.1}
\pscircle*[linecolor=gray](0,-6.3){.1}
\pscircle*[linecolor=gray](0,-7.1){.1}
\pscircle*[linecolor=gray](0,-7.4){.1}
\pscircle*[linecolor=gray](0,-7.7){.1}
\pscircle*[linecolor=gray](0,-2.3){.1}
\pscircle*[linecolor=gray](0,-3.1){.1}
\pscircle*[linecolor=gray](0,-3.4){.1}
\pscircle*[linecolor=gray](0,-3.7){.1}
\pscircle[fillstyle=solid,fillcolor=white](0,-.3){.1}
\pscircle[fillstyle=solid,fillcolor=white](0,-1.1){.1}
\pscircle[fillstyle=solid,fillcolor=white](0,-1.4){.1}
\pscircle[fillstyle=solid,fillcolor=white](0,-1.7){.1}
\pscircle[fillstyle=solid,fillcolor=white](0,-8.3){.1}
\pscircle[fillstyle=solid,fillcolor=white](0,-9.1){.1}
\pscircle[fillstyle=solid,fillcolor=white](0,-9.4){.1}
\pscircle[fillstyle=solid,fillcolor=white](0,-9.7){.1}
%draw some more circles
\pscircle*[linecolor=gray](1.1,-2.8){.1} \pscircle*(0.9,-5.6){.1}
\pscircle*(1.8,-6){.1} \pscircle*[linecolor=gray](1.5,-8){.1}
\pscircle*(4,-7){.1} \psset{dotscale=.5}
\psdots(-0.1,-4.9)(-0.1,-4.7)(-0.1,-4.5)(-0.1,-2.9)(-0.1,-2.7)(-0.1,-2.5)(-0.1,-6.9)(-0.1,-6.7)(-0.1,-6.5)(-0.1,-0.9)(-0.1,-0.7)(-0.1,-0.5)(-0.1,-8.9)(-0.1,-8.7)(-0.1,-8.5)
\psdots(4.5,-7)(4.8,-7)(5.1,-7)(1.5,-2.8)(1.8,-2.8)(2.1,-2.8)(9.5,-5)(9.8,-5)(10.1,-5)
%draw circles at time t
\pscircle*[linecolor=gray](13,-5){.1}
\pscircle*[linecolor=gray](13,-6.5){.1}
\pscircle*[linecolor=gray](13,-7){.1} \pscircle*(13,-3){.1}
\pscircle*(13,-4.5){.1} \pscircle*(13,-5.5){.1}
\pscircle*(13,-6){.1}
\pscircle[fillstyle=solid,fillcolor=white](13,-3.5){.1}
\pscircle[fillstyle=solid,fillcolor=white](13,-4){.1}
\pscircle*(4,-7){.1}
\label{FigD2}
\endpspicture
\vspace{.2cm}

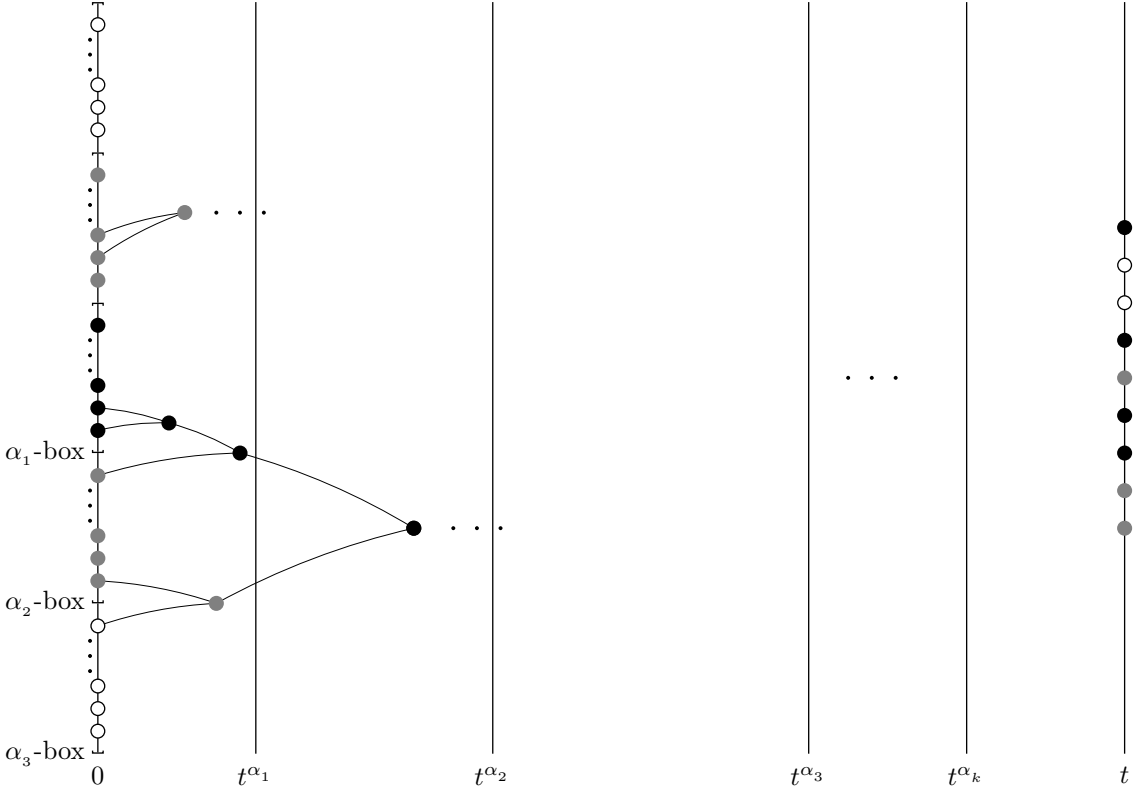
\noindent{\em Figure~{5} illustrates the occurrence of
$\tilde{N}^{\vec{\alpha}}=(\tilde{N}^{\vec{\alpha}}_{\alpha_1},\ldots,
\tilde{N}^{\vec{\alpha}}_{\alpha_m})$
in the limit of the spatial coalescent asymptotics. Notice that
the colors of the particles in Figures~{4} and~{5}
match on purpose to emphasize the correspondence between {\rm space}
(for the spatial coalescent) and {\rm time} (for the Kingman-type coalescent with rebirth).}

In \cite{CoxGri86}, Theorem 6, itis proved that if the coalescent starts
with $n_1,n_2,\cdots, n_m$ particles in
$I_{\alpha_1}(c,t), \cdots, I_{\alpha_m} (c,t)$, then
\be{agre33}
\mathbf{P} (\# C^{\vec{n}}_t =m)
\Tto
P_{n_1,\cdots,n_m;m} (\alpha_1,\cdots,\alpha_m, 1)
\ee
and in (5.3) in \cite{CoxGri86} the r.h.s. is defined by the following
recursive equation (here $t^\beta$ instead of $t$ is considered in
(\ref{agre33})
\be{agre34}
p_{n_1,\cdots,n_k;m} (\alpha_1,\cdots,\alpha_m; \beta)
=
\suml_{i_1,\cdots,i_{m-1}} p_{n_1, i_1} (\alpha_1/\alpha_2) \cdots
p_{n_m + i_{m-1},m} (\alpha_m/\beta),
\ee
with (3.10) in \cite{CoxGri86} defining the input of the recursion for
$m=1$ as
\be{agre35}
p_{n,i} (\alpha) = \mathbf{P} (| K^n_{\log (1/\alpha)} | =i).
\ee

It is straightforward to see that the Theorem 6 in \cite{CoxGri86}
now implies (with a reinterpretation of formula (5.3) and (3.10) in \cite{CoxGri86})
that
the following convergence in distribution
holds for instantaneous coalescence:
\be{agre36}
(\# \rho_{\CI^{\alpha_1}} \diamond C^{\vec{n},t}_t, \cdots , \# \rho_{\CI^{\alpha_m}}
\diamond C^{\vec{n},t}_t
\La
(N^{\vec{n}}_{\alpha_1}, \cdots, N^{\vec{n}}_{\alpha_m}), \mbox{ as } t \to \infty.
\ee
where $N^{\vec{n}}_{\alpha}$ on the r.h.s. is the number of partition elements
added to the system before time $\alpha$ in the following Kingman coalescent
with immigration evaluated at time 0. We start in with $n_1$-individuals
at time $\log \alpha_1$ and evolve until time $\log \alpha_2$ where $n_2$
new individuals are added, then continue evolving until time $\log \alpha_3$ where $n_3$
new individuals are added $\cdots$, and continue until time $\log \alpha_m$ where the last
immigration takes place. Then the coalescent runs until time 0, without
further immigration.

The point here is that the above assumptions ensure that
with overwhelming probability, for each $i=2,\ldots,m$,
none of the particles initially in $\Lambda(t^{\alpha_u})\setminus \Lambda(t^{\alpha_{i-1}})$
coalesce with any other particle
during the time interval $[0,g_{\alpha_i}(t)]$ (see (\ref{ag10}) for the definition of $g_{\alpha}(t)$),
while during the same time interval, on the appropriate time scale, the evolution of the partitions
containing particles with initial positions in $\Lambda(t^{\alpha_{i-1}})$ is approximately that of the
``Kingman coalescent with
immigration'', where at time $\log{\alpha_j}$, $j <i$, a population of size $n_j$ is adjoined
to the existing configuration.
By Lemma \ref{L1CoxGri} the partitions stay sparse with overwhelming probability,
so that the asymptotic exchangeability applies,
and an easy inductive argument
yields the convergence in this finite setting, where the limit is the
described coalescent with immigration (which is different from the limit
on the r.h.s. of (\ref{agr6a}), since here we only have sparse individuals).
Our arguments give hence a convergence statement for {\em instantaneous}
coalescence in the sparse case.

As in Lemma \ref{L3.1mar},
the convergence of \cite{CoxGri86} Theorem 6, extends to the convergence in the delayed
coalescent setting.  Moreover,
using the asymptotic exchangeability
as in Subsection \ref{Sub:convps}, this
can be extended to the convergence in path space.

{\em Step 2}
In the previous step we had finitely many sparse particles, even
as $t \to \infty$, in our problem we have in fact a growing number of
particles as $t \to \infty$ and this will lead to the actual limit in
(\ref{agr6}).

The above mentioned ``immigration'' becomes infinite in the limit as $n_i\to\infty$, $i=2,\ldots,m$.
Indeed, the reasoning of Section~\ref{S:ltsas},
in particular that of the proof of Theorem \ref{T:01},
based on the estimates of
Proposition~\ref{P:tight}
and
Proposition~\ref{P3.1a}
in Section \ref{S:ADPST}
will extend to the current setting and yield (\ref{agr6a}).
The proof is by induction on $m$. We start with $m=2$.

Let  $\rho<\infty$ and consider the joint asymptotics of
$\#\rho_{{\mathcal I}^{\alpha_1}}\diamond C^{1,t,\rho}_{t}$ and
$\#\rho_{{\mathcal I}^{\alpha_2}}\diamond C^{1,t,\rho}_{t}$.
We know that
$\#\rho_{{\mathcal I}^{\alpha_1}}\diamond C^{1,t,\rho}_{t^{\alpha_2}}$
follows approximately the law of
$\#K_{\log(\alpha_2/\alpha_1)}^\infty$, where $K_\cdot^\infty$ is the Kingman coalescent
started with infinitely many particles.
In particular, $\{\#\rho_{{\mathcal I}^{\alpha_1}}\diamond C^{1,t,\rho}_{t^{\alpha_2}}, \, t\geq t_0\}$
is a tight family of random variables.
Moreover, for any $\varepsilon>0$, due to Proposition~\ref{P:tight},
$\{\#\rho_{{\mathcal I}^{\alpha_2}}\diamond C^{1,t,\rho}_{t^{\alpha_2+\varepsilon}}, \, t\geq t_0\}$
is a tight family as well.

Due to (\ref{agr2}), we have that for each $\varepsilon>0$, the total collection of partition
elements
$\rho_{\CI^{\alpha_2}} C^{1,t,\rho}_{t^{\alpha_2+\varepsilon}}$ has positions in $I_{\alpha_2+\varepsilon}(1,t)$
with overwhelming probability, as $t\to \infty$.
Hence the sparse particle convergence of Proposition \ref{P3.1} applies.
By letting $\varepsilon$ to $0$, and using
(\ref{vli1}) and (\ref{agr3}) as in the proof of Theorem \ref{T:01}, we
obtain the statement (a.i) in the case $m=2$.
The induction step is standard now.

Note that, in view of the proof of part (b), one should verify
the estimates analogous to those of
Proposition~\ref{P3.1a}, as well as the extension
of (a.i), in the slightly more general setting
of the coupled spatial coalescents satisfying (\ref{Erhofin}).

(a.ii)
To prove (\ref{agr7}) note that by the construction in
Subsection~\ref{Sub:CSC}, $\#\rho_{{\mathcal I}^\alpha} \diamond C^{1,t}_t$ has monotone
non-decreasing and c\`adl\`ag
(or c${\rm \grave{a}}$gl${\rm \grave{a}}$d) paths in $\alpha$, for all
$t>0$, almost surely.
Furthermore by Theorem \ref{T:01} we know that the family
$\{\#\rho_{{\mathcal I}^\alpha} \diamond C^{1,t}_t;\,t\ge 0\}$ is tight, for each $\alpha < 1$.
Therefore we obtain (\ref{agr7}).

(b) Again the statements can be easily extended to $\rho=\infty$,
for all $\delta>0$ fixed, using monotonicity and the coupling (\ref{agr4b}).
\end{proof}\sm

%Our task is to turn {\em space} (for the spatial coalescent with
%rebirth) into {\em time} (for the Kingman-type coalescent with rebirth).
%In particular, notice that
%the colors of the particles in both Figures~2 and~3
%match on purpose to illustrate this correspondence of
%space and time.

\begin{proof}[{\bf Proof of Theorem~\ref{T:04}}] So far we have shown
  with Proposition~\ref{P.Cofdd1}  the f.d.d.~convergence. It
remains to show the tightness
in path space.
This is now a direct consequence of the monotonicity of the process
$(N_{\alpha})_{\alpha \in (0,1]}$,
as well as  of all the processes
$(\# \rho_{{\mathcal I}^\alpha}\diamond C^{1,t}_t)_{\alpha \in [\alpha_\ell, \alpha_u]}$ and
$(\#\rho_{{\mathcal I}^\alpha}\diamond C^{1,t,\infty,\delta}_t)_{\alpha \in [\alpha_\ell, \alpha_u]}$
in $\alpha$, more precisely, of the fact that their
paths are non-decreasing and bounded from below (by identity $0$), almost surely
and from above by (\ref{agr7}).
\end{proof}\sm

\subsection{Proof of Theorem~\ref{T:05}}
\label{Sub:tasy}
Fix $\alpha\in(0,1)$.
For $m\in\N$, consider the parameters $\alpha<u_1<u_2<\ldots<u_m < 1$.
Recall the definitions (\ref{grev25a}) and (\ref{Nimmalpha}), and as before denote by
 $\vec{u}/\alpha$ the vector $(u_1/\alpha,\ldots,u_m/\alpha)$.
%In this subsection we give the proof of Theorem~\ref{T:05}.
%, i.e., the
%convergence of $N^{\alpha,t,\rho}_u$ to $N^{\mathrm{mer},\vec{u}/\alpha}_{\frac{u}{\alpha}}$,
%as $t\to\infty$.
%, in the sense of the finite-dimensional distributions.

%In the remainder of the section we will abbreviate
%$N^{\mathrm{mer},\vec{u}/\alpha}_{u_i}$
%by
%$N^{\mathrm{mer}}_{u_i}$, $i=1,\ldots,m$.

\begin{proof}[{\bf Proof of Theorem~\ref{T:05}}]
Note that the case $m=1$ is covered by Theorem \ref{T:01}, hence we will assume $m\geq 2$.
The key is to understand the case $m=2$, since then we can conclude
the argument easily by making the induction step from $m$ to $m+1$.
We will concentrate on (\ref{agr6y}), and we comment at the very end on the
extension (\ref{agr6y2}).

Fix a finite $\rho$ and $t\geq t_0$, where as usual $t_0$ is
taken from Proposition~\ref{P:tight}.
For $i=1,\ldots,m$, define
\be{Chati}
   \bar{C}^i
 :=
   \big\{\pi\in C_{t^{u_i}}^{{\mathrm{birth}}}:\,
   L_{t^{u_i}}^{{\mathrm{birth}}}\in\Lambda^{\alpha,t}\big\}.
\ee
We consider the joint evolution
of partition elements $\bar{C}^i$, $i=1,\ldots,m$.
As mentioned before, there are Poisson($\rho$) many partition elements present
at each site of the $\alpha$-box, at all times $s\geq 0$, almost surely.
In particular, $\#\bar{C}^i$ has Poisson($\rho\cdot\#\Lambda^{\alpha,t}$)
distribution.

Note that, for $t$ large, due to (\ref{agr2}) we will have that, with overwhelming probability,
\be{Esparse}
   \pi^1\not\subseteq\pi^2,\quad\forall\,i<j\in\{1,\ldots,m\},
   \pi^1\in\bar{C}^i,\pi^2\in\bar{C}^j.
\ee
In words, it is highly unlikely to have any equivalence class of
$C_{t^{u_i}}^{{\mathrm{birth}}}\cap \bar{C}^i$ reappear (as a subclass) in the $\alpha$-box
at any of the later times $t^{u_l}$, $l\in\{i+1,\ldots,m\}$.
We will henceforth consider our realization on the event (\ref{Esparse})
in the rest of the argument.

Fix $\delta\in(0,1)$ a small quantity, which will be sent to $0$, eventually.
For each $i\in \{1,\ldots,m\}$ and $s\geq t^{u_i}$, denote by $\bar{N}_{s}^i$
the number of equivalence classes of $C_{s}^{{\mathrm{birth}}}$ containing at least one
element of $\bar{C}^i$.
By Theorem \ref{T:01}, $\bar{N}_{t^{u_2}}^1$ follows approximately
the law of $\#K_{\log{(u_2/\alpha)}}^\infty$.
By (\ref{agr2}), the corresponding
 equivalence classes have locations in $I_{u_2}(1,t)$ at time $t^{u_2}$, and stay in
$I_{u_2}(2+\delta,t)$ during the time interval $[t^{u_2},2 t^{u_2}]$, with overwhelming probability.
Note that, similarly, $\bar{N}_{2t^{u_2}}^1 = \bar{N}_{t^{u_2}}^1$ with overwhelming probability,
as $t\to \infty$.

Next consider during the time interval $[t^{u_2},2t^{u_2}]$ the process counting the
number of equivalence classes process for the
 coalescent
$(C_s^{{\mathrm{birth}}},L_s^{{\mathrm{birth}}})_{s\geq t^{u_2}}$ restricted
to the equivalence classes in $\bar{C}^2$.
Due to Theorem \ref{T:01}, the law of the above counting process
is (after appropriate rescaling) approximately that of $(\#K_s^\infty, s \in [0,\log{(u_2/\alpha)}])$,
as $t\to \infty$.

%Moreover, during the time interval $[t^{u_2},t^{u_2}+t^{u_2-\delta}]$,
%with overwhelming probability, there is no coalescence interaction between the equivalence classes
%of $C^{{\mathrm{birth}}}$
%containing at least one element of $\bar{C}^1$ and those containing at least
%one element of $\bar{C}^2$.
%In other words, for any $\delta>0$,
%the two restricted coalescents evolve separately and
%independently during $[t^{u_2},t^{u_2}+t^{u_2-\delta}]$,  asymptotically as $t\to \infty$.
%For $\delta$ very small it is very likely that the event
%$A_\delta^{1,2}:=\{\bar{N}_{2 t^{u_2}}^2= \bar{N}_{t^{u_2}+t^{u_2-\delta}}^2\}$
%happens.
%Due to the above reasoning, the limit law of
%$(\bar N^1_{t^{u_2}}, \bar N^2_{t^{u_2}})$ is $\CL_1 \otimes \CL_2$.

Also note that, on $A_\delta^{1,2}$, the positions of the
\be{agre41}
 \bar{N}_{2 t^{u_2}}^1+ \bar{N}_{2 t^{u_2}}^2= \bar{N}_{t^{u_2}}^1 + \bar{N}_{t^{u_2}+t^{u_2-\delta}}^2
\ee
equivalence classes in $C_{2 t^{u_2}}^{{\mathrm{birth}}}$, that contain at least one element
 either of $\bar{C}^1$ or of $\bar{C}^2$
are contained in $I_{u_2}(2+\delta,t)$.
Therefore the joint evolution of these equivalence classes during the time interval
$[2t^{u_2},t^{u_3}]$ (by Lemma \ref{L3.1mar} and Section \ref{Sub:convps}) is
again well approximated, on the appropriate scale, by that of the
$(K_{s}^{\bar{N}_{2 t^{u_2}}^1+ \bar{N}_{2 t^{u_2}}^2}, s\in [0,\log{(u_3/u_2)}])$, where the
last coalescent process
depends on $\bar{N}_{2 t^{u_2}}^1,\bar{N}_{2 t^{u_2}}$ solely through its initial configuration.

Denote by $\vec{u}^{\,i}/\alpha$ the vector $(u_1/\alpha,\ldots,u_i/\alpha)\in \R^i$.
It is now clear by the above argument that
$(\bar{N}_{t^{u_3}}^1,\bar{N}_{t^{u_3}}^1+ \bar{N}_{t^{u_3}}^2)$
converges in law as $t\to \infty$ as follows
\be{Ebuild}
\begin{aligned}
   &\big(\bar{N}_{t^{u_3}}^1,\bar{N}_{t^{u_3}}^1+
   \bar{N}_{t^{u_3}}^2\big)
  \\
 &\Tto
   \Big(\#\big\{\pi\in K^{\mathrm{mer},\log(\vec{u}^2/\alpha)}_{\log(u_3/\alpha)}:\,
   [(\kappa(\pi)]_{\mbox{mod}(m)} = 0\big\},
\#\big\{\pi\in K^{\mathrm{mer},\log(\vec{u}^2/\alpha)}_{\log(u_3/\alpha)}:\,
   [(\kappa(\pi)]_{\mbox{mod}(m)} \le 1\big\}\Big).
\end{aligned}
\ee
By setting $u_3=1$, one obtains the result for $m=2$.

Moreover, one can use (\ref{Ebuild}) in the induction step
for the argument where $m \geq 3$.
In fact, using induction one first obtains for each $i$, $3\leq i\leq m$
a generalization of (\ref{Ebuild}):
\be{Ebuildind}
\begin{aligned}
   &\big(\bar{N}_{t^{u_i}}^1,\bar{N}_{t^{u_i}}^1+\bar{N}_{t^{u_i}}^2,\ldots,
   \bar{N}_{t^{u_i}}^1+\ldots+\bar{N}_{t^{u_i}}^{i-1}\big)
  \\
 &\Tto
   \Big(\#\big\{\pi\in K^{\mathrm{mer},\log(\vec{u}^{i-1}/\alpha)}_{
   \log(u_i/\alpha)}:\,[(\kappa(\pi)]_{\mbox{mod}(m)}\le 0\big\},\ldots,
  \\
 &\qquad\qquad\qquad
   \#\big\{\pi\in K^{\mathrm{mer},\log(\vec{u}^{i-1}/\alpha)}_{
   \log(u_i/\alpha)}:\,[(\kappa(\pi)]_{\mbox{mod}(m)} \le i-1\big\}\Big),
\end{aligned}
\ee
and from here easily the general statement of part (a).

Note that part (b) will follow
as usual from (\ref{agr6y}) by monotonicity.
Here it suffices to extend the result of (a)
to the two additional settings where: (i) the initial configuration has precisely one particle at each site, and
(ii) the initial configuration has $1+$Poisson($\rho$) particles at each site, i.i.d.~over sites.
All the reasoning above carries through provided that for each $i=1,\ldots,m$,
the configuration $(C_{t^{u_i}}^{{\mathrm{birth}}},L_{t^{u_i}}^{{\mathrm{birth}}})$
satisfies an analogue of (\ref{Erhofin}).
This property is trivially satisfied in the Poisson case, due to stationarity, as mentioned already.
In the above more general settings one can verify, by approximating the infinite system by the systems
on large finite tori, that the expected number of particles at any particular site at any particular time is bounded from above by a fixed constant ($1$ in the first setting, and $1+\rho$ in the second one).
\end{proof}

\section{Proof of the moment bound}
\label{S:proofmobo}
In this section we present the proof of Proposition~\ref{P:tight}
which
follows the proof of a similar statement for the instantaneous
coalescent stated in the proposition on page 615 in
\cite{BramCoxGri86}. In \cite{BramCoxGri86} the particles move
according to the nearest neighbor random walks,
while here the partition elements
move according to more general random walks. Moreover,
coalescence happens with a rate $\gamma$ delay, and it is
therefore possible (often likely)
to have more than $1$
(up to countable many) partition elements per site.

\begin{proof}[{\bf Proof of Proposition~\ref{P:tight}}]
Recall the box $\Lambda(r)$ from (\ref{e:015}), and let
for $A,B\subseteq\Z^2$,
\be{e:013c}
   \big(C_s^A,L_s^A\big)_{s\geq 0},
\ee
be the coalescent started from the configuration (\ref{Econdcoup})
restricted to locations in $A$.
This coalescent was denoted by
$C^{\CI_A}$ in Subsection~\ref{Sub:CSC}.
If $A=\Lambda(t)$ we will in most cases omit the superscript
from the notation.
For $A,B\subseteq \Z^2$ and $s\ge 0$, let
\be{e:013q}
   \# C_s^A(B)
 :=
   \#\big\{\pi\in C_s^A:\,L_s(\pi)\in B\big\}.
\ee
As done before, if $B=\Z^2$ we simply write $\# C_s^A:=\# C_s^A(\Z^2)$.

Following the lines of Section 3 in \cite{BramCoxGri86}, we
introduce an auxiliary spatial coalescing system
$(\widetilde{C},\widetilde{L})$ which follows the spatial
coalescent dynamics over the time interval $[0,2]$, then
keeps coalescing as long as the number of partition elements is
not decreasing too quickly, while otherwise
the ``coalescence is switched off for a while''.
More precisely, we discretize the time on a
logarithmic scale, i.e., set for $T\ge 0$,
\be{e:016}
   m(T)
 :=
   \left\{\begin{array}{cc}0, & \mbox{if }T\le 1, \\
   2^{\left\lfloor \log_2{T}\right\rfloor}, & \mbox{if }T> 1.
   \end{array}\right.
\ee
In this way we have $T\in[m(T),m(2T)\vee 1]$, $T\geq 0$.

Now, let $(\widetilde{C}_0,\widetilde{L}_0)
:=(C^{\Lambda(t)}_0,{L}^{\Lambda(t)}_0)$,
and run the coalescent until time $T=2$.
To define
$(\widetilde{C}_t,\widetilde{L}_t)$, we proceed by induction.
Put
\be{e:017}
   \tau^{\lfloor\log_2{T}\rfloor}
 :=
   m(2T)\wedge\inf\Big\{s\in [m(T),m(2T)]:\,
   {\bf E}\big[\#\widetilde{C}_s\big]
 \le
   \tfrac{1}{2}{\bf E}\big[\#\widetilde{C}_{m(T)}\big]\Big\},
\ee
and start $\widetilde{C}$ at time
$m(T)$ in the spatial configuration given by $\widetilde{C}_{m(T)}$.
The coalescent $(\widetilde{C},\widetilde{L})$ follows the same
dynamics as the spatial coalescent
on $[m(T),\tau^{\lfloor\log_2{T}\rfloor}]$, while
its partition elements perform independent random walks
with kernel $a(x,y)$ on $[\tau^{\lfloor\log_2{T}\rfloor},m(2T)]$
yielding the random configuration $(\widetilde{C}_{m(2T)},
\widetilde{L}_{m(2T)})$.
Now reset $T:=2T$ and repeat the induction step starting at (\ref{e:017}).
Obviously, ${\bf E}\big[\#\widetilde{C}_t\big]\ge
{\bf E}\big[\#C_t\big]$, for all $t\geq 0$.
In fact, one can easily construct a coupling in such a way that
the corresponding inequality for processes
holds for all times, almost surely. Hence it suffices
to prove Proposition~\ref{P:tight} with $(C,L)$ replaced by $(\widetilde{C},
\widetilde{L})$.

Set
\be{grev47}
  Y_T
 :=
   {\bf E}\big[\#\widetilde{C}_T\big]
 =
   {\bf E}\big[\#\widetilde{C}_T^{\Lambda(t)}(\Z^2)\big],
\ee
and note that $Y_T$ also depends on $t$ through the initial configuration
(\ref{Econdcoup}), although this is suppressed from the notation. \sm

We will need a few preliminary lemmas.
We start with a basic fact estimating the ``speed''
of escape from large balls centered at the origin
for a zero mean random walk with finite exponential moments.

\begin{lemma} Let $({\xi}_t)_{t\ge 0}$ be the
unit rate continuous time random walk on $\Z$ with transition kernel
${b}_t(x,y)$. If $\sum_{x\in\Z}x {b}_1(0,x)=0$ and
$\varphi(\lambda):=\sum_{x\in\Z}e^{\lambda x} {b}_1(0,x)<\infty$, for
all $\lambda>0$, then
\label{L:02} there exists a finite constant $c_0=c_0(\xi)$ such that
\be{e:025}
   {\bf P}\big\{{\xi}_t>u\sqrt{t}\big\}
 \le
   e^{-c_0\,u}
\ee
for all $u,t\ge 1$.
\end{lemma}\sm

\begin{proof}[{\bf Proof}]
The argument is based on standard large deviation techniques.
For all $s,t,\lambda>0$,
\be{e:026}
   {\bf P}\big\{{\xi}_t>s t\big\}
 =
   {\bf P}\big\{e^{\lambda {\xi}_t}>e^{\lambda s t}\big\}
 \le
   e^{-\lambda s t}e^{t(\varphi(\lambda)-1)}.
\ee
In particular, if
$I(s)
 :=
   \sup_{\lambda>0}\big\{s\lambda-(\varphi(\lambda)-1)\big\}
$,
then
\be{e:026x}
   {\bf P}\big\{|\xi_t|>s t\big\}
 \le
   e^{-I(s) t}.
\ee
Note that $I(s):[0,\infty)\to[0,\infty)$ is a convex function, such that
$I(s)=0$ if and only if $s=0$.
Therefore, there exists a positive constant $c_0^1$
such that
\be{e:026b}
   I(s)
 \ge
   c_0^1s,\qquad \mbox{ if }s\ge 1.
\ee

Moreover, under our assumptions on exponential moments,
there exists a finite constant $c_0^2$ (without loss of generality can
assume that $c_0^2\geq 1$) such that
$\varphi(\lambda)\le 1+c_0^2\lambda^2$, for all $\lambda\in[0,1]$.
Thus, for all $s\leq 1$,
\be{e:026c}
   I(s)
 \ge
   \sup_{\lambda\in[0,1]}\big\{s\lambda-c_0^2\lambda^2\big\}
 \ge
   \frac{1}{4c_0^2}s^2,
\ee
where we have used the fact that if $s\leq 1$ then
$\lambda^\ast:= \frac{s}{2c_0^2} \leq 1$.

Now set $c_0:=\min\{c_0^1,(4\tilde{c}_0^2)^{-1}\}$, and take $u,t\geq 1$.
If $u\geq \sqrt{t}$ we obtain (\ref{e:025}) from (\ref{e:026x})
by substituting $s=u/\sqrt{t}$ into
(\ref{e:026b}).
Similarly, if $1\leq u\leq \sqrt{t}$ we obtain (\ref{e:025})
by substituting $s=u/\sqrt{t}$ into (\ref{e:026c}).
\end{proof}\sm

The next result states that
if the spatial coalescent starts in $\Lambda(t)$,
then at time $T$ the fraction of partition elements
which lie outside of $\Lambda(t+u\sqrt{T})$ decreases
at least exponentially fast,
as $u\to\infty$.

\begin{lemma}
Fix $t>0$.
Let $\bar{R}:=(\bar{R}^1,\bar{R}^2)$ be the random walk
on $\Z^2$ with kernel $a(x,y)$.
Fix $c_0=c_0 (\bar R)$ such that (\ref{e:025}) holds.
Put $c_1:=2\cdot (2^{\tfrac{5}{(\sqrt{2}-1)}}\wedge e^{c_2})$ where
$c_2=c_2(\bar{R}):=\sqrt{\tfrac{2}{7}}
(c_0(\bar{R}^1)\wedge c_0(\bar{R}^2))$.
Then \label{L:03}
\be{e:018}
   {\bf E}\Big[\#\widetilde{C}_T\big(\Lambda^c(t+u\sqrt{T})\big)\Big]
 \le
   c_1\,e^{-c_2u}Y_T,
\ee
for all $u\ge 0$ and $T\ge 1$.
\end{lemma}\sm

Choosing $a$ large  enough so that $c_1 e^{-c_2 a}\le 1/3$ we obtain
the following:

\begin{corollary} For sufficiently large $a\geq 1$,
\be{e:033}
   {\bf E}\Big[\#\widetilde{C}_T\big(\Lambda^c(t+a\sqrt{T})\big)\Big]
 \le
   \tfrac{1}{3} Y_T,
\ee
for all $T\ge 1$. \label{C:04}
\end{corollary}\sm

\begin{proof}[{\bf Proof of Lemma~\ref{L:03}}]
The proof is by induction over
$\left\lfloor\log_2{T}\right\rfloor$.
First, suppose that $2\le T\le 2^4$ and $u\ge 1$.
By comparison with the independent random walks equal in law to
$\bar{R}:=(\bar{R}^1,\bar{R}^2)$ on $\Z^2$, we obtain
(with $\| \cdot \|$ the maximum norm)
\be{e:029}
\begin{aligned}
   {\bf E}\Big[\#\widetilde{C}_T\big(\Lambda^c(t+u\sqrt{T})\big)\Big]
 &\le
   {\bf E}\big[\# C_0^{\Lambda(t)}\big]\,
   {\bf P}^{(0,0)}\big\{\|\bar{R}_T\|\ge u\sqrt{T}\big\}
  \\
 &\le
   {\bf E}\big[\# C_0^{\Lambda(t)}\big]
   \Big( {\bf P}^0\big\{|\bar{R}^1_T| \ge u\sqrt{T}\big\}
   +{\bf P}^0\big\{|\bar{R}^2_T| \ge u\sqrt{T}\big\}\Big)
  \\
 &\le
   4\cdot {\bf E}\big[\# C_0^{\Lambda(t)}\big]\,
   e^{-(c_0(\bar{R}^1)\wedge c_0(\bar{R}^2)) u}.
\end{aligned}
\ee

By definition, $Y_T\ge Y_{2^4}\ge\frac{1}{2}Y_{2^3}\ge \ldots\ge 2^{-4}
\# C_0$. Moreover the map
$s \mapsto \mathbf{E}\big[\# C_s\big]$ is continuous, and
therefore
\be{e:030}
   {\bf E}\Big[\#\widetilde{C}_T\big(\Lambda^c(t+u\sqrt{T})\big)\big]
 \le
   2^6 \cdot e^{-(c_0(\bar{R}^1)\wedge c_0(\bar{R}^2)) u}\cdot Y_T,
%\end{aligned}
\ee
as required. So (\ref{e:018}) holds in the case $2\le T\le 2^4$, for all
$u\ge 1$, and for $u\in [0,1]$, (\ref{e:018}) holds trivially due to the
fact that $c_1e^{-c_2}\geq 1$.

Suppose now that for some $m\ge 1$, (\ref{e:018}) holds for all $2\le
T\le 2^{m+3}$. Then for $T\in(2^{m+3},2^{m+4}]$,
\be{e:031}
\begin{aligned}
   &{\bf E}\Big[\#
   \widetilde{C}_T\big(\Lambda^c(t+u\sqrt{T})\big)\Big]
  \\
 &\le
   {\bf E}\Big[
   \# \widetilde{C}_{2^m}\big(\Lambda^c(t+\frac{u}{2}\sqrt{T})\big)\Big]
   +Y_{2^m}{\bf P}^{(0,0)}\big\{\|\bar{R}_{T-2^m}\|\ge
   \frac{u}{2}
   \sqrt{T}\big\}
  \\
 &\le
   {\bf E}\Big[\#
   \widetilde{C}_{2^m}\big(\Lambda^c(t+(\sqrt{2}u)2^{\frac{m}{2}})\big)\Big]
   +Y_{2^m}{\bf P}^{(0,0)}\big\{\|\bar{R}_{T-2^m}\|\ge
   \frac{u}{2}
   \sqrt{\tfrac{T}{(T-2^m)}}\sqrt{(T-2^m)}\big\}.
\end{aligned}
\ee
The first inequality above is obtained by the following observation:
each partition element in $\Lambda^c(t+u\sqrt{T})$ at time $T$
corresponds to some partition element,
located either in $\Lambda(t+\frac{u}{2}\sqrt{T})$ or its complement, at time $2^m$.
Applying the induction hypotheses to the first term,
and Lemma~\ref{L:02} to the second term on the right hand
side of (\ref{e:031}), we obtain that
\be{e:032}
\begin{aligned}
   {\bf E}\Big[\#
   \widetilde{C}_T\big(\Lambda^c(t+u\sqrt{T})\big)\Big]
 &\le
   Y_{2^m}\big(c_1 e^{-c_2 \sqrt{2}u}+
   2e^{-c_2 u}\big)
  \\
 &\le
   Y_{T}c_1e^{-c_2 u}\big(2^4e^{-c_2(\sqrt{2}-1)u}+\frac{2^6}{c_1}\big),
\end{aligned}
\ee
where we have used the facts that $\sqrt{T/(T-2^m)}\ge \sqrt{8/7}$,
for all
  $T\in(2^{m+3},2^{m+4}]$, and $Y_T\ge 2^{-4}Y_{2^m}$.

Define $u_0:=\frac{5\ln{2}}{c_2(\sqrt{2}-1)}$. Then an elementary
calculation shows that for all $u\ge u_0$,
\be{e:033a}
   2^4e^{-c_2(\sqrt{2}-1)u}+\frac{2^6}{c_1}
 \le
   2^4e^{-c_2(\sqrt{2}-1)u_0}+\frac{1}{2}
 \le
   2^4 2^{-5}+\frac{1}{2}
 \le
   1,
\ee
while for all $u\in[0,u_0]$,
$c_1e^{-c_2u}\ge c_12^{-\frac{5}{\sqrt{2}-1}}\geq 1$, so
(\ref{e:018}) trivially holds for all $u\in[0,u_0]$.
This completes the induction step and the proof.
\end{proof}\sm

We next provide an estimate of the rate of decrease
for the number of partition elements
during an interval of time, provided that the coalescence dynamics is
switched on.

For two partition elements $\{i\},\{j\}\in C_0$, put
\be{sigma}
   \sigma^{\{i,j\}}
 :=
   \min\big\{u\geq 0:\,L_u(\{i\})=L_u(\{j\})\big\}
\ee
as the waiting time until these particles
share the same location, and set
\be{e:hA}
   h_s^\gamma(A)
 :=
   \inf_{i,j\in{\mathcal I}^A}{\bf P}\big\{\sigma^{\{i,j\}}\le s\big\}.
\ee
One can verify using a last-exit-time decomposition and
the assumption~(\ref{vl3})
(compare Lemma~5 in~\cite{BraGri80})
that for fixed $b>0$,
\be{e:052help}
   h_{r^2}\big(\Lambda(b r)\big)
 \ge
   M(b) \frac{1}{\log(r)},
\ee
for some $M(b)>0$.

Similarly, define
$   \tau^{\{i,j\}}
 :=
   \min\big\{u \geq 0:\, i\sim^u j\big\},
$
and set for $A\subseteq\Z^2$,
\be{e:HA}
   H_s^\gamma(A)
 :=
   \inf_{i,j\in{\mathcal I}^A}{\bf P}\big\{\tau^{\{i,j\}}\le s\big\}.
\ee
We are particularly interested in bounding from below the quantity
\be{grev70}
   H_{4R_T^2}^\gamma(\Lambda(\sqrt{2}R_T)),
\ee
where
\be{e:R}
   R_T=R^{a,t}_T
 :=
   7(1+a)\sqrt{\tfrac{t^2+a T}{Y_T}},
\ee
with $a\geq 1$ chosen according to Corollary~\ref{C:04} such that
(\ref{e:033}) holds.
We will henceforth assume that $T\leq t^3$ (as in (\ref{e:046}) below).
Then, if
\be{EsT}
   s_T
 :=
   4 R_T^2,
\ee
inequality (\ref{e:052help}) implies that
\be{e:052a}
   h_{s_T/2}\big(\Lambda(\sqrt{2}R_T)\big)
 \ge
\frac{M(1)}{\log{R_T}}
\ge
   \frac{M_1}{\log{t}},
\ee
where $M_1\in(0,2 M(1)/3)\subset(0,\infty)$ is chosen depending on $a$.
Recalling inequality (7.48) from
\cite{GreLimWin05}, we obtain that
\be{e:053}
   H_{s_T}^\gamma\big(\Lambda(\sqrt{2}R_T)\big)
 \ge
   \frac{\gamma}{2+\gamma}\Big(1-\exp\big(-\frac{2+\gamma}{2}s_T\big)\Big)
   h_{s_T/2}\big(\Lambda(\sqrt{2}R_T)\big)
 \ge
   \frac{M_2}{\log{t}},
\ee
for some $M_2\in(0,\infty)$, for all $t\ge 2$, where we use
$ s_T
 \ge
   4\cdot 49\cdot(1+a)^2\cdot \frac{t^2}{Y_0}>0$,
since $t\geq 2$.

\begin{lemma}[Rate of decay for the auxiliary coalescent]
  \label{L:04}
Let $2\leq T\le r<r+s\le 2T$.
Suppose that $Y_T\ge 49$, and that $\widetilde{C}$ is
  coalescing during the entire time interval $[T,r+s)$.
Then
\be{e:034}
   Y_{r+s}
 \le
   Y_r\exp\Big[-\frac{1}{3}H_s^\gamma\big(\Lambda(\sqrt{2}R_T)\big)\Big].
\ee
\end{lemma}\sm

\begin{proof}[{\bf Proof}]
Write $C_s^{C}$ for the spatial coalescent
started in the random partition $C$ at time $0$,
and evaluated at time $s$.
For all $T\le
  r<r+s\le 2T$,
\be{e:040}
   Y_{r+s}
 \le
   {\bf E}\big[\# C_s^{\widetilde{C}_r(\Lambda(t+a\sqrt{r}))}\big]
 +
   {\bf E}\big[\# \widetilde{C}_r(\Lambda^c(t+a\sqrt{r}))\big].
\ee

Choose a covering of $\Lambda(t+a\sqrt{r})$ by
\be{e:nT}
%\begin{aligned}
   n_T
 :=
   \left\lfloor 1+[\mbox{Area}(\Lambda(t+a\sqrt{T}))]^{1/2}/R_T\right\rfloor^2
\ee
disjoint  boxes $\{\Lambda_{i,r},\,i=1,\ldots,n_T\}$
of side length
\be{e:lT}
   l_T
 :=
   \left(\frac{\mbox{Area}(\Lambda(t+a\sqrt{r}))}{n_T}\right)^{1/2}
  \le
   \sqrt{2}R_T.
\ee
The last inequality holds since $r \in [T,2T]$.

After ignoring coalescing events between partition elements that are
located in different sub-boxes $\Lambda_{r,i}\cap
\Lambda_{r,j}=\emptyset$
at time $r$,
one can bound from above the first term on the right hand side of
(\ref{e:040}) by
\be{e:041}
   \sum_{i=1}^{n_T}\sum_{C:C(\Lambda_{i,r}^c)=\emptyset}{\bf
     P}\big\{\widetilde{C}_r(\Lambda(t+a\sqrt{r})\cap\Lambda_{i,r})=
         C(\Lambda(t+a\sqrt{r})\cap\Lambda_{i,r})\big\}
   {\bf E}\big[\#C_s^C\big].
\ee
It is straightforward to conclude, as in (7.44)--(7.46)
in~\cite{GreLimWin05}, that for $C$ as above
\be{e:042}
   {\bf E}\big[\# C_s^C\big]
 \le
   \# C-\big(\# C-1\big)H_s^\gamma\big(\Lambda(\sqrt{2}R_T)\big).
\ee

Insert (\ref{e:042}) into (\ref{e:041}) to get
\be{e:043}
\begin{aligned}
    &{\bf E}\big[\#C_s^{\widetilde{C}_r(\Lambda(t+a\sqrt{r}))}\big]
  \\
 &\le
   \sum_{i=1}^{n_T}\sum_{C:C(\Lambda_{i,r}^c)=\emptyset}{\bf
     P}\big\{\widetilde{C}_r(\Lambda(t+a\sqrt{r})\cap\Lambda_{i,r})=
         C(\Lambda(t+a\sqrt{r})\cap\Lambda_{i,r})\big\}\cdot
 \\
 &\qquad
   \cdot\Big(\# C(\Lambda(t+a\sqrt{r})\cap\Lambda_{i,r})-
   \big(\# C(\Lambda(t+a\sqrt{r})\cap\Lambda_{i,r})-1\big)
   H_s^\gamma\big(\Lambda(\sqrt{2}R_T)\big)\Big)
  \\
 &=
   {\bf E}\Big[\#\widetilde{C}_r(\Lambda(t+a\sqrt{r}))\big]    -
   \Big({\bf E}\big[\#\widetilde{C}_r(\Lambda(t+a\sqrt{r}))\big]-n_T\Big)
   H_s^\gamma\big(\Lambda(\sqrt{2}R_T)\big)
  \\
 &\le
   {\bf E}\big[\#\widetilde{C}_r(\Lambda(t+a\sqrt{r}))\big]
   \Big(1-\frac{1}{2}H_s^\gamma\big(\Lambda(\sqrt{2}R_T)\big)\Big).
\end{aligned}
\ee
For the last inequality in (\ref{e:043})
we use (\ref{e:nT}) and the following observations
\begin{itemize}
\item[{(a)}] $Y_u \geq Y_T/2$, for all $u\in [T, r+s)$, and therefore
in particular, $Y_r \geq Y_T/2$, since otherwise
the coalescing would not last during the entire interval $[T,r+s)$,
\item[{(b)}]  for any $r\geq 1$,
\be{EYrequal}
\begin{aligned}
   Y_r
 &=
   {\bf E}\big[\#\widetilde{C}_r(\Lambda(t+a\sqrt{r}))\big]
   +{\bf E}\big[\#\widetilde{C}_r(\Lambda^c(t+a\sqrt{r}))\big]
  \\
 &\leq
   {\bf E}\big[\#\widetilde{C}_r(\Lambda(t+a\sqrt{r}))\big]+
 \frac{Y_r}{3},
\end{aligned}
\ee
by Corollary~\ref{C:04}, and
\be{e:0nt}
\begin{aligned}
   n_T
 &\le
   \big(\frac{2}{7}\sqrt{Y_T}\big)^2
  \leq \frac{4}{49}\cdot 4 Y_r
  \\
 &\le
   \frac{4\cdot 4}{49}\cdot \frac{3}{2}
   {\bf E}\big[\#\widetilde{C}_r(\Lambda(t+a\sqrt{r}))\big]
 <
   \frac{1}{2}{\bf E}\big[\# \widetilde{C}_r(\Lambda(t+a\sqrt{r}))\big].
\end{aligned}
\ee
\end{itemize}
Now by (\ref{e:040}), (\ref{e:043}), (\ref{EYrequal})
and (\ref{e:033}),
we have
\be{e:044}
\begin{aligned}
   Y_{r+s}
 &\le
   {\bf E}\big[\#\widetilde{C}_r(\Lambda(t+a\sqrt{r}))\big]
   \Big(1-\frac{1}{2}H_s^\gamma\big(\Lambda(\sqrt{2}R_T)\big)\Big)
 +
   {\bf E}\big[\#\widetilde{C}_r(\Lambda^c(t+a\sqrt{r}))\big]
  \\
 &=
   Y_r\Big(1-\frac{1}{2}H_s^\gamma\big(\Lambda(\sqrt{2}R_T)\big)\Big)
   +\frac{1}{2}H_s^\gamma\big(\Lambda(\sqrt{2}R_T)\big)
   {\bf E}\big[\#\widetilde{C}_r(\Lambda^c(t+a\sqrt{r}))\big]
  \\
 &\le
   Y_r\Big(1-\frac{1}{2}H_s^\gamma\big(\Lambda(\sqrt{2}R_T)\big)\Big)
   +\frac{1}{6}Y_rH_s^\gamma\big(\Lambda(\sqrt{2}R_T)\big)
  \\
 &=
   Y_r\Big(1-\frac{1}{3}H_s^\gamma\big(\Lambda(\sqrt{2}R_T)\big)\Big)
  \\
 &\le
   Y_r\exp\Big[-\frac{1}{3}H_s^\gamma\big(\Lambda(\sqrt{2}R_T)\big)\Big],
\end{aligned}
\ee
as required.
\end{proof}\sm

\begin{lemma}[Upper bound for the decay rate of partition elements]
Fix $t\ge 2$, and let for $T\ge 2$,
\be{e:044x}
   g(T)
 :=
   \frac{\log{(1+\frac{T}{t^2})}}{\log{t}}\cdot
   Y_T\cdot\Big(1\vee\frac{{\bf E}\big[\#
   C_2^{\Lambda(t)}\big]}{t^2}\Big)^{-1},\qquad T\ge 2.
\ee
Then there exists a finite constant $M$ such that \label{L:05}
\be{e:045}
   g(T)\le M,\qquad 2\le T\le 4,
\ee
and
\be{e:046}
   g(2T)
 \le
   g(T)\vee M,
   \qquad 2\le T\le t^3.
\ee
\end{lemma}\sm

\begin{proof}[{\bf Proof}]
Recall $M_2$ from (\ref{e:053}), and fix $a\ge
1\vee\frac{M_2\log_2{5}}{48}$
suitably large
such that (\ref{e:033}) holds. Put
\be{M}
   M
 :=
   \frac{3\cdot 16\cdot 49\cdot a(1+a)^2}{M_2},
\ee
and notice that $M\ge 49\cdot \log_2{5}$.

Assume first that $2\le T\le 4$. In this case,
since $Y_T/t^2\leq Y_2/t^2\leq 1\vee {\bf E}[C_2^{\Lambda(t)}]/t^2$,
and since $\log(1+x) \leq x$, for all $x>-1$,
\be{e:047}
   g(T)
 \le
   \frac{Tt^2}{t^2\ln{t}}
 \le
   \frac{4}{\ln{2}}
 \le
   M.
\ee
Next assume that $2\le T\le t^3$ and $Y_T\le 49$. Then since
$Y_{2T}\le Y_T\le 49$,
we get
\be{e:048}
   g(2T)
 \le
   49\frac{\ln{(1+2t)}}{\ln{t}}
 \le
   49\log_2{5}\le M.
\ee
It therefore remains to prove (\ref{e:046}) for $Y_T>49$.
%W.l.o.g.\ we may assume that $T\in[2,t^3]$, is a power of $2$.
Without loss of generality we may assume that
\be{e:99}
   \mbox{$\widetilde{C}$ is coalescing during
the entire interval  $(T,\frac{3}{2}T)$}.
\ee
Indeed otherwise we could find an $m\in\N$
such that $\tau^m\in(T,\frac{3}{2}T)$ (recall \ref{e:017}) and therefore since
$Y_{2T}\leq Y_{\tau^m}\leq Y_T/2$, we get
$\frac{g(2T)}{g(T)}\leq\frac{1}{2}\cdot\frac{\log{(1+\frac{2T}{t^2})}}
{\log{(1+\frac{T}{t^2})}}\leq 1$.

However, under (\ref{e:99}), Lemma~\ref{L:04} applies
with any $(r,r+s]\subset(T,\frac{3}{2}T]$,
so that $\lfloor\frac{T}{2s_T}\rfloor$ iterations of (\ref{e:034})
yield that
\be{e:049}
   Y_{2T}
 \le
   Y_T\exp\big[-\frac{1}{3}\left\lfloor\frac{T}{2s_T}\right\rfloor
   H_{s_T}^\gamma\big(\Lambda(\sqrt{2}R_T)\big)\big].
\ee
By (\ref{EsT}),
\be{e:050}
\begin{aligned}
   \left\lfloor\frac{T}{2s_T}\right\rfloor
 \ge
    \frac{T}{4s_T}
 &=
   \frac{Y_T T}{16\cdot 49\cdot (1+a)^2(t^2+a T)}
  \\
 &\ge
   g(T)\frac{T\ln{t}}{16\cdot 49\cdot (1+a)^2(t^2+a
 T)\ln{(1+\frac{T}{t^2})}}.
%  \\
% &\ge
%   M_0 g(T)\frac{T\ln{t}}{(t^2+aT)\ln{(1+\frac{T}{t^2})}},
\end{aligned}
\ee
%for some constant $M_0=M_0(a)\in(0,\infty)$.
Finally,
inserting (\ref{e:050}) and (\ref{e:049}) into (\ref{e:044x}),
and recalling (\ref{e:053}), yields
\be{e:054}
\begin{aligned}
   \frac{g(2T)}{g(T)}
 &\le
   \frac{\log{(1+\frac{2T}{t^2})}}{\log{(1+\frac{T}{t^2})}}
   \exp\Big[-\frac{1}{3}\left\lfloor\frac{T}{2s_T}\right\rfloor
   H_{s_T}^\gamma\big(\Lambda(\sqrt{2}R_T)\big)\Big]
  \\
 &\le
   \exp\Big[\frac{T}{(t^2+T)\log{(1+\frac{T}{t^2})}}-
\frac{1}{3}\left\lfloor\frac{T}{2s_T}   \right\rfloor
   H_{s_T}^\gamma\big(\Lambda(\sqrt{2}R_T)\big)\Big]
  \\
 &\le
   \exp\Big[\frac{T}{(t^2+T)\log{(1+\frac{T}{t^2})}}\big(1-M^{-1}g(T)
   \big)\Big].
\end{aligned}
\ee
We therefore find that either $g(T)\le M$ or if $g(T)>M$ then
$g(2T)\le g(T)$, which proves (\ref{e:046}).
\end{proof}\sm

To finish off the proof of the proposition, note that
Lemma \ref{L:05} readily implies
$g(T)\le M$, for all $t\ge 2$, $0\le T\le t^3$.
Therefore,
\be{e:055}
   Y_T
 \le
   M\frac{\log{t}}{\log(1+\frac{T}{t^2})}
   \Big(1\vee\frac{{\bf E}\big[\#
   C_2^{\Lambda(t)}\big]}{t}\Big),\qquad 2\leq T \leq t^3,
\ee
and after replacing $t$ with $t^{\alpha/2}$ and $T$ with $t^{\beta}$ where
$\beta\in (\alpha, 3\alpha/2]$,
\be{e:055x}
   Y_{t^{\beta}}
 \le
M\frac{\log{t^{\alpha/2}}}{\log(1+\frac{t^\beta}{t^\alpha})}
   \Big(1\vee\frac{{\bf E}\big[\#C_2^{\Lambda(t^{\alpha/2})}\big]}{t^\alpha}\Big)
 \le
   M\Big(1\vee\frac{\alpha}{2(\beta-\alpha)}\vee
   \frac{{\bf E}\big[\# C_2^{\Lambda(t^{\alpha/2})}\big]}{t^\alpha} \Big).
\ee
\end{proof}\sm

{\em Acknowledgment.}
The PStricks coding of figures in this paper was made by Laura Derksen,
as part of her  NSERC USRA training during summer 2005.

\end{document}